\documentclass[10pt]{article}
\usepackage[francais]{babel}
\usepackage[T1]{fontenc}
\usepackage[utf8]{inputenc}

\usepackage{amsmath}
\usepackage{amsfonts}
\usepackage{amssymb}
\usepackage{xypic}
\xyoption{curve}
\usepackage{mathrsfs}
\usepackage{hyperref}
\usepackage[a4paper, margin=3cm]{geometry}

\newcommand{\guil}[1]{\og\ignorespaces #1\unskip\fg}
\newcommand{\cone}{\text{cône}}
\newcommand{\cocone}{\text{cocône}}
\newcommand{\et}{{\text{\'et}}}
\newcommand{\red}{{\text{r\'ed}}}
\newcommand{\Nis}{{\text{Nis}}}
\newcommand{\vers}[1]{\overset{#1}{\longrightarrow}}
\newcommand{\depuis}{\longleftarrow}
\newcommand{\isomto}{\overset{\sim}{\rightarrow}}
\newcommand{\Ab}{\mathrm{Ab}}
\newcommand{\Ens}{\mathrm{Ens}}
\newcommand{\Enspt}{\Ens_\bullet}
\newcommand{\univ}{\mathrm{univ}}
\newcommand{\opp}{\mathrm{opp}}
\newcommand{\PST}[1][S,\Lambda]{\mathbf{PST}_{#1}}
\newcommand{\ST}[1][S,\Lambda]{\mathbf{ST}_{#1}}
\newcommand{\Hyper}{\mathbf{H}}
\newcommand{\HZl}{\mathbf{H}_{\mathbf{Z}/\ell}}
\newcommand{\tw}{\mathrm{tw}}
\newcommand{\Fl}{\mathbf{F}_\ell}
\DeclareMathOperator{\Inv}{Inv}
\DeclareMathOperator{\Gr}{Gr}
\DeclareMathOperator{\Tr}{Tr}
\DeclareMathOperator{\codim}{codim}
\DeclareMathOperator{\degtr}{deg.tr.}
\DeclareMathOperator{\Ex}{Ex}
\DeclareMathOperator{\Sing}{Sing}
\DeclareMathOperator{\coker}{coker}
\DeclareMathOperator{\image}{Im}
\DeclareMathOperator{\Der}{D}
\newcommand{\DerNeg}{\Der_{\geq 0}}
\newcommand{\Dbcoh}{\Der^\bornee_\coh}
\DeclareMathOperator{\Comp}{C}
\newcommand{\CompNeg}{\Comp_{\geq 0}}
\DeclareMathOperator{\Cor}{Cor}
\DeclareMathOperator{\Tor}{Tor}
\DeclareMathOperator{\oub}{oub}
\DeclareMathOperator{\Id}{Id}
\DeclareMathOperator{\hocolim}{hocolim}
\DeclareMathOperator{\colim}{colim}
\DeclareMathOperator{\Sq}{Sq}
\DeclareMathOperator{\Pic}{Pic}
\DeclareMathOperator{\GL}{GL}
\DeclareMathOperator{\Spec}{Spec}
\DeclareMathOperator{\SheafHom}{\mathbf{Hom}}
\DeclareMathOperator{\Hom}{Hom}
\DeclareMathOperator{\End}{End}
\DeclareMathOperator{\Th}{Th}
\newcommand{\Triv}[1][\mathcal X]{\mathrm{Triv}/#1}
\newcommand{\tC}{\widetilde{C}}
\newcommand{\tL}{\widetilde{L}}
\newcommand{\tMnn}{\widetilde{M}^{\mathrm{nn}}}
\newcommand{\tM}{\widetilde{M}}
\newcommand{\tH}{\widetilde{H}}
\newcommand{\DMNeg}[1][S]{DM^{\text{eff}}_{\leq 0}(#1)}
\newcommand{\DMeffgm}[1][k]{DM^{\text{eff}}_{\gm}(#1)}
\newcommand{\DMmoins}[1][S]{DM^{\text{eff}}_-(#1)}
\newcommand{\SH}[1][k]{\mathcal{SH}(#1)}
\newcommand{\SHnaive}[1][k]{\mathcal{SH}_{\text{na\"ive}}(#1)}
\newcommand{\Hosimppt}[1][S]{\mathcal{H}_{\text{s},\bullet}(#1)}
\newcommand{\Hosimp}[1][S]{\mathcal{H}_{\text{s}}(#1)}
\newcommand{\Hopt}[1][S]{\mathcal{H}_\bullet(#1)}
\newcommand{\Ho}[1][S]{\mathcal{H}(#1)}
\newcommand{\Sc}[1][S]{\mathrm{Sc}_\bullet/#1}
\newcommand{\Sm}[1][S]{Sm/#1}
\newcommand{\Smqp}[1][S]{Sm^{\text{qp}}/#1}
\newcommand{\SmCor}[1][S]{SmCor/#1}
\newcommand{\OO}{{\mathscr{O}}}
\newcommand{\cequi}{c_{\text{equi}}}
\newcommand{\tZtr}{\widetilde{\mathbf{Z}_{\text{tr}}}}
\newcommand{\Ztr}{\mathbf{Z}_{\text{tr}}}
\newcommand{\Ltr}{\Lambda_{\text{tr}}}
\newcommand{\Sym}[1]{\mathfrak{S}_{#1}}
\newcommand{\AlgSym}{\mathbf{S}}
\renewcommand{\L}{\mathbf{L}}
\newcommand{\R}{\mathbf{R}}
\newcommand{\bornee}{\mathrm{b}}
\newcommand{\coh}{\mathrm{coh}}
\newcommand{\gm}{\text{gm}}
\newcommand{\Bgm}{\mathbf{B}_\gm}
\newcommand{\Egm}{\mathbf{E}_\gm}
\newcommand{\Bet}{\mathbf{B}_\et}
\newcommand{\Gm}{\mathbf{G}_\text{m}}
\newcommand{\abs}[1]{\left|#1\right|}
\newcommand{\floor}[1]{\left\lfloor#1\right\rfloor}
\newcommand{\acc}[2]{\left<#1,#2\right>}
\newcommand{\acccro}[2]{\left[#1,#2\right]}

\newcommand{\cercle}[1]{\textcircled{\raisebox{-0.8pt}{#1}}}

\title{Opérations de Steenrod motiviques}
\author{Joël Riou}
\date{13 juillet 2012}

\newtheorem{theoreme}[subsubsection]{Théorème}
\newtheorem{definition}[subsubsection]{Définition}
\newtheorem{lemme}[subsubsection]{Lemme}
\newtheorem{proposition}[subsubsection]{Proposition}
\newtheorem{remarque}[subsubsection]{Remarque}

\newtheorem{corollaire}[subsubsection]{Corollaire}

\newtheorem{theoreme1}[subsection]{Théorème}

\newtheorem{proposition1}[subsection]{Proposition}
\newtheorem{remarque1}[subsection]{Remarque}

\newtheorem{corollaire1}[subsection]{Corollaire}

\begin{document}

\maketitle

\begin{abstract}
Cet article apporte des corrections à la construction par Voevodsky des
opérations de Steenrod agissant sur la cohomologie motivique à coefficients
$\mathbf{Z}/\ell\mathbf{Z}$ des espaces de la catégorie homotopique de
Morel-Voevodsky d'un corps parfait de caractéristique différente de $\ell$.
En outre, comme conséquence de la méthode de démonstration d'un théorème de
Voevodsky portant sur les opérations cohomologiques stables, on montre que
sur un corps de caractéristique zéro, le spectre représentant
la cohomologie motivique à coefficients $\mathbf{Z}/\ell\mathbf{Z}$ n'admet
pas d'endomorphisme \guil{stablement fantôme} non nul.
\end{abstract}

\begin{abstract}
This article fills some gaps in Voevodsky's construction of the Steenrod
operations acting on the motivic cohomology with coefficients in
$\mathbf{Z}/\ell\mathbf{Z}$ of motivic spaces in the sense of Morel and
Voevodsky over a perfect field of characteristic different from $\ell$.
Moreover, as a consequence of the method of proof of a theorem by Voevodsky
on stable cohomology operations, we show that the spectrum that represents
motivic cohomology with coefficients $\mathbf{Z}/\ell\mathbf{Z}$ has no
nonzero ``superphantom'' endomorphism.
\end{abstract}

\tableofcontents

Cet article se veut avant tout un texte donnant une construction solide des
opérations de Steenrod agissant sur la cohomologie motivique modulo $\ell$
des variétés lisses ou plus généralement des objets de la catégorie
homotopique $\Ho[k]$ de Morel-Voevodsky associée à un corps parfait $k$ de
caractéristique différente de $\ell$.

Les énoncés seront essentiellement ceux de \cite{voevodsky-reduced}. Si les
principes généraux de \cite{voevodsky-reduced} sont évidemment corrects, un
certain nombre de détails en sont imparfaitement rédigés et surtout il
semble qu'à quelques endroits il faille utiliser des arguments
substantiellement différents. Vu l'importance des résultats obtenus en
utilisant les opérations de Steenrod motiviques (conjecture de Milnor
\cite{voevodsky-milnor}, conjecture de Bloch-Kato
\cite{voevodsky-bloch-kato}), il paraissait important d'y remédier. Il a
semblé préférable de rédiger le présent article comme un ensemble cohérent
plutôt que comme une liste d'\emph{errata} à \cite{voevodsky-reduced}. Cet
article couvre l'essentiel du sujet traité dans \cite{voevodsky-reduced}, à
deux exceptions près pour lesquelles nous ferons référence à
\cite{voevodsky-reduced} : l'annulation de $P^i$ pour $i<0$
\cite[Proposition~3.6]{voevodsky-reduced}
(cf.~\S\ref{subsection-construction-operation-p-i-b-i})
et la formule $u^2=\tau v+\rho u$ dans la cohomologie
modulo $2$ du classifiant de $\mu_2$
\cite[Theorem~6.10]{voevodsky-reduced} (cf.
proposition~\ref{proposition-formule-u-carre}).
Pour le reste, le présent article
pourra se lire indépendamment de \cite{voevodsky-reduced}.

Le point-clef de la construction des opérations de Steenrod motiviques
réside dans la définition d'une transformation naturelle pour tous $r\geq
0$ et $\mathcal X\in\Ho[k]$ :
\[H^{2r,r}(\mathcal X)\to H^{2rn,rn}(\mathcal X\times (G\backslash
U))\;\text{,}\]
une telle opération étant associée à un $k$-schéma lisse $U$ sur lequel
agit librement un sous-groupe $G$ du groupe symétrique $\Sym n$.
La principale idée nouvelle développée dans cet article est de ne pas
considérer que la donnée est celle de $U\in\Sm[k]$ muni d'une action libre
de $G$, mais de considérer un schéma $S\in \Sm[k]$ muni d'un $G$-torseur
étale $U$ (ainsi $G\backslash U=S$). Toutes les constructions seront alors
relatives au schéma de base $S$ et en particulier la transformation
naturelle cherchée peut se réinterpréter comme un morphisme $K_r\to
K_{rn}$ dans la catégorie homotopique pointée $\Hopt[S]$ au-dessus de $S$
où les espaces $K_i$ sont les espaces d'Eilenberg-Mac Lane motiviques
représentant la cohomologie motivique en bidegré $(2i,i)$.

À vrai dire, en tant que préfaisceau d'ensembles pointés, l'espace
d'arrivée ne sera pas exactement $K_{rn}$, mais une version tordue de
celui-ci. Plus précisément, à chaque fibré vectoriel $E$ sur $S$ est
associé un préfaisceau $K_E$ sur $\Sm[S]$ dont le type d'homotopie ne
dépend que du rang de $E$. C'est un isomorphisme de Thom relatif qui fait
l'objet de la section~\ref{section-thom-relatif}. Par ailleurs, à une
action du groupe fini $G$ sur un ensemble fini $A$ de cardinal $n$ est
associée \emph{via} le $G$-torseur étale $U$ un fibré vectoriel $\xi$ de
rang $n$ sur $S$. À ces données, nous associerons dans la
section~\ref{section-operation-totale} un morphisme de préfaisceaux $K_r\to
K_{\xi^r}$, ce qui donnera lieu à l'opération cohomologique annoncée plus
haut. Une partie de l'intérêt est qu'en fait, à tout fibré vectoriel $E$
sur $S$ est associée un morphisme $K_E\to K_{E\otimes \xi}$ dans
$\Hopt[S]$, le cas précédent étant celui où $E$ est le
fibré trivial de rang $r$ sur
$S$. Il faut penser à cette opération comme à l'\guil{élévation à la
puissance $n$ tordue par le $G$-torseur $U$}. Si on dispose d'un autre
groupe fini $G'$, d'un $G'$-ensemble
$A'$, un $G'$-torseur étale $U'$
donnant lieu à un autre fibré vectoriel $\xi'$ sur $S$, il
devient alors loisible de composer l'opération associée à $A$ et celle
associée à $A'$ : \[K_E\to K_{E\otimes \xi}\to K_{E\otimes
\xi\otimes\xi'}\;\text{.}\] Le morphisme composé correspondra à la
construction associée au $(G\times G')$-ensemble produit $A\times A'$ et au
$(G\times G')$-torseur $U\times_S U'$
(cf.~proposition~\ref{proposition-compatibilite-composition}).
Cette observation nous permettra de
réinterpréter le théorème de symétrie (\ref{theoreme-symetrie}) ayant pour
conséquence les relations d'Adem (\S{}\ref{subsection-adem}). Les
relations d'Adem modulo $2$ données dans \cite[Theorem~10.2]{voevodsky-reduced}
comportaient quelques erreurs ; elles seront corrigées dans le
théorème~\ref{theoreme-adem-modulo-2}.

Les propriétés générales des opérations totales $H^{2r,r}(\mathcal X)\to
H^{2rn,rn}(\mathcal X\times (G\backslash U))$ pour $X\in\Ho[k]$ avec $k$ un
corps parfait seront étudiées dans la
section~\ref{section-proprietes-operation-totale}. Le
paragraphe~\S\ref{subsection-comparaison-voevodsky} fera le lien avec la
construction de Voevodsky et expliquera en quoi le fait qu'elle soit bien
définie était incorrectement justifié dans \cite{voevodsky-reduced}. Il est
tentant de penser que cette erreur provenait d'une contrainte d'écriture
consistant à n'utiliser d'autre catégorie homotopique $\Ho[S]$ que celle du
corps de base $k$. L'utilisation des catégories $\Ho[S]$ pour $S\in\Sm[k]$
permet de mieux saisir ce qui est véritablement important. Cela sera
également utile dans le \S\ref{subsection-annulation-bockstein} où nous
montrerons l'annulation du Bockstein évalué sur l'opération modulo $\ell$
associée au groupe symétrique $\Sym \ell$. La démonstration de ce fait
était incorrecte dans \cite{voevodsky-reduced}. Nous montrerons qu'il est
possible de le montrer en suivant d'assez près les arguments originaux de
Steenrod \cite{steenrod}.

Dans la section~\ref{section-operations}, nous définissons les opérations
de Steenrod proprement dites $P^i$ et $B^i$ qui sont des opérations
cohomologiques stables de bidegrés respectifs $(2i(\ell-1),i(\ell-1))$ et
$(2i(\ell-1)+1,i(\ell-1))$. Certains énoncés concernant les classifiants de
$\mu_\ell$ et $\Sym \ell$ ont été formulés en utilisant le langage des
motifs plutôt que celui de la cohomologie motivique utilisé dans
\cite{voevodsky-reduced}. On s'exonère en effet ici de ce qui semble avoir
été une autre contrainte d'écriture pour Voevodsky : ne pas utiliser les
catégories triangulées de motifs comme $\DMmoins[k]$ et \emph{a fortiori}
les catégories $\DMNeg[S]$ que nous utilisons dans la
section~\ref{section-thom-relatif}. Nous donnerons deux démonstrations
d'une formule-clef (cf.~proposition~\ref{proposition-p-c-1}) décrivant
l'action des opérations de Steenrod sur les classes de Chern de fibrés en
droites.

Dans la section~\ref{section-algebre-steenrod}, nous étudions l'algèbre de
Steenrod $A^{\star,\star}$ engendrée par les opérations de Steenrod, ainsi
que sa duale $A_{\star,\star}$. Quelques arguments ont été précisés par
rapport à \cite{voevodsky-reduced} ; la présentation que nous avons choisie
met aussi l'accent de façon peut-être plus systématique sur les différentes
structures utilisées. J'espère que la tâche du lecteur en sera facilitée.

Enfin, dans la section~\ref{section-endomorphismes}, on obtient un
résultat \guil{nouveau} exprimant que le spectre $\HZl$ représentant la
cohomologie motivique modulo $\ell$ n'admet pas d'endomorphisme stablement
fantôme non nul dans la catégorie homotopique stable $\SH[k]$, quand $k$
est un corps de caractéristique zéro. Ceci utilise de façon essentielle le
théorème principal de \cite{voevodsky-eilenberg-maclane} qui énonce qu'à
opérations stablement fantômes près, l'anneau bigradué des endomorphismes
du spectre $\HZl$ s'identifie à $A^{\star,\star}$. Grâce à ce résultat, on
peut maintenant identifier $A^{\star,\star}$ à l'anneau des endomorphismes
de $\HZl$.

\bigskip

Je remercie Denis-Charles Cisinski, Frédéric Déglise, Vincent Maillot et
Jörg Wildeshaus de m'avoir donné l'occasion d'approfondir ce sujet en me
demandant de faire un cours sur les opérations de Steenrod lors de l'école
d'été \emph{Motives and Milnor conjecture} à l'Institut de Mathématiques de
Jussieu en juin 2011. Je remercie aussi Vladimir Voevodsky d'avoir répondu
à certaines de mes questions.

\section{Isomorphisme de Thom relatif}
\label{section-thom-relatif}
\setcounter{subsubsection}{0}
Dans ce texte, on supposera implicitement que les schémas sont séparés. Un
schéma régulier sera pour nous un schéma noethérien (séparé) dont les
anneaux locaux sont réguliers au sens habituel. Si $S$ est un schéma
noethérien, la catégorie $\Sm$ désignera la catégorie des $S$-schémas
lisses (et séparés) de type fini. Une composante irréductible d'un schéma
sera toujours munie de la structure de sous-schéma fermé réduit. Le point
générique d'un schéma intègre $S$ sera noté $\eta_S$. Si $S$ est un schéma
et $s\in S$, on note $S_{(s)}=\Spec \OO_{S,s}$ le localisé (Zariski) de $S$
en $s$. L'ensemble des points de codimension $i$ d'un schéma noethérien $S$
est noté $S^{(i)}$.

\begin{definition}\cite[Corollary~3.4.6, Chapter~II]{livre-orange}
Soit $S$ un schéma régulier. Soit $\Lambda$ un groupe abélien de coefficients.
Soit $X\in \Sm$. On note $\cequi(X/S,0)$ (ou $\cequi(X/S,0)_\Lambda$
s'il y a lieu de préciser le groupe de coefficients) le
$\Lambda$-module libre sur l'ensemble des sous-schémas fermés intègres $Z$ de
$X$ tels que $Z\to S$ soit un morphisme fini et dont l'image (ensembliste)
soit une composante connexe de $S$.
\end{definition}

\begin{definition}
Soit $S$ un schéma régulier. Soit $E$ un fibré vectoriel sur $S$ de rang
$r$. On note $K\tM(\Th_S E)$
le préfaisceau d'ensembles pointés sur $\Sm$
qui à $U\in\Sm$ associe le $\Lambda$-module quotient 
\[\cequi(U\times_S E/U,0)/\cequi(U\times_S (E-\{0\})/U,0)\]
où $E-\{0\}$ est le fibré vectoriel épointé, c'est-à-dire le complémentaire
de la section nulle de $E$.
On rappellera dans le \S\ref{subsection-changement-de-base-cequi} comment
est définie la structure de préfaisceau.
\end{definition}

Si $S$ est lisse sur un corps parfait, la cohomologie motivique en bidegré
$(2r,r)$ est représentée dans $\Hopt$ par $K\tM(\Th_S \mathbf{A}^r)$ (pour
ainsi dire par définition). Un des buts principaux de cette section est de
montrer qu'au-dessus d'une telle base, $K\tM(\Th_S E)$ ne dépend que du
rang du fibré vectoriel $E$. On obtiendra ainsi différents préfaisceaux
représentant la cohomologie motivique et on les utilisera dans la
construction des opérations de Steenrod.

\subsection{Le changement de base pour $\cequi(X/S,0)$ avec $S$ régulier}
\label{subsection-changement-de-base-cequi}

\begin{definition}
Soit $X$ un schéma noethérien.
On note $\Dbcoh(X)$ la catégorie dérivée bornée de la catégorie abélienne
des faisceaux cohérents sur $X$. Soit $Z$ un fermé de $X$.
On note $\Dbcoh(X)_Z$ la sous-catégorie pleine de $\Dbcoh(X)$ formée des
complexes dont les faisceaux de cohomologie sont supportés par $Z$.
Soit $C$ une composante irréductible de $Z$ dont on note $\eta=\eta_C$ le
point générique. Le foncteur de restriction $\Dbcoh(X)_Z\to
\Dbcoh(X_{(\eta)})_\eta$ induit après passage aux groupes de Grothendieck
(cf.~SGA~5~VIII~2) une application $K(\Dbcoh(X)_Z)\to
K(\Dbcoh(X_{(\eta)})_\eta)\simeq \mathbf{Z}$ (notée $\chi_C$
ou $\chi_\eta$), l'isomorphisme de droite
étant donné par la somme alternée des longueurs des objets de cohomologie.
\end{definition}

\begin{definition}\label{definition-cb-cequi}
Soit $T\vers f S$ un morphisme de type fini entre schémas réguliers. Soit
$X\in\Sm$. On pose $X_T=X\times_S T$. On note $f^\star\colon
\cequi(X/S,0)\to \cequi(X_T/T,0)$ le morphisme de groupes abéliens qui à
$Z$ (avec $Z\subset X$ un sous-schéma fermé intègre tel que $Z\to S$ soit
un morphisme fini et surjectif (au-dessus d'une composante connexe de $S$)
associe
\[f^\star Z=\sum_{\eta\in (Z\times_S T)^{(0)}}
\chi_\eta(\OO_Z\otimes^\L_{\OO_S} \OO_T) \cdot C\]
\end{definition}

Il convient de vérifier que le cycle ainsi défini appartient bien à
$\cequi(X_T/T,0)$ :

\begin{lemme}\label{lemme-cb-cequi}
Les notations étant les mêmes que dans la définition précédente, le cycle
$f^\star(Z)$ est bien un élément de $\cequi(X_T/T,0)$. Plus précisément,
pour toute composante irréductible $C$ de $Z\times_S T$
le morphisme de projection $C\to T$ est fini et $C$ se
surjecte sur une composante connexe de $T$.
\end{lemme}

On rappelle (cf.~\cite{stroh-pilloni})
qu'une fonction de dimension (Zariski) sur un schéma noethérien
$U$ est une fonction $\delta\colon U\to\mathbf{Z}$ telle que pour toute
spécialisation immédiate (Zariski) $x\leadsto y$ (c'est-à-dire que $x$ et
$y$ sont deux points de $U$, que $y$ appartient à l'adhérence de $x$ et que
l'adhérence du point correspondant à $x$ dans le localisé $X_{(y)}$ est de
dimension $1$), on a $\delta(y)=\delta(x)-1$. Sur un schéma noethérien
connexe, deux fonctions de dimension diffèrent d'une constante. Si $U$ est
régulier, $-\codim$ est une fonction de dimension. Si $V\vers f U$ est un
morphisme de type fini, que $U$ est universellement caténaire et que
$\delta_U$ est une fonction de dimension sur $U$, alors on peut définir une
fonction de dimension $\delta_V$ sur $V$ en posant
$\delta_V(v)=\delta_U(f(v))+\degtr (v/f(v))$ où $\degtr$ est le degré de
transcendance\;\footnote{Pour vérifier cette affirmation,
on peut se ramener au cas où
$V=\mathbf{A}^1_U$ et $U$ est intègre universellement caténaire.}. (En
particulier, la restriction d'une fonction de dimension à un sous-schéma
est encore une fonction de dimension.)

Démontrons le lemme. On munit $S$ d'une fonction de dimension $\delta_S$.
Tout $S$-schéma de type fini sera muni de la fonction de dimension déduite
de $\delta_S$ par la règle ci-dessus. La propriété à établir est de nature
locale sur $S$, donc on peut supposer que $f=\pi\circ i$ avec $i\colon T\to
S'$ une immersion fermée et $\pi\colon S'\to S$ un morphisme lisse. On peut
aussi supposer que $S$, $S'$ et $T$ sont connexes. On note $c$ la
codimension de l'immersion fermée (régulière) $i$. Notons $C$ une
composante irréductible de $Z_T=Z_{S'}\cap X_T$, où l'intersection
schématique est
prise dans le schéma régulier $X_{S'}$. Notons $d'$ la dimension relative
du morphisme lisse $S'\to S$ et $C'$ une composante
irréductible de $Z_{S'}$ contenant $C$. Bien sûr, $C$ est une composante
irréductible de l'intersection $C'\cap X_T$ dans $X_{S'}$. En appliquant
\cite[Théorème~3, \S{}V.B.6]{algebre-locale} à cette intersection, on
obtient :
\[\delta_C(\eta_C)=\delta_{X_{S'}}(\eta_C)\geq \delta_{X_{S'}}(\eta_{C'})-c
=\delta_Z(\eta_Z)+d'-c=\delta_S(\eta_S)+d'-c\;\text{.}\]
Par ailleurs, si on note $t\in T$ l'image de $\eta_C$ par la projection
$C\to T$ qui est un morphisme fini, on a :
\[\delta_C(\eta_C)=\delta_T(t)=\delta_T(\eta_T)-\codim_T
t=\delta_S(\eta_S)+d'-c-\codim_T t\leq \delta_S(\eta_S)+d'-c\;\text{.}\]
En combinant ces inégalités, on obtient que
$\delta_C(\eta_C)=\delta_S(\eta_S)+d'-c$ et que $\codim_T t=0$,
c'est-à-dire que $t=\eta_T$. Autrement dit, le sous-schéma fermé intègre
$C$ de $X\times_S T$ définit bien un élément de $\cequi(X_T/T,0)$.

La proposition suivante montre que la construction du changement de base
est fonctorielle :

\begin{proposition}\label{proposition-fonctorialite-cb-cequi}
Considérons un diagramme de deux morphismes de type fini composables
$U\vers g T\vers f S$ entre schémas réguliers.
Soit $X\in \Sm$.
Alors, on a l'égalité $(f \circ g)^\star=g^\star\circ f^\star\colon
\cequi(X/S,0)\to \cequi(X_U/U,0)$.
\end{proposition}

\begin{lemme}
Soit $T\vers f S$ un morphisme de type fini entre schémas réguliers.
Soit $X\in\Sm$. Soit $Z$ un sous-schéma fermé
de $X$ dont les composantes irréductibles $(Z_i)_{i\in I}$ soient
telles que pour tout $i\in I$,
$Z_i\to S$ soit un morphisme fini et surjectif au-dessus d'une
composante connexe de $S$. Alors, pareillement, pour toute composante
irréductible $C$ de $Z_T=Z\times_S T$, le morphisme $C\to T$ est fini et
surjectif au-dessus d'une composante irréductible de $T$,
et le diagramme suivant commute, où les morphismes verticaux notés $\chi$
s'identifient à des sommes de morphismes $\chi_C\cdot C$
pour $C$ parcourant les composantes irréductibles de $Z$ (resp. $Z_T$) :
\[\xymatrix{\ar[d]^\chi K(\Dbcoh(X)_Z)\ar[r]^-{-\otimes^\L_{\OO_S}\OO_T }
& K(\Dbcoh(X_T)_{Z_T})\ar[d]^\chi \\
\cequi(X/S,0)\ar[r]^-{f^\star} & \cequi(X_T/T,0)}\]
\end{lemme}

L'assertion concernant les composantes irréductibles de $Z_T$ a été établie
dans le lemme~\ref{lemme-cb-cequi}. Par construction, la compatibilité à
établir pour tout élément $x\in K(\Dbcoh(X)_Z)$ est satisfaite pour les
éléments $[\OO_{Z_i}]$. Ceci permet de supposer que $x$ est dans le noyau de l'application
$K(\Dbcoh(X)_Z)\to \cequi(X/S,0)$. 

Autrement dit, il s'agit de montrer que si $x$ appartient à
$\ker(K(\Dbcoh(X)_Z)\vers\chi\cequi(X/S,0))$, alors l'image de $x$
\emph{via} le foncteur $-\otimes^{\L}_{\OO_S} \OO_X$ appartient à
$\ker(K(\Dbcoh(X_T)_{Z_T})\vers\chi\cequi(X_T/T,0))$.
De façon évidente,
cette condition sur $f$ se comporte bien par composition. Comme
on peut raisonner localement pour la topologie de Zariski, on peut utiliser
une factorisation de $f$ sous la forme d'une immersion fermée suivie d'un
morphisme lisse. Il suffit de traiter chacun de ces deux cas.

Par un dévissage sur les groupes de Grothendieck des faisceaux cohérents,
on se ramène à montrer que si $D$ est une composante irréductible de $Z$
(c'est-à-dire un des $Z_i$) et
$Y\subset D$ un sous-schéma fermé intègre strictement contenu dans $D$,
alors $\chi_C(\OO_Y\otimes^\L_{\OO_S}\OO_T)=0$ pour toute composante
irréductible $C$ de $Z_T$. Si $f\colon T\to S$ est lisse, c'est évident. Si
$f\colon T\to S$ est une immersion fermée, cela résulte de
\cite[Théorème~1~(a), \S{}V.C.1]{algebre-locale} dans le cas où $S$ et $T$
sont des schémas de type fini sur un corps, et de
\cite[Corollary~5.6]{gillet-soule} dans le cas général.

\medskip

Pour démontrer la proposition~\ref{proposition-fonctorialite-cb-cequi}, on
peut appliquer le lemme précédent à $f$, $g$ et $f\circ g$.
L'isomorphisme évident $(-\otimes^\L_{\OO_S} \OO_T)
\otimes^{\L}_{\OO_T}\OO_U\simeq -\otimes^\L_{\OO_S} \OO_U$ de foncteurs
$\Dbcoh(X)_Z\to \Dbcoh(X_U)_{Z_U}$ pour un certain choix de $Z$ permet
d'obtenir l'identité voulue $g^\star(f^\star(x))=(f\circ g)^\star(x)$ pour
tout élément $x$ appartenant à l'image du morphisme $\chi\colon
K(\Dbcoh(X)_Z)\to \cequi(X/S,0)$. Il est bien évidemment possible
d'agrandir $Z$ pour obtenir le résultat pour un élément arbitraire de
$\cequi(X/S,0)$.

\subsection{La catégorie $\DMNeg[S]$ pour un schéma régulier $S$}

Le but de cette sous-section est de définir une catégorie $\DMNeg[S]$ de
\guil{motifs} au-dessus de $S$. Dans le cas où $S=\Spec k$ est le spectre
d'un corps parfait, la catégorie $\DMNeg[k]$ s'identifiera à une
sous-catégorie pleine de la catégorie $\DMmoins[k]$ définie par Voevodsky
\cite[Chapter~V]{livre-orange} (cf.
\cite[Theorem~1.15]{voevodsky-eilenberg-maclane}). On a choisi une
définition qui permette d'obtenir de façon aussi simple que possible
l'adjonction entre cette catégorie $\DMNeg$ et la catégorie homotopique
$\Hopt$ (cf. sous-section~\ref{subsection-adjonction}). Il serait également
possible de procéder à ces constructions en utilisant la théorie radditive
(sic) de \cite{voevodsky-radditive} comme dans
\cite[\S{}1.2]{voevodsky-eilenberg-maclane}.

On présente ici la théorie d'une façon aussi concise que possible afin
d'établir le théorème d'isomorphisme de Thom. Le lecteur intéressé trouvera
de plus amples détails sur des constructions voisines dans
\cite{deglise-regular} et \cite{cisinski-deglise}.

\begin{definition}
Soit $S$ un schéma régulier. Si $X$ et $Y$ sont des objets de $\Sm$,
on note $\Cor_S(X,Y)=\cequi(X\times_S Y/X,0)$.
Supposons que $X$, $Y$ et $Z$ soient des objets de $\Sm$,
que $\alpha\in \Cor_S(X,Y)$ et que $\beta\in\Cor_S(Y,Z)$. Notons
$p_{23},p_{13},p_{12}$ les morphismes de projection respectifs de
$X\times_S Y\times_S Z$ sur $Y\times_S Z$, $X\times_S Z$ et $X\times_S Y$.
On note $\beta\circ \alpha=p_{13\star}(p_{12}^\star \alpha\cdot
p_{23}^\star\beta)\in \Cor_S(X,Z)$. Ceci définit la loi de composition pour
une catégorie additive notée $\SmCor$ dont les objets sont ceux de
$\Sm$ et dont les groupes d'homomorphismes sont les groupes $\Cor_S(X,Y)$
pour $X,Y\in\Sm$. On notera parfois $[X]$ l'objet de $\SmCor$ associé à un
objet $X$ de $\Sm$.
\end{definition}

Il convient d'indiquer pourquoi cette définition a un sens. Plus
précisément, le cycles $p_{12}^\star \alpha$ (resp. $p_{23}^\star \beta$) sur
$X\times_SY\times_S Z$ est définis par image inverse par le morphisme
plats $p_{12}$ (resp. $p_{23}$) de la façon usuelle. Il s'agit aussi de
l'élément de $\cequi(X\times_S Y\times_S Z/X\times_S Z,0)$ (resp.
$\cequi(X\times_S Y\times_S Z/X\times_S Y,0)$) obtenu par le changement de
base du \S\ref{subsection-changement-de-base-cequi}.
Dans la situation où $X$,
$Y$, $Z$ et $S$ sont connexes (cas auquel on se ramène facilement), on
montre facilement que si $G$ (resp. $F$) est une composante irréductible
d'un sous-schéma fermé de $X\times_S Y\times_S Z$ intervenant dans
l'écriture de $p_{12}^\star \alpha$ (resp. $p_{23}^\star \beta$), alors la
codimension de $G$ (resp. $F$) est la dimension relative de $Y\to S$ (resp.
$Z\to S$). En procédant comme dans la démonstration du
lemme~\ref{lemme-cb-cequi}, on montre ensuite que $G$ et
$F$ s'intersectent proprement et
que toutes les composantes irréductibles $C$ de $G\cap F$ dominent $X$ (en
fait, $C\to X$ est fini surjectif). On définit le produit d'intersection
$p_{12}^\star \alpha\cdot p_{23}^\star\beta$ en utilisant pour
multiplicités les caractéristiques d'Euler de modules $\Tor$ comme dans
\cite[\S{}V.C.1]{algebre-locale}. Le fait que les schémas $C$ considérés plus haut
soient finis sur $C$ permet de définir l'image directe de ce cycle, et
cette image est bien dans $\Cor_S(X,Y)$. L'associativité de cette loi de
composition provient de l'associativité (non triviale) du produit
d'intersection défini comme ci-dessus, et à vrai dire, comme les points
génériques des fermés intervenant dans l'écriture des cycles sont tous dans
la fibre des schémas considérés au-dessus du point générique de $S$, les
résultats énoncés dans \cite{algebre-locale} pour les schémas réguliers
d'égale caractéristique sont suffisants pour construire cette catégorie
$\SmCor$.

\bigskip

En associant à un morphisme $f\colon X\to Y$ dans $\Sm$ la
correspondance finie donnée par le cycle du graphe de $f$, on obtient
trivialement un foncteur $\Sm\to\SmCor$.

\begin{definition}
On note $\PST[S]$ la catégorie des foncteurs additifs $\SmCor^\opp\to
\Ab$. C'est la catégorie des préfaisceaux avec transferts. On dispose d'un
foncteur d'oubli des transferts $\PST[S]\to \Sm^\opp\Enspt$ vers la
catégorie des préfaisceaux d'ensembles (pointés) sur $\Sm$. Un
préfaisceau avec transferts qui devient un faisceau pour la topologie de
Nisnevich après application de ce foncteur d'oubli est appelé
\guil{faisceau (Nisnevich) avec tranferts}.
On note $\ST[S]$ la sous-catégorie pleine
de $\PST[S]$ formée des faisceaux avec transferts.
\end{definition}

On rappelle que la topologie de Nisnevich est la topologie
\guil{hensélienne} définie dans \cite{kato-saito}.

\begin{definition}
Pour tout $X\in\Sm$ et tout groupe abélien $\Lambda$, on note
$\Ltr(X)\in\PST[S]$ le préfaisceau avec transferts défini par :
\[\Ltr(X)(U)=\Hom_{\SmCor}([U],[X])=\Cor_S(U,X)=\cequi(U\times_S X/U,0)\]
pour tout $U\in\SmCor$.
\end{definition}

On vérifie facilement que $\Ltr(X)$ appartient à $\ST[S]$. (En fait, le
préfaisceau sur $\Sm$ induit par $\Ltr(X)$ est un faisceau pour la
topologie étale.)

\begin{lemme}
Le foncteur d'inclusion $\ST[S]\to\PST[S]$ admet un adjoint à gauche
$a_\Nis\colon \PST[S]\to \ST[S]$ qui commute avec l'oubli des tranferts.
Pour $\mathcal P\in \PST[S]$, on appelle $a_\Nis\mathcal P$ le faisceau
avec tranferts associé à $\mathcal P$.
\end{lemme}

Le fait que cet adjoint à gauche commute à l'oubli des transferts justifie que
l'on n'introduise pas une notation supplémentaire. Pour établir ce lemme, il
s'agit essentiellement de munir de transferts le faisceau associé au
préfaisceau sur $\Sm$ induit par $\mathcal P$ après oubli des tranferts. On
note $a_{\Nis}\mathcal P$ ce faisceau associé. On présente ici une variante des
arguments de \cite[Lemma~3.1.6, Chapter~V]{livre-orange} et
\cite[\S{}2.3]{deglise-regular}.

Pour $X\in\Sm$, tout élément $x\in (a_\Nis \mathcal P)(X)$ est induit
par un élément $\tilde x\in \ker (p^\star-q^\star\colon \mathcal P(U)\to
\mathcal P(V))$ où $(a\colon U \to  X,(p,q)\colon V\to U\times_X U)$ sont
deux morphismes étales couvrants pour la topologie de Nisnevich (autrement
dit, $U\to X$ et $V\to U\times_X U$ constituent un $1$-hyper-recouvrement
Nisnevich de $X$). On montre comme dans
\cite[Proposition~3.1.3, Chapter~V]{livre-orange} que le complexe de
préfaisceaux avec tranferts
\[\Ztr(V)\vers{p-q}\Ztr(U)\to \Ztr(X)\to 0\]
induit après oubli des transferts une suite exacte de faisceaux abéliens
sur $\Sm_\Nis$. L'élément $\tilde{x}$ définit un morphisme de
préfaisceaux avec transferts $\Ztr(U)/\image \Ztr(V)\to \mathcal P$, qui
après oubli des transferts et passage au faisceau associé, définit d'après
ce qui précède un morphisme $\Ztr(X)\to a_{\Nis}\mathcal P$
de faisceaux abéliens sur $\Sm_\Nis$. On vérifie facilement que ce
morphisme ne dépend pas du représentant $\tilde{x}$ de $x$.

La donnée de ce morphisme $\Ztr(X)\to a_\Nis\mathcal P$ décrit l'action des
tranferts sur $x$. On observe que ceci munit $a_\Nis\mathcal P$ d'une
structure de faisceau avec transferts et que $\mathcal P\to a_\Nis\mathcal
P$ est bien le morphisme universel de $\mathcal P$ vers un faisceau avec
transferts.

\bigskip

Il résulte aussitôt de ce lemme que non seulement $\PST[S]$ mais aussi
$\ST[S]$ sont des catégories abéliennes de Grothendieck, tout comme leurs
versions $\PST$ et $\ST$ à coefficients dans un
anneau commutatif $\Lambda$.

\begin{definition}\label{definition-categories-de-modeles}
Soit $\mathcal C$ une catégorie abélienne de Grothendieck. D'après
\cite{hovey}, on peut définir une structure de catégorie de modèles sur la
catégorie (non bornée) des complexes $\Comp(\mathcal C)$ à valeurs dans
$\mathcal C$ de façon à ce que les équivalences faibles soient les
quasi-isomorphismes, les cofibrations les monomorphismes et les fibrations
les morphismes ayant la propriété de relèvement à droite par rapport aux
cofibrations triviales.

On note $\CompNeg(\mathcal C)$ la sous-catégorie pleine de $\Comp(\mathcal C)$
formée des complexes situés en degrés \emph{homologiques} positifs,
c'est-à-dire des complexes de la forme $\dots \depuis 0 \depuis K_0\depuis K_1
\depuis K_2 \depuis \dots$. Le foncteur d'inclusion $\CompNeg(\mathcal C)\to
\Comp(\mathcal C)$ admet pour adjoint à droite le foncteur de troncature
$\tau_{\geq 0}$ qui à un complexe $K\in\Comp(\mathcal C)$ associe son
sous-complexe $\dots\depuis 0\depuis \ker (d\colon K_0\to K_{-1})\depuis
K_1\depuis K_2\depuis \dots$.

On munit $\CompNeg(\mathcal C)$ de la structure de catégorie de modèles dont
les cofibrations sont les monomorphismes, les équivalences faibles les
quasi-isomorphismes et les fibrations les morphismes ayant la propriété de
relèvement à droite par rapport aux cofibrations triviales (de
$\CompNeg(\mathcal C)$). En outre, le couple de foncteurs adjoints formé de
l'inclusion de $\CompNeg(\mathcal C)$ dans $\Comp(\mathcal C)$ et le foncteur
de troncature $\tau_{\geq 0}$ est une adjonction de Quillen.

On note $\Der\mathcal C$ et $\DerNeg\mathcal C$ les catégories homotopiques
de ces catégories de modèles.
\end{definition}

La structure de catégorie de modèles sur $\CompNeg(\mathcal C)$ se déduit
presqu'immédiatement de celle définie par Hovey \cite{hovey} sur
$\Comp(\mathcal C)$. Le seul axiome qui pourrait poser une difficulté est
le fait que dans $\CompNeg(\mathcal C)$ les cofibrations ont bien la
propriété de relèvement à gauche par rapport aux fibrations triviales. Ceci
peut se vérifier en utilisant l'astuce de Joyal
(cf.~\cite[p.~64-65]{jardine}). (On
prendra garde au fait qu'une fibration de $\CompNeg(\mathcal C)$ n'est pas
forcément une fibration de $\Comp(\mathcal C)$.)

\bigskip

Pour $\mathcal C$ parmi les catégories $\PST$, $\ST$
et $\Ab$ (la catégorie des groupes abéliens), on munit les catégories
$\Comp(\mathcal C)$ et $\CompNeg(\mathcal C)$ des structures de catégories
de modèles de la définition~\ref{definition-categories-de-modeles}.

\begin{definition}
Si $f\colon K\to L$ est un morphisme dans $\Comp(\Ab)$, on note
$\cocone(f)=\cone(f)[-1]$. (Le cône est un modèle de la cofibre homotopique
de $f$. Dualement, le cocône est un modèle de la fibre homotopique de $f$ :
il s'insère dans un triangle distingué \[\cocone(f)\to K\to L\to
\cocone(f)[1]\;\text{.}\]
\end{definition}

\begin{definition}
Soit
un carré commutatif dans
$\Comp(\Ab)$ :
\[\xymatrix{A\ar[r] \ar[d]& B\ar[d]^f \\ B' \ar[r]_{f'}& C}\]
On dit que ce carré est homotopiquement cartésien si le
morphisme évident $A\to \cocone(B\oplus B'\vers{f-f'} C)$ est un
quasi-isomorphisme.

Supposons que $A$, $B$, $B'$ et $C$ soient dans $\CompNeg(\Ab)$. On dit que
le carré est homotopiquement cartésien dans $\CompNeg(\Ab)$ si le morphisme
évident $A\to \tau_{\geq 0}\cocone(B\oplus B'\vers{f-f'} C)$ est un
quasi-isomorphisme.
\end{definition}

On vérifie facilement que les catégories de modèles $\Comp(\Ab)$ et
$\CompNeg(\Ab)$ sont propres (à gauche et à droite) et que les notions
axiomatiques des carrés homotopiquement cartésiens (cf.
\cite[\S{8}, Chapter~II]{goerss-jardine}) dans l'une et l'autre de ces
catégories de modèles coïncident avec celles définies ici. Le foncteur de
Quillen à droite $K\colon \CompNeg(\Ab)\simeq
\Delta^\opp\Ab\to\Delta^\opp\Ens$ donné par l'équivalence de Dold-Kan
(\cite[\S{2}, Chapter~III]{goerss-jardine}) préserve et reflète les
équivalences faibles. Ainsi, un carré dans $\CompNeg(\Ab)$ est
homotopiquement cartésien si et seulement si le carré obtenu en appliquant
le foncteur $K$ est homotopiquement cartésien dans $\Delta^\opp\Ens$.

Dans \cite[Definition~1.13]{morel-voevodsky}, la propriété de Brown-Gersten
est définie pour les préfaisceaux simpliciaux sur $\Sm$. On la définit ici
pour des préfaisceaux à valeurs dans $\CompNeg(\Ab)$ :

\begin{definition}
Soit $K$ un préfaisceau sur $\Sm$ à valeurs dans $\CompNeg(\Ab)$. On dit
que $K$ vérifie la propriété de Brown-Gersten si $K(\emptyset)$ est
acyclique et si pour tout carré distingué Nisnevich
\[\xymatrix{p^{-1}(U)\ar[r]\ar[d] & V \ar[d]^p \\ U\ar[r] & X}\]
c'est-à-dire que $U$ est un sous-schéma ouvert de $X\in\Sm$ et que
$p\colon V\to X$ est un morphisme étale induisant un isomorphisme au-dessus
du fermé complémentaire $Z=(X-U)_\red$, alors le carré commutatif
\[\xymatrix{K(X)\ar[r]\ar[d] & K(U)\ar[d] \\ K(V)\ar[r]& K(p^{-1}(U))}\]
qui s'en déduit dans $\CompNeg(\Ab)$ est homotopiquement cartésien.

Si $K\in\CompNeg(\PST)$,
on dit qu'il vérifie la propriété de Brown-Gersten si le
préfaisceau sur $\Sm$ déduit de $K$ par composition avec le foncteur
$\Sm\to \SmCor$ la vérifie.
\end{definition}

\begin{proposition}
Le foncteur $a_\Nis\colon \CompNeg(\PST)\to
\CompNeg(\ST)$ est un
foncteur de Quillen à gauche. On note $\Hyper\colon \DerNeg(\ST)\to
\DerNeg(\PST)$ le foncteur dérivé total à droite du foncteur d'inclusion 
$\CompNeg(\ST)\to \CompNeg(\PST)$ (c'est l'adjoint à droite du
foncteur $\DerNeg(\PST)\to \DerNeg(\ST)$ induit par $a_\Nis$).

Si $\mathcal P\in \CompNeg(\PST)$, alors le morphisme d'adjonction
$\mathcal P\to \Hyper a_\Nis\mathcal P$ est un isomorphisme dans
$\DerNeg(\PST)$ si et seulement si $\mathcal P$ vérifie la propriété de
Brown-Gersten.
\end{proposition}

Il est évident que $a_\Nis$ préserve les cofibrations et les équivalences
faibles. Le foncteur d'inclusion
$\CompNeg(\ST)\to \CompNeg(\PST)$ et son adjoint à droite constituent donc
une adjonction de Quillen.

Pour la deuxième assertion (qui est homologue à
\cite[Proposition~1.16, p.~100]{morel-voevodsky}),
on commence par observer que les objets fibrants
$\mathcal F$ de $\CompNeg(\ST)$ vérifient la propriété de Brown-Gersten.
Pour cela, on peut utiliser le critère \cite[Remark~1.15,
page~100]{morel-voevodsky} : comme $\mathcal F$ est un faisceau, il suffit
de vérifier que pour toute immersion ouverte $U\to V$ dans $\Sm$, le
morphisme de restriction $\mathcal F(V)\to \mathcal F(U)$ induit une fibration
d'ensembles simpliciaux après application du foncteur $K\colon
\CompNeg(\Ab)\to \Delta^\opp\Ens$.
Cette condition résulte d'un jeu d'adjonctions
standard à partir de l'observation que le morphisme $\Ltr(U)\to \Ltr(V)$
est un monomorphisme dans $\ST$. Pour conclure, on utilise une
résolution fibrante $\mathcal P\to \mathcal F$ et on applique la variante
Nisnevich du théorème de Brown-Gersten
\cite[Lemma~1.18, p.~101]{morel-voevodsky}.

\begin{corollaire}\label{corollaire-brown-gersten}
La catégorie $\DerNeg(\ST)$ est équivalente à la
catégorie homotopique
de la localisation à la Bousfield de la catégorie de modèles
$\CompNeg(\PST)$ dans laquelle les objets locaux
seraient ceux vérifiant la propriété de Brown-Gersten.
\end{corollaire}

On notera que d'après \cite[Proposition~A.4]{riou-sh}, les catégories de
modèles que nous considérons sont cellulaires ; il est donc possible de leur
appliquer les techniques de localisation développées dans
\cite{hirschhorn}.

\begin{definition}
Soit $\Lambda$ un anneau commutatif. La localisation à la Bousfield de la
catégorie de modèles $\CompNeg(\ST)$ par rapport aux morphismes
$\Ltr(X\times\mathbf{A}^1)\to\Ltr(X)$ définit une nouvelle catégorie de
modèles ($\mathbf{A}^1$-localisée)
dont les équivalences faibles et les fibrations seront appelées
$\mathbf{A}^1$-équivalences faibles et $\mathbf{A}^1$-fibrations.

On note $\DMNeg$ la catégorie homotopique de la structure de
catégorie de modèles $\mathbf{A}^1$-localisée de $\CompNeg(\ST)$.
\end{definition}

D'après le corollaire~\ref{corollaire-brown-gersten}, $\DMNeg$
est aussi équivalente à la catégorie homotopique de la localisation à la
Bousfield de $\CompNeg(\PST)$ dont les objets locaux $\mathcal
F$ seraient ceux qui vérifient la propriété de Brown-Gersten et tels que
$\mathcal F(X)\to \mathcal F(\mathbf{A}^1\times X)$ soit un
quasi-isomorphisme pour tout $X\in\Sm$.

\subsection{L'adjonction entre $\DMNeg[S]$ et $\Hopt$}
\label{subsection-adjonction}

\begin{definition}
On note $\Hosimppt$ (resp. $\Hopt$) la catégorie homotopique de la
structure de catégorie de modèles sur $\Delta^\opp\Sm^\opp\Enspt$ des
préfaisceaux simpliciaux (pointés) sur $\Sm$ munie de la structure définie
dans \cite{jardine} pour le site $\Sm_\Nis$ (resp. muni de la structure
$\mathbf{A}^1$-localisée définie dans \cite{morel-voevodsky}). (Dans tous
les cas, les cofibrations sont les monomorphismes.)
\end{definition}

\begin{proposition}\label{proposition-adjonction-non-a-1-localise}
On définit un foncteur d'oubli des transferts $\oub\colon
\CompNeg\ST\to\Delta^\opp\Sm^\opp\Enspt$ en combinant l'équivalence
de Dold-Kan $\CompNeg\ST\simeq \Delta^\opp\ST$
(cf.~\cite[\S{2}, Chapter~III]{goerss-jardine}) et le
foncteur d'oubli $\ST\to\Sm^\opp\Enspt$. Ce foncteur envoie les
quasi-isomorphismes sur des équivalences faibles (locales pour la topologie
de Nisnevich), ce qui induit un foncteur noté $\oub\colon \DerNeg(\ST)\to
\Hosimppt$. Ce foncteur admet un adjoint à gauche :
\[\tM\colon \Hosimppt\to \DerNeg(\ST)\]
\end{proposition}

Seule l'existence de l'adjoint à gauche est à justifier. Il est formel que
le foncteur d'oubli $\ST\to \Sm^\opp\Enspt$ admette un adjoint à gauche
$\tL\colon \Sm^\opp\Enspt\to\ST$. En effet, pour $(X,x)$ un objet pointé
$\Sm$, l'image par $\tL$ du préfaisceau représenté par $X$ et pointé par
$x$ est $\coker(\Ltr(S)\vers x\Ltr(X))$. Dans le cas général, on utilise
qu'un préfaisceau pointé est colimite de préfaisceaux représentables. On
étend $\tL$ en un foncteur $\Delta^\opp\Sm^\opp\Enspt\to \Delta^\opp\ST$ et
on note encore $\tL\colon \Delta^\opp\Sm^\opp\Enspt\to\CompNeg(\ST)$ le
foncteur obtenu par composition avec l'équivalence de Dold-Kan.

Pour $\mathcal X\in\Delta^\opp\Sm^\opp\Enspt$, un candidat
naturel pour
$\tM\mathcal X$ est $\tL\mathcal X$. On procède alors comme dans
\cite[p. 63-64]{morel-voevodsky} (voir aussi \cite[\S{}2.1]{riou-sh}). On
dit qu'un objet $\mathcal X$ est admissible si pour tout objet fibrant
$\mathcal F\in \CompNeg(\ST)$ et toute résolution fibrante $\oub \mathcal F\to
\mathcal G$ dans $\Delta^\opp\Sm^\opp\Enspt$, le morphisme évident
d'ensembles simpliciaux suivant est une équivalence faible :
\[\hom(\tL\mathcal X,\mathcal F)\to \hom(\mathcal X,\mathcal G)\;\text{.}\]
Ici, $\hom$ désigne les bifoncteurs $\Hom$ simpliciaux (i.e. à valeurs dans
$\Delta^\opp\Ens$) provenant des
structures simpliciales évidentes sur $\CompNeg\ST\simeq \Delta^\opp\ST$
et $\Delta^\opp\Sm^\opp\Enspt$. Du fait de l'adjonction (simpliciale) entre
$\oub$ et $\tL$, cela s'énonce aussi en disant que
\[\hom(\mathcal X,\oub \mathcal F)\to \hom(\mathcal X,\mathcal G)\]
est une équivalence faible d'ensembles simpliciaux.

En passant au $\pi_0$ l'équivalence faible ci-dessus, on obtient aussitôt
que si $\mathcal X$ est admissible, alors le foncteur adjoint $\tM$ est
défini en $\mathcal X$ et $\tM(\mathcal X)=\tL(\mathcal X)$. Montrer que
l'adjoint à gauche $\tM$ existe revient à montrer que tout objet de
$\Hosimppt$ peut être représenté par un préfaisceau simplicial pointé
admissible, ce que nous allons faire (et la démonstration montre qu'en fait
$\tM$ est le foncteur dérivé total à gauche de $\tL$).

En utilisant des lemmes de résolutions (voir 
\cite[Lemma~1.16, p.~52]{morel-voevodsky}), on peut supposer
que $\mathcal X$ est
un préfaisceau simplicial tel que pour tout $n\in\mathbf{N}$, ${\mathcal X}_n$
soit une somme directe de préfaisceaux (pointés) représentables. Pour
montrer qu'un tel préfaisceau $\mathcal X$ est admissible, en utilisant les
propriétés des objets admissibles (cf.
\cite[Lemma~1.53, p. 64]{morel-voevodsky} et
\cite[Lemmes~2.7 et 2.8]{riou-sh}), on peut se
ramener au cas où $\mathcal X$ est induit par un préfaisceau d'ensembles
représenté par un objet pointé $(X,x)$ de $\Sm$. Il s'agit de
montrer que l'application d'ensembles simpliciaux
\[\hom((X,x),\oub \mathcal F)\to \hom((X,x),\mathcal G)\]
est une équivalence faible. L'ensemble simplicial
$\hom((X,x),\oub \mathcal F)$ est la fibre du morphisme d'ensembles
simpliciaux pointés $\oub\mathcal F(X)\vers {x^\star}\oub\mathcal F(S)$
(de même pour $\mathcal G$).
Pour conclure, il suffit donc de montrer que $\oub\mathcal
F(X)\to \mathcal G(X)$ et $\oub\mathcal F(S)\to \mathcal G(S)$ sont des
équivalences faibles et que les morphismes $\oub\mathcal F(X)\to
\oub\mathcal F(S)$
et $\mathcal G(X)\to\mathcal G(S)$ induits par $x$ sont des fibrations. La
première assertion résulte du théorème de Brown-Gersten puisqu'aussi bien
$\mathcal F$ que $\mathcal G$ vérifient la propriété de Brown-Gersten. La
deuxième assertion se déduit facilement du fait que les objets $\mathcal
F$ et $\mathcal G$ sont fibrants dans $\CompNeg(\ST)$ et
$\Delta^\opp\Sm^\opp\Enspt$ respectivement.

\begin{lemme}
Soit $U\to X$ une immersion ouverte dans $\Sm$. Alors, le préfaisceau
pointé $X/U$ est admissible au sens introduit dans la démonstration de la
proposition~\ref{proposition-adjonction-non-a-1-localise}.
Autrement dit, on a un isomorphisme
$\tM(X/U)\simeq \Ltr(X)/\Ltr(U)$.
\end{lemme}

Ceci repose sur le principe \cite[Lemma~1.53(3), p.~65]{morel-voevodsky}
qui s'applique puisque $U\to X$ et $\Ltr(U)\to\Ltr(X)$ sont des
monomorphismes.

\bigskip

Pour qu'il n'y ait pas d'ambiguité sur le sens de l'énoncé de la
proposition~\ref{proposition-oubli-transferts}, on précise la définition
d'objet $\mathbf{A}^1$-local utilisée ici :

\begin{definition}\label{definition-complexe-a-1-local}
Si $\mathcal F\in\CompNeg(\ST)$, on dit que $\mathcal F$ est
$\mathbf{A}^1$-local si pour tout morphisme $\mathcal
F\vers{\varphi}\mathcal F'$
tel que $\mathcal F'$ soit fibrant et $\varphi$ une
équivalence faible (pour la structure non
$\mathbf{A}^1$-localisée), le morphisme $\mathcal F'(X)\to\mathcal
F'(\mathbf{A}^1\times X)$ est un quasi-isomorphisme pour tout $X\in\Sm$.
On dit que $\mathcal F\in\CompNeg(\PST)$ est $\mathbf{A}^1$-local si
$a_\Nis\mathcal F$ est $\mathbf{A}^1$-local au sens ci-dessus. (Noter que la
condition ne dépend que de la classe d'isomorphisme de l'image de $\mathcal
F$ dans la catégorie $\DerNeg(\ST)$. Noter aussi que cela a un sens de dire
d'un morphisme dans $\DerNeg(\ST)$ que c'est une $\mathbf{A}^1$-équivalence
faible.)
\end{definition}

On utilise une définition semblable des objets $\mathbf{A}^1$-locaux dans
$\Hosimppt$, voir \cite[Definition~3.1, p.~86]{morel-voevodsky}.

\begin{proposition}\label{proposition-oubli-transferts}
Le foncteur $\oub\colon \DerNeg(\ST)\to \Hosimppt$ envoie les objets
$\mathbf{A}^1$-locaux sur des objets $\mathbf{A}^1$-locaux.
Un morphisme $f$ dans $\DerNeg(\ST)$ est une $\mathbf{A}^1$-équivalence
faible si et seulement si $\oub f$ est une $\mathbf{A}^1$-équivalence
faible.
\end{proposition}

Si $\mathcal F$ est un objet fibrant de $\CompNeg(\ST)$ (pour la structure
non $\mathbf{A}^1$-localisée), alors $\mathcal F$ vérifie la propriété de
Brown-Gersten. Dire que $\mathcal F$ est $\mathbf{A}^1$-local revient à
dire que $\mathcal F(X)\to \mathcal F(\mathbf{A}^1\times X)$ est un
quasi-isomorphisme. Après application du foncteur $\oub$, $\oub \mathcal F$
vérifie encore bien sûr la propriété de Brown-Gersten. L'hypothèse que
$\mathcal F$ est $\mathbf{A}^1$-local montre alors que $\oub\mathcal
F(X)\to\oub\mathcal F(\mathbf{A}^1\times X)$ est une équivalence faible, ce
qui conjointement avec la propriété de Brown-Gersten, implique que
$\oub\mathcal F$ est $\mathbf{A}^1$-local.

Observons qu'il découle de ceci que si $\mathcal F\to \mathcal G$ est une
$\mathbf{A}^1$-équivalence faible entre objets $\mathbf{A}^1$-locaux, c'est
une équivalence faible, donc $\oub \mathcal F\to \oub\mathcal G$ est une
équivalence faible, \emph{a fortiori} une $\mathbf{A}^1$-équivalence faible.

Pour montrer que $\oub$ préserve les $\mathbf{A}^1$-équivalences faibles,
il suffit donc de montrer que pour tout $\mathcal F\in\CompNeg(\ST)$, pour
une certaine $\mathbf{A}^1$-résolution fonctorielle $\mathcal F'$, à savoir
que $\mathcal F\to\mathcal F'$ est une $\mathbf{A}^1$-équivalence faible et
$\mathcal F'$ est $\mathbf{A}^1$-local, alors $\oub \mathcal
F\to\oub\mathcal F'$ est une $\mathbf{A}^1$-équivalence faible.

Pour cela, on utilise une variante des constructions $\Sing=
\Sing^{\mathbf{A}^1}$
de \cite[p.~85--91]{morel-voevodsky} et $C_\star$ de
\cite[Proposition~3.2.3, Chapter~5]{livre-orange}.
On note $\mathbf{\Delta}^\bullet$ l'objet cosimplicial standard dans $\Sm$
formé d'espaces affines.
Soit $\mathcal X$ un préfaisceau simplicial sur $\Sm$.
Pour $U\in \Sm$, on
note $\SheafHom(U,\mathcal X)$ le préfaisceau simplicial qui à $V\in \Sm$
associe $\mathcal X(U\times_S V)$. L'application qui à $n$ associe
$\SheafHom(\mathbf{\Delta}^n,\mathcal X)$ définit un objet simplicial dans
la catégorie des préfaisceaux simpliciaux, autrement dit un préfaisceau
bisimplicial. On note $\Sing\mathcal
X\in\Delta^\opp\Sm^\opp\Enspt$ la diagonale de cet objet bisimplicial.
D'après \cite[Corollary~3.8, p. 89]{morel-voevodsky}, le morphisme évident
$\mathcal X\to \Sing\mathcal X$ est une
$\mathbf{A}^1$-équivalence faible.

La construction $\SheafHom(U,\mathcal X)$ mentionnée plus haut pour
$U\in\Sm$ admet une variante quand on remplace le préfaisceau simplicial
$\mathcal X$ par un objet $\mathcal F\in\PST$ (ceci est lié au produit
tensoriel sur $\SmCor$ étendant le produit dans $\Sm$). Cette construction
$\SheafHom(U,-)$ commute à l'oubli des tranferts. Ainsi, si $\mathcal
F\in\CompNeg(\ST)$, modulo l'identification $\CompNeg(\ST)\simeq
\Delta^\opp \ST$, on obtient un objet $\Sing\mathcal
F\in\Delta^\opp\ST\simeq \CompNeg(\ST)$. Pour aller plus loin, nous allons
utiliser le lemme suivant :

\begin{lemme}\label{lemme-a-1-homotopie}
Soit $\mathcal F\in \PST$. Alors, le morphisme évident
$\SheafHom(\mathbf{A}^1,\mathcal F)\to\mathcal F$ est une
$\mathbf{A}^1$-équivalence faible.
\end{lemme}

La démonstration de ce fait est étonnamment délicate. On note
$\Ltr(\mathbf{A}^1)\otimes-\colon \PST\to\PST$ le foncteur adjoint du
foncteur $\SheafHom(\mathbf{A}^1,-)\colon \PST\to\PST$. Il envoie $\Ltr(X)$ sur
$\Ltr(\mathbf{A}^1\times X)$ et par la compatibilité que ce foncteur doit
vérifier par rapport aux colimites, on l'étend tautologiquement à $\PST$.
Si $\mathcal X\in\PST$ est une somme directe de préfaisceaux de la forme
$\Ltr(X)$, alors $\Ltr(\mathbf{A}^1)\otimes \mathcal X\to\mathcal
X$ est tautologiquement une $\mathbf{A}^1$-équivalence faible parce qu'une
somme directe de $\mathbf{A}^1$-équivalences faibles est une
$\mathbf{A}^1$-équivalence faible. En utilisant une variante appropriée de
\cite[Proposition~2.14, page~74]{morel-voevodsky}, on obtient que si
$\mathcal X\in\CompNeg(\PST)$ est tel que pour tout $n\in\mathbf{N}$,
$\mathcal X_n$ soit une somme directe de préfaisceaux $\Ltr(X)$, alors
$\Ltr(\mathbf{A}^1)\otimes \mathcal X\to\mathcal X$ est une
$\mathbf{A}^1$-équivalence faible. Il en résulte que pour un tel objet
$\mathcal X$, les
deux morphismes $\mathcal X\to \Ltr(\mathbf{A}^1)\otimes \mathcal X$
définis par les $S$-points $0$ et $1$ de $\mathbf{A}^1_S$ deviennent égaux dans
la catégorie $\DMNeg$. Ceci vaut en fait pour tout $\mathcal
X\in\CompNeg(\PST)$ comme on le déduit facilement d'un
quasi-isomorphisme $\mathcal
X'\to\mathcal X$ dans $\CompNeg(\PST)$ avec $\mathcal X'$ de la
forme particulière ci-dessus.

Pour $\mathcal F\in\CompNeg(\PST)$, on considère ensuite le morphisme 
\[\varphi\colon \Ltr(\mathbf{A}^1)\otimes \SheafHom(\mathbf{A}^1,\mathcal F)\to
\SheafHom(\mathbf{A}^1,\mathcal F)\]
qui correspond par adjonction au morphisme
\[\SheafHom(\mathbf{A}^1,\mathcal F)\to
\SheafHom(\mathbf{A}^1,\SheafHom(\mathbf{A}^1,\mathcal
F))=\SheafHom(\mathbf{A}^2,\mathcal F)\]
qui provienne par fonctorialité du morphisme 
$\mathbf{A}^2\to\mathbf{A}^1$ qui à $(x,y)$ associe $xy$.

En composant ce morphisme $\varphi$ avec les deux morphismes 
$\SheafHom(\mathbf{A}^1,\mathcal F)\to
\Ltr(\mathbf{A}^1)\otimes \SheafHom(\mathbf{A}^1,\mathcal F)$ correspondant
à $0$ et $1$, on obtient que le morphisme composé évident
$\SheafHom(\mathbf{A}^1,\mathcal F)\vers {0^\star} \mathcal F\to
\SheafHom(\mathbf{A}^1,\mathcal F)$ devient l'identité dans $\DMNeg$, ce
qui permet de finir la démonstration du lemme~\ref{lemme-a-1-homotopie}.

\medskip

On peut maintenant finir la démonstration de la
proposition~\ref{proposition-oubli-transferts}. Grâce au lemme précédent,
en procédant comme dans le cas des préfaisceaux simpliciaux, on obtient que
pour tout $\mathcal F\in \CompNeg(\ST)$, le morphisme évident $\mathcal
F\to \Sing\mathcal F$ est une $\mathbf{A}^1$-équivalence
faible (c'est-à-dire devient un isomorphisme dans $\DMNeg$).

Notons $\Ex$ un foncteur de résolution fibrante sur $\CompNeg(\ST)$ (pour
la structure non $\mathbf{A}^1$-localisée) et $F=\Sing\circ
\Ex\colon \CompNeg(\ST)\to \CompNeg(\ST)$.
En utilisant la transformation naturelle $\Id\to F$ pour définir des
morphismes $F^n\to F\circ F^{n}=F^{n+1}$, on obtient un
système inductif $(F^n)_{n\geq 0}$ de foncteurs dont on note
$F^{\infty}$ la colimite. Pour tout $\mathcal F\in \CompNeg(\ST)$, le
morphisme $\mathcal F\to F^{\infty}\mathcal F$ est une
$\mathbf{A}^1$-équivalence faible et $F^{\infty}\mathcal F$ est
$\mathbf{A}^1$-local d'après l'argument de
\cite[Lemma~2.6, p.~107]{morel-voevodsky}. En observant que $\Sing$ commute
à l'oubli des transferts, on peut utiliser \cite[p.~107]{morel-voevodsky}
pour obtenir que le morphisme
$\oub\mathcal F \to \oub F^{\infty}\mathcal F$ est une
$\mathbf{A}^1$-équivalence faible dans $\Hosimppt$.

Pour montrer que $\oub$ préserve et réflète les $\mathbf{A}^1$-équivalences
faibles, en utilisant le résultat du paragraphe précédent, on peut supposer
que la source et le but de la flèche considérée dans $\CompNeg(\ST)$ sont
$\mathbf{A}^1$-locaux. Dans ce cas, les notions de
$\mathbf{A}^1$-équivalences faibles et d'équivalences faibles (locales pour
la topologie de Nisnevich) coïncident. On est ainsi ramené à vérifier que
si $f\colon \mathcal F\to\mathcal G$ est un morphisme dans $\CompNeg(\ST)$,
alors $f$ est un quasi-isomorphisme si et seulement si $\oub f$ est une
équivalence faible (locale pour la topologie de Nisnevich), ce qui est
évident.

\medskip

On déduit aussitôt des
propositions~\ref{proposition-adjonction-non-a-1-localise} et
\ref{proposition-oubli-transferts} le corollaire suivant :

\begin{corollaire}
Le foncteur $\oub\colon \CompNeg(\ST)\to \Delta^\opp\Sm^\opp\Enspt$
préservant les $\mathbf{A}^1$-équivalences faibles, il définit par
localisation un foncteur $K\colon \DMNeg\to \Hopt$. Le foncteur
$\tM\colon \Hosimppt\to \DerNeg(\ST)$ préservant les
$\mathbf{A}^1$-équivalences faibles, on note aussi $\tM\colon \Hopt\to
\DMNeg$ le foncteur qu'il induit et qui est adjoint à gauche de $K$.
\end{corollaire}

\subsection{La fonctorialité élémentaire de la catégorie $\DMNeg[S]$}

Dans cette sous-section, on se donne un morphisme $f\colon T\to S$ entre
schémas réguliers. On définit un foncteur $f^\star\colon \SmCor[S]\to
\SmCor[T]$ qui sur les objets envoie $X\in\Sm[S]$ sur $X_T=X\times_S
T\in\Sm[T]$ et qui au niveau des morphismes, pour $X$ et $Y$ dans $\Sm[S]$,
soit donné par l'application de changement de base (cf.
définition~\ref{definition-cb-cequi}) :
\[\cequi(X\times_S Y/X,0)\to \cequi(X_T\times_T Y_T/X_T,0)\;\text{.}\]
Il n'est pas bien difficile de vérifier que cette définition est compatible
à la composition des correspondances finies, c'est-à-dire que ceci définit
bien un foncteur $f^\star\colon \SmCor[S]\to \SmCor[T]$.

La composition avec le foncteur additif $f^\star\colon \SmCor[T]\to
\SmCor[S]$ définit un foncteur $f_\star\colon
\PST[T,\Lambda]\to\PST[S,\Lambda]$ (qui applique la sous-catégorie
$\ST[T,\Lambda]$ dans $\ST[S,\Lambda]$). Ce foncteur $f_\star\colon
\PST[T,\Lambda]\to\PST[S,\Lambda]$ admet un adjoint à gauche $f^\star\colon 
\PST[T,\Lambda]\to\PST[S,\Lambda]$ (qui étend par passage à la limite
inductive le foncteur $f^\star\colon \SmCor[T]\to
\SmCor[S]$ au niveau des objets représentables).

La catégorie abélienne
$\PST[S,\Lambda]$ admettant suffisamment d'objets projectifs,
il est permis d'introduire le foncteur dérivé total à gauche $\L
f^\star\colon \DerNeg(\PST[S,\Lambda])\to \DerNeg(\PST[T,\Lambda])$.

\begin{proposition}
Si $\rho\colon \CompNeg(\PST[T,\Lambda])\to\CompNeg(\PST[T,\Lambda])$
est une résolution $\mathbf{A}^1$-fibrante, le foncteur $f_\star\circ
\rho\colon \CompNeg(\PST[T,\Lambda])\to \CompNeg(\PST[S,\Lambda])$ préserve
les $\mathbf{A}^1$-équivalences faibles ; il induit donc un foncteur
$\R f_\star \colon \DMNeg[T,\Lambda]\to\DMNeg[S,\Lambda]$.

Le foncteur $\L f^\star\colon
\CompNeg(\PST[S,\Lambda])\to\CompNeg(\PST[T,\Lambda])$ préserve les
$\mathbf{A}^1$-équivalences faibles. Le couple de foncteurs $(\L f^\star,\R
f^\star)$ ainsi défini entre les catégories $\DMNeg[T,\Lambda]$ et
$\DMNeg[S,\Lambda]$ est un couple de foncteurs adjoints.
\end{proposition}

La démonstration est essentiellement la même que celle de
\cite[Proposition~2.8, p. 108]{morel-voevodsky}.
Cela vaut aussi pour l'énoncé
suivant, homologue de \cite[Proposition~2.9, p. 108]{morel-voevodsky} :

\begin{proposition}
Supposons que le morphisme $f\colon T \to S$ soit lisse. Le foncteur
$f^\star \colon \SmCor[S] \to\SmCor[T]$ admet alors un adjoint à
gauche $f_\sharp\colon \SmCor[T]\to \SmCor[S]$ (qui au
niveau des objets représentables
est tel que pour tout $X\in \Sm[T]$, $f_\sharp [X]=[X]$
où, à droite, $X$ est considéré comme un objet de $\Sm[S]$). Le
foncteur $f^\star\colon \PST[S,\Lambda]\to \PST[T,\Lambda]$ s'identifie
à la composition avec $f_\sharp\colon
\SmCor[T]\to \SmCor[S]$. Le couple de foncteurs
$(f^\star,f_\star)$ entre $\CompNeg(\PST[T,\Lambda])$ et
$\CompNeg(\PST[S,\Lambda])$ (resp. $\CompNeg(\ST[T,\Lambda])$ et
$\CompNeg(\ST[S,\Lambda])$) est une adjonction de Quillen (pour les
structures $\mathbf{A}^1$-localisées ou non) ;
on peut alors identifier $f^\star$ et $\L f^\star$.
Enfin, le foncteur $f^\star\colon \DMNeg[S,\Lambda]\to\DMNeg[T,\Lambda]$
admet un adjoint à gauche $\L f_\sharp$.
\end{proposition}

\subsection{Isomorphisme de Thom relatif}
\label{subsection-isomorphisme-thom-relatif}

Soit $k$ un corps parfait. Soit $\Lambda$ un anneau commutatif de coefficients.
Le foncteur d'inclusion de
$\CompNeg(\ST[k,\Lambda])$ dans
la catégorie des complexes bornés supérieurement (pour la numérotation
cohomologique) induit un foncteur pleinement fidèle
$\DMNeg[k,\Lambda]\to\DMmoins[k,\Lambda]$
(ceci découle aussitôt du fait que le foncteur de
$\mathbf{A}^1$-localisation $C_\star$ de
\cite[Proposition~3.2.3, Chapter~V]{livre-orange} préserve
$\CompNeg(\ST[k,\Lambda])$).

Pour $M\in\DMmoins[k,\Lambda]$ et $(p,q)\in\mathbf{Z}^2$, on note
$H^{p,q}(M)=\Hom_{\DMmoins[k,\Lambda]}(M,\mathbf{Z}(q)[p])$ où
$\mathbf{Z}(q)=0$ si $q<0$ et
$\mathbf{Z}(q)[2q]=\tM(\mathbf{A}^q/\mathbf{A}^q-\{0\})\in
\DMNeg[k,\Lambda]$.

Pour $\mathcal X\in \Hopt[k]$, on note aussi $\tH^{p,q}(\mathcal
X)=H^{p,q}(\tM(\mathcal X))$ et pour $\mathcal X\in\Ho[k]$,
$H^{p,q}(\mathcal X)=\tH^{p,q}(\mathcal X_+)$.

\begin{definition}\label{definition-classe-thom}
Soit $S\in\Sm[k]$. Soit $L$ un fibré en droites sur $S$. On note $c_1(L)$
la classe de $L$ dans $H^{2,1}(S)\simeq \Pic(S)\otimes \Lambda$ (cf.
\cite[Corollary~3.4.3, Chapter~V]{livre-orange}).
Soit $E$ un fibré vectoriel de rang $r$, on note $c_1(E),\dots,c_r(E)$ les
uniques éléments de $H^{2\star,\star}(S)$ tels que
\[\xi^r+c_1(E)\xi^{r-1}+\dots+c_r(E)\]
appartienne au noyau de $H^{2r,r}(\mathbf{P}(E\oplus
\OO_X))\to H^{2r,r}(\mathbf{P}(E))$ où $\xi=c_1(\OO(1))$\;\footnote{Ici,
$\mathbf{P}(E)$ est le Proj de l'Algèbre symétrique du dual de $E$.}.

On note $t_E\in \tH^{2r,r}(\Th_S E)\simeq
\ker(H^{2r,r}(\mathbf{P}(E\oplus
\OO_X))\to H^{2r,r}(\mathbf{P}(E)))$ (voir
\cite[Proposition~3.5.1, Chapter~V]{livre-orange}) la classe correspondant à
$\xi^r+c_1(E)\xi^{r-1}+\dots+c_r(E)$ \emph{via} l'isomorphisme canonique
$\Th_S E=E/(E-\{0\})\simeq \mathbf{P}(E\oplus \OO_X)/\mathbf{P}(E)$ dans
$\Hopt[k]$ (cf.~\cite[Proposition~2.17, p.~112]{morel-voevodsky}).
\end{definition}

\begin{proposition}\label{proposition-isomorphisme-thom-relatif}
Soit $S\in\Sm[k]$. Soit $E$ un fibré vectoriel de rang $r$ sur $S$. Alors,
le morphisme
\[\tM(\Th_S E)\to \tM(\Th_S \mathbf{A}^r)\]
dans $\DMNeg[S,\Lambda]$ induit par la classe $t_E\in H^{2r,r}(\Th_S E)$ est un
isomorphisme dans $\DMNeg[S,\Lambda]$. En particulier, $t_E$ induit un
isomorphisme $K\tM(\Th_S E)\isomto  K\tM(\Th_S \mathbf{A}^r)$ dans $\Hopt$.
\end{proposition}

On peut considérer $\Th_S E$ comme un objet aussi bien de $\Hopt[S]$ que de
$\Hopt[k]$, mais on a bien sûr $a_\sharp \Th_S E=\Th_S E$ et
$a_\sharp\tM(\Th_S E)=\tM(a_\sharp \Th_S E)$. On a des bijections :
\begin{eqnarray*}
H^{2r,r}(\Th_S E)&=&\Hom_{\DMNeg[k,\Lambda]}(a_\sharp \tM(\Th_S
E),\tM(\Th_k \mathbf{A}^r))\\&=&\Hom_{\DMNeg}(\tM(\Th_S E),a^\star
\tM(\Th_k\mathbf{A}^r)\\
&=&\Hom_{\DMNeg}(\tM(\Th_S E),\tM(\Th_S\mathbf{A}^r))
\end{eqnarray*}
Ainsi, à $t_E$ correspond bien un morphisme
$\tM(\Th_S E)\to\tM(\Th_S\mathbf{A}^r)$ dans $\DMNeg$.
Pour montrer que c'est un
isomorphisme, il suffit de le vérifier après application
des foncteurs de restriction évidents
$\DMNeg[S,\Lambda]\to\DMNeg[U,\Lambda]$ où $U$ parcourt un recouvrement
de $S$ par des ouverts au-dessus desquels $E$ serait trivial. On peut ainsi
supposer que $E=\mathbf{A}^r_S$. Le fibré $E$ provient alors par image
inverse d'un fibré (trivial) sur $\Spec k$. On se ramène ainsi au cas où
$S=\Spec k$ et cela résulte alors de
\cite[Proposition~3.5.1, Chapter~V]{livre-orange}.

\begin{remarque}
Dans le cas où le fibré $E$ est $\mathbf{A}^r_S$, l'isomorphisme
de la proposition~\ref{proposition-isomorphisme-thom-relatif} est bien
évidemment l'identité, pourvu qu'on ait fait le bon choix de signe dans
l'isomorphisme $H^{2,1}(S,\mathbf{Z})\simeq \Pic(S)$.
\end{remarque}

\subsection{Classes tautologiques}

On suppose ici que $k$ est un corps parfait et que $S\in\Sm[k]$.

\begin{definition}\label{definition-eilenberg-maclane}
Pour tout fibré vectoriel $E$ sur $S$, on note $K_E\in
\Delta^\opp\Sm^\opp\Enspt$ le préfaisceau $K\tM(\Th_S E)$. En particulier,
pour tout entier naturel $r$, on note $K_r=K_{\mathbf{A}^r_S}$.
\end{definition}

D'après l'isomorphisme de Thom relatif (cf.
proposition~\ref{proposition-isomorphisme-thom-relatif}), on a un
isomorphisme canonique $K_E\simeq K_r$ dans $\Hopt$ si $E$ est un fibré
vectoriel de rang $r$ sur $S$.

\begin{proposition}\label{proposition-classes-tautologiques}
Soit $E$ un fibré vectoriel de rang $r$ sur $S$. On note $\tau_E$ le
morphisme tautologique (d'adjonction) $\tau_E\colon \Th_S E\to
K_E=K\tM(\Th_S E)$ dans
$\Sm^\opp\Enspt$. Alors, le diagramme suivant est commutatif dans $\Hopt$ :
\[\xymatrix{\ar[rd]_-{t_E}\Th_S E\ar[r]^-{\tau_E} & K_E\ar[d]^\sim \\
& K_r}\]
\end{proposition}

En vertu de l'adjonction de la sous-section~\ref{subsection-adjonction},
cette compatibilité est tautologique puisque les
deux morphismes $t_E\colon \Th_S E\to K_r$ et $K_E\simeq K_r$ sont d'une
manière ou d'une autre induits par la classe de Thom $t_E$.

\begin{proposition}
On suppose donnés deux fibrés vectoriels $E$ et $F$ de rangs respectifs $r$
et $s$ sur $S$. Alors, le diagramme suivant est commutatif dans $\Hopt$,
où les morphismes horizontaux sont les morphismes \guil{produits} évidents
et les morphismes verticaux sont induits par l'isomorphisme de Thom relatif
appliqué aux fibrés vectoriels $E$, $F$ et $E\oplus F$ :
\[\xymatrix{K_E\wedge K_F \ar[r]\ar[d] & K_{E\oplus F}\ar[d] \\
K_r\wedge K_s\ar[r] & K_{r+s}}\]
\end{proposition}

En utilisant l'adjonction de la sous-section~\ref{subsection-adjonction},
cela se déduit du diagramme commutatif suivant
qui énonce une compatibilité facile
des classes de Thom à la somme directe des fibrés :
\[\xymatrix{
\ar[d]^{t_E\otimes t_F}
\tM(\Th_S E)\otimes \tM(\Th_S F)\ar[r]^-\sim & \ar[d]^{t_{E\oplus F}}
\tM(\Th_S(E\oplus F))\\
\tM(\Th_S \mathbf{A}^r_S)\otimes \tM(\Th_S \mathbf{A}^r_S)\ar[r]^-\sim &
\tM(\Th_S \mathbf{A}^{r+s}_S)
}\]

\begin{remarque}
On a utilisé ici le produit tensoriel sur la catégorie $\DMNeg[S,\Lambda]$.
S'il n'en sera pas fait ici un usage intensif,
il convient cependant d'en esquisser une construction.
Ce produit tensoriel
provient d'une structure monoïdale symétrique sur $\SmCor$ qui est telle
que $[X]\otimes[Y]=[X\times_S Y]$. Cette structure monoïdale s'étend sans
difficulté à la sous-catégorie pleine $\mathcal S$ de $\CompNeg\PST[S,\Lambda]$
formée des complexes de préfaisceaux avec transferts qui terme à terme sont
sommes directes de préfaisceaux représentés par des objets de $\SmCor$.
Comme à équivalences faibles près, on peut remplacer (fonctoriellement)
tout objet de $\CompNeg\PST[S,\Lambda]$ par un objet de $\mathcal S$, il suffit
de justifier que $\otimes\colon \mathcal S\times \mathcal S\to\mathcal S$
préserve les équivalences faibles. Pour cela, on se ramène à montrer que
pour tout $[U]\in\SmCor$, le foncteur $-\otimes[U]\colon \mathcal S\to
\mathcal S$ préserve les ($\mathbf{A}^1$-)équivalences faibles. Ceci
s'obtient facilement en considérant le foncteur $\SheafHom([U],-)\colon
\CompNeg\PST[S,\Lambda]\to\CompNeg\PST[S,\Lambda]$ qui à $K$ associe
$[X]\longmapsto K([U]\otimes[X])$. On peut en effet dériver à droite ce
foncteur pour obtenir $\R\SheafHom([U],-)\colon
\DMNeg[S,\Lambda]\to\DMNeg[S,\Lambda]$, l'adjoint à gauche de ce foncteur
étant donné sur les objets de $\mathcal S$ par le foncteur $-\otimes [U]$,
lequel préserve donc les ($\mathbf{A}^1$-)équivalences faibles.
\end{remarque}

\section{Construction de l'opération totale}
\label{section-operation-totale}

On fixe un corps de base parfait $k$. Soit $\Lambda$ un anneau de
coefficients. On se donne une action (à gauche) d'un groupe fini $G$ sur un
ensemble fini (non vide) $A$ à $n$ éléments. On
fixe aussi $S\in\Sm[k]$ et $\pi\colon U\to S$ un $G$-torseur étale (à
gauche) au-dessus de $S$ (ainsi, $S=G\backslash U$). (Le seul cas
véritablement intéressant est celui où l'action du groupe $G$ sur $A$ est
fidèle.)

\subsection{Construction de l'opération}
\label{subsection-construction-operation}

L'action de $G$ sur $A$ et le torseur $U$ permettent de définir un fibré
vectoriel de rang $n$ sur $S$ par descente fidèlement plate :

\begin{definition}\label{definition-xi}
On note $\xi$ le fibré vectoriel sur $S$ défini à isomorphisme unique près
par la donnée d'un isomorphisme de fibrés $\Phi\colon
\pi^\star \xi\isomto \OO_U^A$ tel que pour tout $g\in G$, le diagramme
évident de fibrés vectoriels sur $U$ commute :
\[\xymatrix{\ar[d]^{\Phi}\pi^\star\xi\ar@{=}[r] & (\pi\circ g)^\star \xi \ar@{-}[r]^-\sim&
g^\star \pi^\star \xi\ar[d]^{g^\star\Phi}\\
\OO_U^A &\ar[l]_-{\circ g} \OO_U^A&\ar[l]_-\sim
g^\star \OO_U^A}\]
où $\circ g$ est l'automorphisme (de permutation)
de $\OO_U^A$ induit par la bijection $A\to A$ correspondant à l'action de $g$.
\end{definition}

Comme les matrices de permutation sont des matrices orthogonales, on peut
remarquer qu'il existe un isomorphisme canonique de fibrés vectoriels
$\xi\simeq \xi^\vee$ où $\xi^\vee$ est le dual de $\xi$.

\begin{remarque}
Plus généralement, si $V$ est une représentation $k$-linéaire de dimension
finie de $G$, on peut remplacer $\OO_U^A$ par $\Hom_k(V,\OO_U)$ dans la
définition ci-dessus pour associer à une telle représentation $V$ un fibré
vectoriel sur $S$. La définition ci-dessus correspond alors
au cas particulier des représentations de permutation.
\end{remarque}

\begin{definition}\label{definition-operation-totale}
Pour tout fibré vectoriel $E$ sur $S$, on définit un morphisme 
\[P_E\colon K_E \to K_{E\otimes\xi}\]
dans $\Delta^\opp\Sm^\opp\Enspt$. Il est induit par le morphisme de
préfaisceaux (aussi noté $P_E$) $\Ltr(E)\to\Ltr(E\otimes \xi)$ sur
$\Sm^\opp\Enspt$ faisant commuter le diagramme suivant pour tout $X\in\Sm$
:
\[\xymatrix{\cequi(X\times_S E/X,0)\ar[r]^{P_E}\ar[dd]^{x\longmapsto
x^n}_{\cercle{1}} &
\cequi(X\times_S(E\otimes\xi)/X,0)\ar[d]^{\pi^\star}_{\cercle{4}\sim} \\
& \cequi(X\times_S(E^A\times_S U)/(X\times_S U),0)^G
\ar@{^(->}[d]_{\cercle{3}}\\
\cequi(X\times_S E^A/X,0)\ar[r]_-{\cercle{2}} &\cequi(X\times_S(E^A\times_S
U)/(X\times_S U),0)
}\]
\end{definition}

Le morphisme \cercle{1} est celui d'\guil{élévation à la puissance $n$}
(qui induit un morphisme $K_E\to K_{E^A}$) : il est induit par la structure
tensorielle sur $\SmCor$ et la diagonale de $E$. Le morphisme \cercle{2}
provient du changement de base (étale) par $U\to S$, tout comme le
morphisme \cercle{4}. Il convient de préciser ce fait. Tout d'abord,
l'action de $G$ sur $U$ induit une action sur $(E\otimes\xi)\times_S U$.
Par construction, on a un isomorphisme canonique $(E\otimes\xi)\times_S
U\simeq E^A\times_S U$ et cet isomorphisme est $G$-équivariant (où
$E^A\times_S U$ est muni du produit de l'action par permutation sur $E^A$
et de celle donnée sur $U$). Ceci définit les morphismes \cercle{3} et
\cercle{4}, \cercle{4} étant un isomorphisme par descente étale des cycles.
L'existence de $P_E\colon \Ltr(E)\to\Ltr(E\otimes\xi)$
provient du fait évident que le composé de \cercle{1} et
\cercle{2} se factorise par \cercle{3}.

Si deux éléments $x$ et $x'$ de $\cequi(X\times_S E/X,0)$ diffèrent d'un
élément $y=x'-x\in \cequi(X\times_S (E-\{0\})/X,0)\subset \cequi(X\times_S
E/X,0)$, alors la différence $(x+y)^n-x^n$ de leurs images par \cercle{1}
est une somme d'orbites sous l'action du groupe symétrique de termes
$y^px^{n-p}$ avec $p\geq 1$, ce qui montre que cette différence appartient
à $\cequi(X\times_S (E^A-\{0\})/X,0)\subset\cequi(X\times_S (E^A)/X,0)$. La
considération des variantes des morphismes de groupes \cercle{2},
\cercle{3} et \cercle{4} où l'on aurait épointé les fibrés vectoriels $E^A$
et $E\otimes\xi$ montre que $P_E(x')-P_E(x)\in \cequi(X\times_S
(E\otimes\xi-\{0\})/X,0)\subset \cequi(X\times_S (E\otimes\xi)/X,0)$.
Ainsi, le morphisme $P_E\colon \Ltr(E)\to\Ltr(E\otimes \xi)$ passe bien au
quotient pour définir un morphisme $P_E\colon K_E\to K_{E\otimes\xi}$.

Les deux propositions suivantes sont évidentes :

\begin{proposition}\label{proposition-compatibilite-produit}
Si $E$ et $F$ sont des fibrés vectoriels sur $S$, alors on dispose d'un
diagramme commutatif (où les morphismes verticaux sont les morphismes
\guil{produits} évidents) :
\[\xymatrix{\ar[d]K_E\wedge K_F\ar[rr]^-{P_E\wedge P_F} & &
K_{E\otimes\xi}\wedge K_{F\otimes\xi}\ar[d] & \\
K_{E\oplus F} \ar[r]^-{P_{E\otimes F}} & K_{(E\oplus F)\otimes \xi}
\ar[r]^\sim & K_{E\otimes\xi\oplus F\otimes \xi}}\]
\end{proposition}

\begin{proposition}\label{proposition-compatibilite-composition}
Si outre $(G,A,U)$ on se donne un groupe fini $G'$, un $G'$-ensemble fini
$A'$ et un $G'$-torseur étale $U'$ au-dessus de $S$ et que l'on note $\xi'$
le fibré vectoriel sur $S$ associé à $(G',A',U')$ comme $\xi$ l'était à
$(G,A,U)$, on a un diagramme commutatif :
\[\xymatrix@C=2cm{
& K_{E\otimes \xi} \ar[dr]^{~~~P_{E\otimes\xi,(G',A',U')}} \\
K_E
\ar[ru]^{P_{E,(G,A,U)~~~}}
\ar[rr]_-{P_{E,(G\times G',A\times A',U\times_S U')}} & &
K_{E\otimes \xi\otimes \xi'}
}\]
\end{proposition}

\begin{definition}
Pour tout fibré vectoriel $E$ de rang $r$ sur $S$,
on a défini un morphisme $P_E\colon
K_E \to K_{E\otimes \xi}$ dans $\Delta^\opp\Sm^\opp\Enspt$. En utilisant
l'isomorphisme de Thom relatif pour $E$ et $E\otimes\xi$
(cf.~\S\ref{subsection-isomorphisme-thom-relatif}), on en déduit un
morphisme dans $\Hopt$, noté aussi $P_E$ :
\[P_E\colon K_r\simeq K_E \vers {P_E} K_{E\otimes \xi} \simeq
K_{rn}\;\text{.}\]
Ainsi, pour tout $\mathcal X\in \Hopt$, $P_E$ définit une application :
\[P_E\colon \tH^{2r,r}(\mathcal X)\to \tH^{2rn,rn}(\mathcal X)\]
\end{definition}

Nous verrons au \S{}\ref{subsection-independante-fibre} que cette action
$P_E$ ne dépend que du rang du fibré $E$ (et bien sûr aussi de $(G,A,U)$).
Notons cependant le calcul suivant :

\begin{proposition}\label{proposition-compatibilite-triviale}
Soit $E$ un fibré vectoriel de rang $r$ sur $S$. On note $t_E\in
\tH^{2r,r}(\Th_S E)$ la classe de Thom, on a alors :
\[P_E(t_E)=\delta^\star t_{E\otimes \xi}\in \tH^{2rn,rn}(\Th_S
(E\otimes \xi))\]
où $\delta\colon \Th_S E\to \Th_S E\otimes \xi$ est le morphisme induit
l'inclusion $E\to E\otimes \xi$ déduite de l'inclusion évidente
$\mathbf{A}^1_S\to \xi$.
\end{proposition}

En effet, d'après la proposition~\ref{proposition-classes-tautologiques},
il suffit de vérifier que, strictement, au niveau des cycles, on a
$P_E(\tau_E)=\delta^\star (\tau_{E\otimes \xi})$, ce qui est immédiat.

\subsection{Raffinement de l'opération totale}
\label{subsection-operation-p-tilde}

\begin{definition}
On note $\Sc$ la sous-catégorie pleine de la catégorie des faisceaux
pointés sur $\Sm$ formée des faisceaux pointés isomorphes à une somme directe
d'objets de la forme $X_+$ pour $X\in\Smqp$ (où $\Smqp$ est la
sous-catégorie pleine de $\Sm$ formée des $S$-schémas quasi-projectifs).
\end{definition}

\begin{definition}
Soit $\mathcal X\in \Delta^\opp\Sm^\opp\Enspt$. On note $\Triv$ la catégorie
des morphismes $\mathcal X'\to\mathcal X$ dans $\Delta^\opp\Sm^\opp\Enspt$
qui sont des fibrations triviales locales et
tels que $\mathcal X'\in\Delta^\opp\Sc$.
\end{definition}

On encourage le lecteur à lire \cite[\S{1}]{jardine} où est définie la
notion de fibration triviale locale. On rappelle que dans la situation
considérée ici une fibration (triviale) locale est un morphisme de
préfaisceaux simpliciaux induisant des fibrations (triviales) d'ensembles
simpliciaux après application d'un ensemble conservatif de foncteurs fibres
pour le site $\Sm_\Nis$.

\begin{proposition}
Soit $\mathcal Y\in\Delta^\opp\Sm^\opp\Enspt$ un préfaisceau tel que pour
tout $X\in\Sm$, $\mathcal Y(X)$ soit un ensemble simplicial fibrant. Alors,
pour tout $\mathcal X\in\Delta^\opp\Sm^\opp\Enspt$ :
\[\Hom_{\Hopt}(\mathcal X,\mathcal Y)\simeq \underset{\mathcal X'
\in \Triv^\opp}{\colim} \pi_0\hom(\mathcal X',\mathcal Y)\]
\end{proposition}

Ceci résulte de \cite[p.~55]{jardine}.
On pourra noter que le système inductif se factorise par la catégorie
homotopique $\pi\Triv$ de $\Triv$. Il n'y a évidemment pas de différence
entre la colimite indexée par $\Triv^\opp$ et par $\pi\Triv^\opp$.
L'avantage est que la catégorie $\pi\Triv^\opp$ est filtrante. Cependant, ni
$\Triv$ ni $\pi\Triv$ ne sont des petites catégories, mais dans
\cite[Lemma~1.12, p.~51]{morel-voevodsky}, il est montré que $\pi\Triv^\opp$
contient une petite sous-catégorie pleine cofinale.

\begin{corollaire}\label{corollaire-calcul-hom-hopt}
Pour tout fibré vectoriel $E$ sur $S$ et $\mathcal
X\in\Delta^\opp\Sm^\opp\Enspt$, on a une bijection canonique :
\[\Hom_{\Hopt}(\mathcal X,K_E)\simeq \underset{\mathcal
X'\in\Triv}{\colim}\pi_0\hom(\mathcal X',\Sing K_E)\]
\end{corollaire}

Tout d'abord, les sections de $\Sing K_E$ sur les objets de $\Sm$ sont des
ensembles simpliciaux induits par des groupes abéliens simpliciaux ; ils
sont donc fibrants. En outre,
$\Sing K_E$ est un objet $\mathbf{A}^1$-local de $\Hopt$ (on peut supposer
que $E$ est trivial puisque la vérification est locale sur $S$, mais alors
$K_E$ est induit par un préfaisceau simplicial sur $\Sm[k]$ qui est
$\mathbf{A}^1$-local d'après \cite[\S{}3.2, Chapter~V]{livre-orange}).

\begin{proposition}\label{proposition-quotient}
On considère le foncteur $\Smqp\to\Smqp$ qui à $X$ associe
$G\backslash(X^A\times_S U)$. On en déduit un foncteur $\Sc\to\Sc$ noté
$\mathcal X\longmapsto G\backslash (\mathcal X^{\wedge A}\times_S U)$ qui à
$X_+$ associe $G\backslash(X^A\times_S U)_+$ et qui commute aux limites
inductives filtrantes. Le foncteur $\Delta^\opp\Sc\to\Delta^\opp\Sc$
que l'on en déduit préserve les équivalences faibles (locales pour la
topologie de Nisnevich) et les $\mathbf{A}^1$-équivalences faibles.
\end{proposition}

Ceci résulte des résultats principaux de \cite[\S5.1,
5.2]{voevodsky-deligne} concernant les foncteurs $X\longmapsto X^A$ et
$X\longmapsto G\backslash X$. On notera que les catégories homotopiques
$\Hopt$ et $\Hosimppt$ définies à partir de la catégorie $\Sm$ (ou $\Smqp$,
cela revient au même) sont des sous-catégories pleines des catégories
homotopiques homologues définies en considérant la catégorie des
$S$-schémas de type fini (quasi-projectifs) utilisée dans
\cite{voevodsky-deligne} ; cela s'obtient par exemple en utilisant des
arguments semblables à ceux de \cite[\S{}1.3]{voevodsky-eilenberg-maclane}.

\begin{definition}
La catégorie $\Hopt$ étant équivalente à la catégorie obtenue en inversant
formellement les $\mathbf{A}^1$-équivalences faibles dans $\Delta^\opp\Sc$,
la proposition~\ref{proposition-quotient} permet de définir un foncteur
$\Hopt\to \Hopt$ déduit du foncteur précédemment défini sur
$\Delta^\opp\Sc$ et que nous noterons $\mathcal X\longmapsto
G\backslash_{\L}(\mathcal X^{\wedge A}\times_S U)$.
\end{definition}

\bigskip

Comme dans le \S\ref{subsection-construction-operation}, $E$ sera un fibré
vectoriel sur $S$.

\begin{definition}
Soit $X\in\Smqp$. On définit une application
\[\tilde{P}_E\colon K_E(X)\to K_{E\otimes \xi}(G\backslash(X^A\times_S U))\] de
façon à faire commuter le diagramme suivant :
\[\xymatrix{
\cequi(X\times_S E/X,0)\ar[dd]^{x\longmapsto x^{\otimes
n}}\ar[r]^-{\tilde{P}_E} & 
\cequi(G\backslash(X^A\times_S U)\times_S (E\otimes \xi)/G\backslash(X^A\times U),0)\ar[d]^{\sim}\\
& \cequi(G\backslash(X^A\times_S U\times_S
E^A)/G\backslash(X^A\times U),0)\ar@{^(->}[d]\\
\cequi(X^A\times_S E^A/X^A,0)\ar[r]& \cequi(X^A\times_S U\times_S
E^A/X^A\times U,0)
}\]
\end{definition}

Le principe est le même que dans la
définition~\ref{definition-operation-totale} à ceci près que l'on n'a pas
utilisé la diagonale $E\to E^A$. 
Cette construction raffine donc
la première qui se déduit de celle-ci en prenant une image inverse par la
diagonale $X=G\backslash(X\times_S U)\to G\backslash(X^A\times_S U)$.

\medskip

Par passage à la limite inductive, la définition de
$\tilde{P}_E$ s'étend en une application
\[\tilde{P}_E\colon \Hom(\mathcal X,K_E)\to
\Hom(G\backslash(\mathcal X^{\wedge A}\times_S U),K_{E\otimes \xi})\]
pour tout $\mathcal X\in\Sc$.

Si $D\in\Smqp$, on déduit de cette construction une application
\[\tilde{P}_E\colon \Hom(\mathcal X,\SheafHom(D,K_E))\to
\Hom(G\backslash(\mathcal X^{\wedge A}\times_S U),\SheafHom(D,
K_{E\otimes \xi}))\]
pour tout $\mathcal X\in\Sc$. On peut en effet appliquer la construction
précédente à $D_+\wedge \mathcal X$ et utiliser le morphisme déduit de
façon évidente de la diagonale de $D$ :
\[D_+\wedge (G\backslash (\mathcal X^{\wedge A}\times_S U))
\to G\backslash (D^A_+\wedge \mathcal X^{\wedge A}\times_S U)\]

En faisant parcourir à $D$
les espaces affines de l'objet cosimplicial standard
$\mathbf{\Delta}^\opp$, on obtient un morphisme :
\[\tilde{P}_E\colon \Hom(\mathcal X,\Sing K_E)\to
\Hom(G\backslash(\mathcal X^{\wedge A}\times_S U),\Sing K_{E\otimes \xi})\]
pour tout $\mathcal X\in\Delta^\opp\Sc$.

\begin{definition}\label{definition-operation-totale-raffinee}
En appliquant cette construction à $\Delta^n_+\wedge \mathcal X$ pour tout
$n\in \mathbf{N}$ et en utilisant la diagonale de $\Delta^n$ (comme plus
haut pour $D$), on peut observer que l'application précédente est
l'application induite au niveau des $0$-simplexes par un morphisme
d'ensembles simpliciaux, aussi noté $\tilde{P}_E$,
pour tout $\mathcal X\in\Delta^\opp \Sc$ :
\[\tilde{P}_E\colon \hom(\mathcal X,\Sing K_E)\to
\hom(G\backslash(\mathcal X^{\wedge A}\times_S U),\Sing K_{E\otimes
\xi})\;\text{.}\]
\end{definition}

\begin{definition}
Soit $\mathcal X\in \Delta^\opp\Sc$. Pour tout $\mathcal X'\in \Triv$, la
construction précédente de $\tilde{P}_E$ et la
proposition~\ref{proposition-quotient}
permettent de définir la flèche de gauche sur le diagramme suivant :
\[\xymatrix{\pi_0\hom(\mathcal X',\Sing K_E)\ar[r]^-{\tilde{P}_E}\ar[d] &
\pi_0 \hom(G\backslash(\mathcal X'^{\wedge A\times_S U}),\Sing
K_{E\otimes\xi}) \ar[d]\\
\Hom_{\Hopt}(G\backslash(\mathcal X^{\wedge A}\times_S
U),K_{E\otimes\xi})\ar[r]^-\sim
& \Hom_{\Hopt}(G\backslash(\mathcal X'^{\wedge A}\times_S
U),K_{E\otimes\xi}) \\
}\]
En passant à la colimite sur $\mathcal X'$, on obtient d'après le
corollaire~\ref{corollaire-calcul-hom-hopt} une application
\[\tilde{P}_E\colon \Hom_{\Hopt}(\mathcal X,K_E)\to
\Hom_{\Hopt}(G\backslash(\mathcal X^{\wedge A}\times_S
U),K_{E\otimes\xi})\;\text{.}\]
Plus généralement, pour $\mathcal X\in\Hopt$, on en déduit une application
fonctorielle :
\[\tilde{P}_E\colon \Hom_{\Hopt}(\mathcal X,K_E)\to
\Hom_{\Hopt}(G\backslash_{\L}(\mathcal X^{\wedge A}\times_S
U),K_{E\otimes\xi})\;\text{.}\]

\end{definition}

\begin{proposition}\label{proposition-raffinement}
La construction $\tilde{P}_E$ raffine la construction $P_E$, c'est-à-dire
que le diagramme suivant est commutatif pour tout $\mathcal X\in\Hopt$ :
\[\xymatrix{\Hom_{\Hopt}(\mathcal X,K_E)\ar[r]^-{\tilde{P}_E}\ar[rd]_-{P_E}&
\Hom_{\Hopt}(G\backslash_{\L}(\mathcal X^{\wedge A}\times_S
U),K_{E\otimes\xi})\ar[d] \\
& \Hom_{\Hopt}(\mathcal X,K_{E\otimes\xi})}\]
où le morphisme de droite, pour $\mathcal X\in\Delta^\opp\Sc$, est induit
par la diagonale de $\mathcal X$.
\end{proposition}

\subsection{Compatibilité au changement de base}
\label{subsection-compatibilite-cb}

Dans cette section, on suppose que l'on s'est donné un morphisme de type
fini $f\colon T\to S$ entre schémas réguliers.

\begin{definition}\label{definition-cb-k-e}
Soit $F$ un fibré vectoriel sur $S$. On note $f^\star\colon K_F\to f_\star
K_{f^\star F}$ le morphisme dans $\Sm^\opp\Enspt$ induit par la
construction de la définition~\ref{definition-cb-cequi}. On note aussi
$f^\star\colon f^\star K_F\to K_{f^\star F}$ le morphisme adjoint à ce
morphisme.
\end{definition}

\begin{proposition}\label{proposition-compatibilite-maclane-cb}
S'il existe un schéma régulier $B$ tel que $f\colon T\to S$
soit un morphisme dans $\Sm[B]$, alors on a des isomorphismes :
\[\L f^\star K_F\isomto f^\star K_F\underset{f^\star}\isomto K_{f^\star F}\]
dans $\Hopt[T]$.
\end{proposition}

Ces isomorphismes vallent en fait déjà dans la catégorie homotopique
$\Hosimppt[T]$ (avant $\mathbf{A}^1$-localisation). Le résultat est évident
si $f$ est un morphisme lisse puisqu'alors $f^\star$ préserve les
équivalences faibles (voir \cite[Corollary~1.24, p.~104]{morel-voevodsky})
et on a tautologiquement $f^\star K_F\simeq K_{f^\star F}$ dans
$\Sm[T]^\opp\Enspt$.

Si le fibré $F$ provient par image inverse d'un fibré vectoriel sur la base
$B$, en utilisant la compatibilité à la composition des foncteurs $f^\star$ et
$\L f^\star$ et le fait que l'énoncé de la proposition vaille pour $S\to B$
et $T\to B$, on peut l'obtenir pour $T\to S$.

Dans le cas général, il suffit de s'assurer que les morphismes $\L f^\star
K_F\to f^\star K_F$ et $f^\star K_F\to K_{f^\star F}$ deviennent des
équivalences faibles après restriction à des ouverts d'un recouvrement de
$T$. Si $U$ est un ouvert de $S$ tel que $F_{|U}$ est trivial, les
arguments précédents montrent que le résultat vaut après restriction à
l'ouvert $f^{-1}(U)$ de $T$. Ceci fournit un recouvrement ouvert
convenable de $T$.

\begin{remarque}
Sous les hypothèses de la proposition précédente (avec $B=\Spec k$ et $k$
un corps parfait), si $F$ est un fibré de rang $r$, on a bien sûr un
diagramme commutatif 
dans $\Hopt[T]$ énonçant la compatibilité au
changement de base de l'isomorphisme de Thom :
\[\xymatrix{\L f^\star K_F\ar[r]^\sim\ar[d]^\sim & \L f^\star
K_r\ar[d]^\sim \\
K_{f^\star F}\ar[r]^\sim & K_r }\]
(Les morphismes verticaux sont ceux de la
définition~\ref{definition-cb-k-e} et les morphismes horizontaux
proviennent de l'isomorphisme de Thom relatif appliqué aux fibrés $F$ et
$f^\star F$.)
\end{remarque}

\begin{proposition}\label{proposition-compatibilite-p-tilde-cb}
Soit $G$ un groupe fini. Soit $A$ un $G$-ensemble fini de cardinal $n$.
Soit $U$ un $G$-torseur étale sur $S$. Soit $E$ un fibré vectoriel
sur $S$. Notons $U_T$ et $f^\star E$ les données homologues définies
sur $T$ par changement de base par $f$. Notons $\tilde{P}_E$ et
$\tilde{P}_{f^\star E}$ les opérations associées à ces données dans la
sous-section~\ref{subsection-operation-p-tilde}. Alors, le diagramme
suivant est commutatif pour tout $\mathcal X\in\Hopt$ :
\[\xymatrix{
\Hom_{\Hopt}(\mathcal X,K_E)\ar[r]^-{\tilde{P}_E}\ar[d]^{f^\star} &
\Hom_{\Hopt}(G\backslash_\L(\mathcal X^{\wedge A}\times_S
U),K_{E\otimes\xi})\ar[d]^{f^\star}\\
\Hom_{\Hopt[T]}(\L f^\star \mathcal X,K_{f^\star E})\ar[r]^-{\tilde{P}_{f^\star
E}} &
\Hom_{\Hopt[T]}(G\backslash_\L((\L f^\star \mathcal X)^{\wedge A}\times_T
U_T),
K_{f^\star E\otimes f^\star \xi})
}\]
\end{proposition}

Pour établir cette compatibilité, il suffit d'obtenir la commutativité du
diagramme semblable où l'on aurait remplacé le $\Hom$ dans les catégories
$\Hopt$ et $\Hopt[T]$ par ceux dans les catégories de préfaisceaux
simpliciaux pointés sur $\Sm$ et $\Sm[T]$ et où $\mathcal X$ serait un objet
de $\Delta^\opp\Sc$. En raisonnant degré par degré et par passage à la
colimite, on se ramène à montrer la commutativité du diagramme suivant pour
tout $X\in\Sm$ :
\[\xymatrix{K_E(X)\ar[r]^-{\tilde{P}_E} \ar[d]^{f^\star} &
K_{E\otimes\xi}(G\backslash(X^A\times_S U))\ar[d]^{f^\star} \\
K_{f^\star E}(X_T)\ar[r]^-{\tilde{P}_{f^\star E}}& K_{f^\star E\otimes
f^\star \xi}(G\backslash(X_T^A\times_T U_T))}\]
Cette compatibilité est triviale.

\begin{remarque}
Dans le cas particulier où $f\colon T\to S$ est un morphisme dans $\Sm[k]$
avec $k$ un corps parfait, on peut récrire la compatibilité de la
proposition~\ref{proposition-compatibilite-p-tilde-cb} sous la forme
du carré commutatif suivant (où l'on a noté $r$ le rang de $E$),
pour tout $\mathcal X\in \Hopt$ :
\[\xymatrix{
\ar[d]^{f^\star}\tH^{2r,r}(\mathcal X)\ar[r]^-{\tilde{P}_E} &
\tH^{2rn,rn}(G\backslash_\L(\mathcal X^{\wedge A}\times_S
U))\ar[d]^{f^\star} \\
\tH^{2r,r}(\L f^\star \mathcal X)\ar[r]^-{\tilde{P}_{f^\star E}} &
\tH^{2rn,rn}(G\backslash_\L((\L f^\star \mathcal X)^{\wedge A}\times_T U_T))
}\]
\end{remarque}

\begin{corollaire}\label{corollaire-compatibilite-cb}
Avec les hypothèses et notations de la
proposition~\ref{proposition-compatibilite-p-tilde-cb}, on a aussi le carré
commutatif suivant :
\[\xymatrix{
\Hom_{\Hopt}(\mathcal X,K_E)\ar[r]^-{P_E}\ar[d]^{f^\star} &
\Hom_{\Hopt}(\mathcal X,K_{E\otimes\xi})\ar[d]^{f^\star}\\
\Hom_{\Hopt[T]}(\L f^\star \mathcal X,K_{f^\star E})\ar[r]^-{P_{f^\star
E}} &
\Hom_{\Hopt[T]}(\L f^\star \mathcal X,K_{f^\star E\otimes f^\star \xi})
}\]
\end{corollaire}

\begin{remarque}
Appliqué avec $\mathcal X=K_E$ et modulo la
proposition~\ref{proposition-compatibilite-maclane-cb},
ce corollaire énonce en particulier que le
morphisme $P_E\colon K_E\to K_{E\otimes \xi}$ est envoyé sur le morphisme
$P_{f^\star E}\colon K_{f^\star E}\to K_{f^\star E\otimes f^\otimes \xi}$
par le foncteur $\L f^\star\colon \Hopt[S]\to\Hopt[T]$.
\end{remarque}

\begin{corollaire}\label{corollaire-compatibilite-cb-f-lisse}
Avec les hypothèses et notations de la
proposition~\ref{proposition-compatibilite-p-tilde-cb}, et sous l'hypothèse
supplémentaire que $f$ est lisse, on a aussi le carré
commutatif suivant, pour tout $\mathcal Y\in \Hopt[T]$ :
\[\xymatrix{
\Hom_{\Hopt}(\L f_\sharp \mathcal Y,K_E)\ar[r]^-{P_E}\ar[d]^{\sim} &
\Hom_{\Hopt}(\L f_\sharp \mathcal Y,K_{E\otimes\xi})\ar[d]^{\sim}\\
\Hom_{\Hopt[T]}(\mathcal Y,K_{f^\star E})\ar[r]^-{P_{f^\star
E}} &
\Hom_{\Hopt[T]}(\mathcal Y,K_{f^\star E\otimes f^\star \xi})
}\]
\end{corollaire}

Il suffit d'appliquer le corollaire~\ref{corollaire-compatibilite-cb} avec
$\mathcal X=\L f_\sharp \mathcal Y$ et d'utiliser le morphisme d'adjonction
$\mathcal Y\to f^\star \L f_\sharp\mathcal Y$.

\subsection{Compatibilité aux classes de Thom}

\begin{proposition}
Soit $F$ un fibré vectoriel sur $S$. On a un isomorphisme canonique dans
$\Hopt$ :
\[G\backslash_{\L}(\Th_S F\times_S U)\simeq \Th_S
(F\otimes\xi)\;\text{.}\]
\end{proposition}

Ceci résulte de \cite[Example~5.2.8]{voevodsky-deligne}.

\medskip

Ainsi, pour deux fibrés vectoriels $E$ et $F$ de rang $r$ sur $S$,
la construction $\tilde{P}_E$ définit
une application :
\[\tilde{P}_E\colon
\tH^{2r,r}(\Th_S F)\to\tH^{2rn,rn}(\Th_S (F\otimes \xi))\]

\begin{proposition}\label{proposition-compatibilite-tilde-thom}
On suppose donnés deux fibrés vectoriels de rang $r$ sur $S$.
On note $t_E\in
\tH^{2r,r}(\Th_S E)$ la classe de Thom de $E$. Compte tenu de
l'identification ci-dessus, on a alors :
\[\tilde{P}_E(t_F)=t_{F\otimes\xi}\in H^{2rn,rn}(\Th_S (F\otimes
\xi))\;\text{.}\]
\end{proposition}

Dans le cas où $E=F$, la démonstration est semblable à celle de la
proposition~\ref{proposition-compatibilite-triviale}. Dans le cas général,
on peut observer que si $S$ est connexe (cas auquel on peut se ramener), le
groupe $\tH^{2rn,rn}(\Th_S (F\otimes \xi))$ dans lequel s'énonce l'égalité
est canoniquement isomorphe à l'anneau de coefficients $\Lambda$, le
générateur canonique étant donné par $t_{F\otimes\xi}$. Pour vérifier
l'égalité, il est donc permis de remplacer $S$ par un ouvert dense
au-dessus duquel les deux fibrés $E$ et $F$ sont isomorphes (par exemple
parce qu'ils seraient tous les deux triviaux).

\begin{corollaire}\label{corollaire-compatibilite-thom}
On suppose donnés deux fibrés vectoriels $E$ et $F$ de rang $r$ sur $S$.
On a alors :
\[P_E(t_F)=\delta^\star t_{F\otimes \xi}\]
dans $\tH^{2rn,rn}(\Th_S (F))$
où $\delta\colon \Th_S F\to \Th_S (F\otimes \xi)$ est le morphisme induit
par l'inclusion (admissible) $F\to F\otimes \xi$ de fibrés vectoriels
déduite du sous-fibré vectoriel de rang $1$ évident dans $\xi$.
\end{corollaire}

Cela se déduit immédiatement de la
proposition~\ref{proposition-compatibilite-tilde-thom} et de la
proposition~\ref{proposition-raffinement}.

\begin{corollaire}\label{corollaire-compatibilite-thom-general}
Soit $E$ un fibré vectoriel de rang $r$ sur $S$. Soit $X\in\Sm[S]$. Soit
$F$ un fibré vectoriel de rang $r$ sur $X$.
\[P_E(t_F)=\delta^\star t_{F\otimes \xi}\]
dans $\tH^{2rn,rn}(\Th_X (F))$
où $\delta\colon \Th_X F\to \Th_X (F\otimes \xi)$ est le morphisme induit
par l'inclusion (admissible) $F\to F\otimes \xi$ de fibrés vectoriels
déduite du sous-fibré vectoriel de rang $1$ évident dans $\xi$.
\end{corollaire}

La compatibilité au changement de base énoncée dans le
corollaire~\ref{corollaire-compatibilite-cb-f-lisse} permet de supposer que
$X=S$, ce qui nous ramène à l'énoncé du
corollaire~\ref{corollaire-compatibilite-thom}.

\subsection{Indépendance en le fibré $E$}
\label{subsection-independante-fibre}

\begin{proposition}
\label{proposition-independance-tilde-de-l-operation-en-le-fibre}
L'opération
$\tilde{P}_E\colon
\tH^{2r,r}(\mathcal X)\to \tH^{2rn,rn}(G
\backslash_{\L}(\mathcal X^{\wedge A}\times_S U))$ ne dépend
que du rang du fibré vectoriel $E$.
\end{proposition}

\begin{corollaire}
\label{corollaire-independance-de-l-operation-en-le-fibre}
L'opération cohomologique 
$P_E\colon \tH^{2r,r}(\mathcal X)\to \tH^{2rn,rn}(\mathcal X)$ ne dépend
que du rang du fibré vectoriel $E$.
\end{corollaire}

\begin{lemme}
Si $\mathcal X$ et $\mathcal Y$ sont deux objets de
$\Delta^\opp\Sc$, le morphisme canonique
\[G\backslash(\mathcal X^{\wedge A}\times U)\wedge
G\backslash(\mathcal Y^{\wedge A}\times U)\to
G\backslash((\mathcal X\wedge \mathcal Y)^{\wedge A}\times U)\]
est un isomorphisme dans $\Delta^\opp\Sc$.
\end{lemme}

On peut se ramener au cas particulier où $\mathcal X=X_+$ et $\mathcal
Y=Y_+$ avec $X$ et $Y$ des objets de $\Sm$. Le lemme
provient alors du fait que le morphisme canonique suivant dans $\Sm$ soit
un isomorphisme dans $\Sm$ :
\[G\backslash(X^A\times Y^A\times U)\isomto
G\backslash(X^A\times_S U)\times_S G\backslash (Y^A\times_S U)\;\text{.}\]
Cet isomorphisme peut se vérifier localement sur $S$ pour la topologie
étale. On peut donc supposer que $U$ est le $G$-torseur trivial,
auquel cas le morphisme est isomorphe à l'identité de $X^A\times Y^A$.

\begin{lemme}
On suppose que $\mathcal X$ et $\mathcal Y$ sont deux objets de
$\Delta^\opp\Sc$ et que $E$ et $F$ sont des fibrés vectoriels de rangs
respectifs $r$ et $s$ sur $S$. Alors, le diagramme évident qui suit
est commutatif :
\[
\xymatrix{\tH^{2r,r}(\mathcal X)\times\tH^{2s,s}(\mathcal Y)\ar[r]^-{\cup}
\ar[d]^{\tilde{P}_E\times \tilde{P}_F} &
\tH^{2(r+s),r+s}(\mathcal X\wedge \mathcal Y) \ar[d]^-{\tilde{P}_{E\oplus
F}} \\
\tH^{2rn,rn}(G\backslash(\mathcal X^{\wedge A}\times U))\times
\tH^{2sn,sn}(G\backslash(\mathcal Y^{\wedge A}\times U)) \ar[r]^-{\cup}&
\tH^{2(r+s)n,(r+s)n}
(G\backslash((\mathcal X\wedge \mathcal Y)^{\wedge A}\times U))}
\]
Autrement dit, si $\mathcal X$ et $\mathcal Y$ sont des objets de $\Hopt$,
si $x\in \tH^{2r,r}(\mathcal X)$ et $y\in
\tH^{2s,s}(\mathcal Y)$, alors on a une égalité :
\[\tilde{P}_E(x)\cdot \tilde{P}_F(y)=\tilde{P}_{E\oplus F}(x\cdot y)\]
dans $\tH^{2(r+s)n,(r+s)n}(G\backslash((\mathcal X\wedge \mathcal
Y)^{\wedge A}\times U))$, ou, ce qui revient au même, dans
$\tH^{2(r+s)n,(r+s)n}(G\backslash(\mathcal X^{\wedge A}\times U)\wedge
G\backslash(\mathcal Y^{\wedge A}\times U))$.
\end{lemme}

C'est évident.

\medskip

Nous pouvons maintenant démontrer la
proposition~\ref{proposition-independance-tilde-de-l-operation-en-le-fibre}.
Nous allons montrer que $\tilde{P}_E$ coïncide avec $\tilde{P}_{\OO_S^r}$
où $r$ est le rang de $E$.
Quitte à utiliser l'astuce de Jouanolou (cf.~\cite[Lemme~1.5]{jouanolou} et
\cite[Proposition~4.4]{weibel}), on peut supposer qu'il existe un
fibré vectoriel $E'$ et un entier $N$ tel que $E\oplus E'\simeq \OO_S^N$.
Notons $r$ et $s$ les rangs respectifs de $E$ et $E'$. Appliquons la
formule du lemme précédent aux décompositions $\OO_S^N=E\oplus E'$
et $\OO_S^N=\OO_S^r\oplus \OO_S^s$. Ici, $\mathcal X$ est arbitraire, mais on
pose $\mathcal Y=\Th_S E'$. Pour tout $x\in\tH^{2r,r}(\mathcal X)$,
on a :
\[\tilde{P}_E(x)\cdot \tilde{P}_{E'}(t_{E'})=\tilde{P}_{\OO_S^N}(x\cdot
t_{E'})) = \tilde{P}_{\OO_S^r}(x)\cdot \tilde{P}_{\OO_S^s}(t_{E'})
\]
Or, on sait déjà que
$\tilde{P}_{E'}(t_{E'})=\tilde{P}_{\OO_S^s}(t_{E'})=t_{E'\otimes \xi}$
par la proposition~\ref{proposition-compatibilite-tilde-thom}. La
multiplication par cette classe définissant un isomorphisme
\[\tH^{\star,\star}(\mathcal Z)\isomto \tH^{\star,\star}(\mathcal Z\wedge
\Th_S (E'\otimes\xi))\]
pour tout $\mathcal Z\in \Hopt$, en particulier $\mathcal Z=
G\backslash_{\L}
(\mathcal X^{\wedge A}\times_S U)$, on peut simplifier
l'identité $\tilde{P}_E(x)\cdot t_{E'\otimes \xi}=
\tilde{P}_{\OO_S^r}(x)\cdot t_{E'\otimes \xi}$
afin d'obtenir l'égalité $\tilde{P}_E(x)=\tilde{P}_{\OO_S^r}(x)$ voulue.

\section{Propriétés de l'opération totale}
\label{section-proprietes-operation-totale}
\setcounter{subsubsection}{0}

On fixe un corps parfait $k$. Dans cette section, on va s'intéresser plus
spécifiquement aux opérations suivantes :

\begin{definition}\label{definition-p-sur-k}
Soit $G$ un groupe fini agissant sur un ensemble fini $A$ de cardinal $n$.
Soit $S\in\Sm[k]$ (on note $a_S\colon S\to \Spec k$ le morphisme structual).
Soit $U$ un $G$-torseur étale au-dessus de $S$. Soit
$r\geq 0$. On note $P$ le morphisme dans $\Sm[k]^\opp\Ens$
\[K_r\to a_{S\star}K_{\xi^r}\]
déduit par adjonction de $P_E\colon K_E\to K_{E\otimes\xi}$
dans $\Sm[S]^\opp\Ens$ (cf. définition~\ref{definition-operation-totale})
dans le cas particulier
où $E$ est le fibré trivial de rang $r$
sur $S$ et $\xi$ le fibré vectoriel sur $S$ de la
définition~\ref{definition-xi}. (Quand le contexte l'exigera, on ajoutera
implicitement à la notation $P$ quelques indices pour préciser les données
utilisées.)
Ceci induit un morphisme $P\colon K_r\to \R
a_{S\star}K_{\xi^r}\overset{\text{Thom}}{\simeq} \R
a_{S\star}K_{rn}\simeq \R\SheafHom(S,K_{rn})$ dans $\Hopt[k]$ qui
correspond par adjonction à un morphisme $K_r\wedge S_+\to K_{rn}$ dans
$\Hopt[k]$ aussi noté $P$.
On note encore $P\colon \tH^{2r,r}(\mathcal X)\to \tH^{2rn,rn}(\mathcal
X\wedge S_+)$ l'opération induite pour tout $\mathcal X\in\Hopt[k]$.
\end{definition}

\subsection{Fonctorialité en $U$ et en $G$}

\begin{proposition}\label{proposition-fonctorialite-u}
Soit $G$ un groupe fini agissant sur un ensemble fini $A$ de cardinal $n$.
Soit $f\colon S'\to S$ un morphisme dans $\Sm[k]$. On suppose que $S$ et
$S'$ sont munis respectivement de $G$-torseurs étales $U$ et $U'$ et qu'il
existe un isomorphisme de $G$-torseurs $f^\star U\simeq U'$.

Alors, le diagramme suivant est commutatif pour tout $r\geq 0$ et $\mathcal
X\in \Hopt$ :
\[
\xymatrix{\tH^{2r,r}(\mathcal X)\ar[r]^-{P_U} \ar[rd]_-{P_{U'}}
& \tH^{2rn,rn}(\mathcal X\wedge S_+)
\ar[d]^{f^\star} \\
& \tH^{2rn,rn}(\mathcal X\wedge S'_+)}
\]
\end{proposition}

Ceci résulte essentiellement du corollaire~\ref{corollaire-compatibilite-cb}.

\begin{proposition}\label{proposition-fonctorialite-g}
Soit $G$ un groupe fini agissant sur un ensemble fini $A$ de cardinal $n$.
Soit $H$ un sous-groupe de $G$. Soit $S\in\Sm[k]$ un schéma muni d'un
$G$-torseur étale $U$. On note $S'=H\backslash U$. On peut considérer aussi
$U$ comme un $H$-torseur étale au-dessus de $S'$. On note $p\colon S'\to S$
le morphisme évident. Alors, le diagramme suivant est commutatif pour tous
$r\geq 0$ et $\mathcal X\in\Hopt[k]$ :
\[\xymatrix{\tH^{2r,r}(\mathcal X)\ar[r]^-{P_G}\ar[rd]_-{P_H} &
\tH^{2rn,rn}(\mathcal X\wedge S_+) \ar[d] ^{p^\star}\\
& \tH^{2rn,rn}(\mathcal X\wedge S'_+)}\]
\end{proposition}

Ceci résulte de l'évidente commutativité du diagramme suivant :
\[\xymatrix{K_r\ar[r]^-{P_G} \ar[rd]_-{P_H}
& a_{S\star}K_{\xi^r}\ar[d]\\
& a_{S'\star}K_{\xi'^r}}\]
où $a_S\colon S\to \Spec k$ et $a_{S'}\colon S'\to \Spec k$ sont les
morphismes structuraux, $\xi'=p^\star \xi$ et où le morphisme de droite
s'identifie au morphisme $a_{S\star}K_{\xi^r}\to a_{S\star}p_\star
K_{\xi'^r}$ obtenu à partir du morphisme d'adjonction
$K_{\xi^r}\to p_\star p^\star K_{\xi^r}=p_\star K^{\xi'^r}$.

En combinant les propositions~\ref{proposition-fonctorialite-u}
et \ref{proposition-fonctorialite-g}, on obtient la compatibilité plus
générale :

\begin{corollaire}\label{corollaire-fonctorialite-u-g}
Soit $G'\to G$ un morphisme injectif de groupes\;\footnote{Si $G'\to G$
n'est pas injectif, la compatibilité est vraie aussi, mais elle ne fait que
montrer que ce cas, et donc celui d'une action infidèle de $G'$ sur $A$,
n'est pas très intéressant.}. On suppose que l'on s'est 
donné une action de $G$ sur un
ensemble fini $A$ de cardinal $n$. Par restriction, $A$ est également muni
d'une action de $G'$.
On suppose donnés un morphisme $f\colon S'\to S$ dans $\Sm[k]$, un
$G$-torseur étale $U$ sur $S$, un $G'$-torseur étale $U'$ sur $S'$
et un $S$-morphisme
$U'\to U$ qui soit $G'$-équivariant. Autrement dit, on dispose d'un
diagramme :
\[\xymatrix@C=1.5cm{U'\ar[d]_{G'\text{-torseur}}\ar[r]^{G'\text{-équiv.}} &
U\ar[d]^{G\text{-torseur}}\\
S'\ar[r]^f & S}\]
Alors, le diagramme suivant commute pour tous $r\geq 0$ et $\mathcal
X\in\Hopt[k]$ :
\[\xymatrix{\tH^{2r,r}(\mathcal X)\ar[r]^-{P_G} \ar[dr]_-{P_{G'}}&
\tH^{2rn,rn}(\mathcal X\wedge S_+)\ar[d]^{f^\star} \\
& \tH^{2rn,rn}(\mathcal X\wedge S'_+)
} \] 
\end{corollaire}

On peut identifier le carré ci-dessus à une composition de deux carrés :
\[\xymatrix@C=1.5cm{U'\ar[d]_{G'\text{-torseur}} \ar[r]^{G'\text{-équiv.}}&
U\ar[d]_{G'\text{-torseur}}\ar@{=}[r]^{G'\text{-équiv.}} &
U\ar[d]^{G\text{-torseur}}\\
S'\ar[r] & G'\backslash U \ar[r] & S}\]
Pour établir la compatibilité énoncée dans le corollaire, on
se ramène ainsi aux
propositions~\ref{proposition-fonctorialite-u} et
\ref{proposition-fonctorialite-g} qui correspondent respectivement
aux carrés de gauche et de droite.

\subsection{Compatibilité à la multiplication}

\begin{proposition}\label{proposition-compatibilite-multiplication}
Soit $G$ un groupe fini agissant sur un ensemble fini $A$ de cardinal $n$.
Soit $S\in \Sm[k]$ un schéma muni d'un $G$-torseur étale $U$.
On note $P_G$
l'opération associée $\tH^{2r,r}(\mathcal X)\to \tH^{2rn,rn}(\mathcal
X\wedge S_+)$ pour $r\geq 0$ et $\mathcal X\in \Hopt[k]$.
Pour tous $\mathcal X\in\Hopt[k]$, $\mathcal Y\in \Hopt[k]$, $r,s\geq 0$, $x\in
\tH^{2r,r}(\mathcal X)$, $y\in \tH^{2s,s}(\mathcal Y)$, on a l'identité
suivante dans $\tH^{2(r+s),r+s}(\mathcal X\wedge \mathcal Y\wedge S_+)$ :
\[P_G(x\cdot y)=\Delta^\star(P_G(x)\cdot P_G(y))\]
où $\Delta\colon \mathcal X\wedge \mathcal Y\wedge S_+\to
(\mathcal X\wedge S_+)\wedge (\mathcal Y\wedge S_+)$
est le morphisme canonique provenant de la diagonale de $S$.
\end{proposition}

Ceci résulte trivialement de la
proposition~\ref{proposition-compatibilite-produit}.

(Plus loin, on se dispensera d'écrire $\Delta^\star$ pour ne pas alourdir
les notations.)

\subsection{Action sur les classes de Thom}

\begin{proposition}\label{proposition-action-classes-de-thom}
Soit $G$ un groupe fini agissant sur un ensemble fini $A$ de cardinal $n$.
Soit $S\in \Sm[k]$ un schéma muni d'un $G$-torseur étale $U$. On note $P_G$
l'opération associée $\tH^{2r,r}(\mathcal X)\to \tH^{2rn,rn}(\mathcal
X\wedge S_+)$ pour $r\geq 0$ et $\mathcal X\in \Hopt[k]$.

Soit $X\in \Sm[k]$. Soit $F$ un fibré vectoriel de rang $r$ sur $X$. On
note $t_F\in \tH^{2r,r}(\Th_X F)$ la classe de Thom de $F$. Alors,
\[P_G(t_F)=t_F\cdot c_{r(n-1)}(F\boxtimes (\xi/\OO))\in \tH^{2rn,rn}(\Th_X
F\wedge S_+)\;\text{.}\]
En particulier, si $F$ est le fibré trivial de rang $r$, on a :
$P(t_F)=t_F\cdot c_{n-1}(\xi)^r$.
\end{proposition}

D'après le corollaire~\ref{corollaire-compatibilite-thom-general}, $P(t_F)$
s'identifie à l'image inverse de $t_{F\boxtimes \xi}$ par le morphisme
\[\Th_X F\wedge S_+\simeq \Th_{X\times S}(F\boxtimes \OO)\to \Th_{X\times
S}(F\boxtimes \xi)\]
associé à l'inclusion $F\boxtimes \OO \to F\boxtimes \xi$ de fibrés
vectoriels sur $X\times S$ déduite du sous-fibré trivial de rang $1$
évident dans $\xi$. La formule de la proposition résulte donc du lemme
suivant :

\begin{lemme}
Soit $i\colon E\to F$ un monomorphisme admissible entre fibrés vectoriels
sur $X\in \Sm[k]$. Cette inclusion $i$ induit un morphisme $\Th i\colon \Th_X
E\to \Th_X F$.

Alors, dans $\tH^{\star,\star}(\Th_X E)$, on a l'identité
\[(\Th i)^\star (t_F)=t_E\cdot c_r(F/E)\]
où $r$ est le rang de $F/E$ et $(\Th i)^\star$ le morphisme de restriction
$\tH^{\star,\star}(\Th_X F)\to \tH^{\star,\star}(\Th_X E)$.
\end{lemme}

Les classes $t_E$ et $t_F$ peuvent être identifiées respectivement à des
éléments de $H^{\star,\star}(\mathbf{P}(F\oplus \OO_X))$ et de
$H^{\star,\star}(\mathbf{P}(E\oplus \OO_X))$. On note
$\xi=c_1(\OO(1))$ la première classe de Chern du fibré en droites
fondamental de ces deux fibrés projectifs.
Notons $e$ le rang de $E$. La multiplicativité des polynômes de Chern de
$E$ et de $F/E$ fait que dans $H^{\star,\star}(\mathbf{P}(F\oplus
\OO_X))$, on a :
\[t_F = (\xi^e+c_1(E)\xi^{e-1}+\dots +c_e(E))\cdot
(\xi^r+c_1(F/E)\xi^{r-1}+\dots+c_r(F/E))\]
En appliquant $(\Th i)^\star$, on obtient une égalité dans
$H^{\star,\star}(\mathbf{P}(E\oplus \OO_X))$ :
\[(\Th i)^\star (t_F)=(t_E \xi)\cdot
(\xi^{r-1}+c_1(F/E)\xi^{r-2}+\dots+c_{r-1}(F/E))+t_E\cdot c_r(F/E)\;\text{.}\]
On peut conclure en utilisant que
$t_E\xi=\xi^{e+1}+c_1(E)\xi^{e}+\dots+c_r(E)\xi=0\in
H^{\star,\star}(\mathbf{P}(E\oplus \OO_X))$, fait qui résulte par exemple
de la définition des classes de Chern de $E\oplus \OO_X$ et de l'identité
$c_i(E\oplus \OO_X)=c_i(E)$ pour tout $i\geq 1$.

\subsection{Cas du classifiant $\Bet G$}
\label{subsection-bet-g}

Soit $G$ un groupe fini agissant sur un ensemble fini $A$ à $n$ éléments.
On va montrer ici que, cette action de $G$ étant fixée, les opérations
$P$ associées aux différents $G$-torseurs $U$ sur $S\in\Sm[k]$ se déduisent
toutes d'une opération universelle au niveau du classifiant étale $\Bet G$
(cf. \cite[p.~130]{morel-voevodsky}).

On rappelle que cet objet $\Bet G\in \Hosimppt$ a la vertu que pour tout
préfaisceau simplicial $\mathcal X$ sur $\Sm[k]$,
$\Hom_{\Hosimp[k]}(\mathcal X, \Bet G)\simeq H^1_{\et}(a_\et \mathcal X,G)$
où $H^1_{\et}(a_\et\mathcal X,G)$ désigne l'ensemble des classes
d'isomorphismes de $G$-torseurs étales au-dessus du faisceau (simplicial)
étale $a_\et \mathcal X$ associé à $\mathcal X$. (Dans le cas favorable où
l'ordre de $G$ est inversible dans $k$, on a une bijection
$H^1_{\et}(X,G)\simeq H^1_{\et}(\mathbf{A}^1\times X,G)$ pour tout
$X\in\Sm$ et l'objet $\Bet G$ est $\mathbf{A}^1$-local (cf.
\cite[Proposition~4.1, p.~137]{morel-voevodsky}), ce qui permet de
remplacer ci-dessus les $\Hom$ de la catégorie $\Hosimp[k]$ par ceux de
la catégorie $\Ho[k]$.)

Dans cette sous-section, on identifiera un $k$-espace vectoriel de
dimension finie $W$ à l'espace affine $\Spec \AlgSym^\star W^{\vee}$ qui
lui est associé ($\AlgSym^\star W^{\vee}$
est l'algèbre symétrique du dual de $W$).
La représentation régulière de $G$ est notée $k[G]$. Pour
tout $d\geq 1$, on peut ainsi noter $U_d$ l'ouvert (dense) de l'espace
affine $k[G]^{\oplus d}$ sur lequel $G$ agit librement. On définit le
schéma $S_d=G\backslash U_d$ qui est muni du $G$-torseur étale $U_d$. En
utilisant les injections évidentes $k[G]^d\to k[G]^{d+1}$, on fait des
$U_d$ et $S_d$ un système inductif pour $d\geq 1$. Les colimites $U_\infty$
et $S_\infty$ de ces systèmes, calculées dans la catégorie des préfaisceaux
sur $\Sm[k]$ sont notées respectivement $\Egm G$ et $\Bgm G$. Bien entendu,
$\Egm G$ est un $G$-torseur étale au-dessus de $\Bgm G$ et la propriété
caractéristique de $\Bet G$ énoncée plus haut associe à ce $G$-torseur
$\Egm G$ un morphisme canonique $\Bgm G\to \Bet G$ dans $\Hosimp[k]$.

\begin{theoreme}[{\cite[Proposition~2.6, p. 135]{morel-voevodsky}}]%
\label{theoreme-bg}
Le morphisme canonique $\Bgm G\to \Bet G$ est une
$\mathbf{A}^1$-équivalence faible. Autrement dit, dans $\Ho[k]$, on a un
isomorphisme $\Bgm G\simeq \Bet G$.
\end{theoreme}

Pour tout $d\geq 1$, on peut appliquer la construction de la
définition~\ref{definition-p-sur-k} à l'action de $G$ sur $A$ et au
$G$-torseur étale $U_d$ sur $S_d$, ce qui fournit
un morphisme dans $\Ho[k]$ pour tout $r\geq 0$ :
\[P_{G,d}\colon K_r\to \R\SheafHom(S_d,K_{rn})\;\text{.}\]
D'après la proposition~\ref{proposition-fonctorialite-u},
ces morphismes sont compatibles avec les morphismes de transition du
système inductif $(S_d)_{d\geq 1}$. Il existe donc un morphisme $K_r\to
\R\SheafHom(\Bgm G,K_{rn})$ qui induise tous les morphismes $P_{G,d}$ après
composition des $\R\SheafHom$ par les différents morphismes évidents $S_d\to
\Bgm G$ (cela résulte de la suite exacte de Milnor, cf.
\cite[Proposition~2.15, Chapter~VI]{goerss-jardine}). Compte tenu de
l'isomorphisme $\Bgm G \simeq \Bet G$ dans $\Hopt$, ceci permet de procéder
à la définition suivante :

\begin{definition}\label{definition-p-univ}
Pour tout $r\geq 0$, on choisit un morphisme
$P_G^{\univ}\colon K_r\to \R\SheafHom(\Bet
G,K_{rn})$ dans $\Ho$ tel que pour tout $d\geq 1$, $P_{G,d}$ soit le
composé $K_r\vers{P_G^{\univ}} \R\SheafHom(\Bet
G,K_{rn})\vers{[U_d]^\star}\R\SheafHom(S_d,K_{rn})$.
\end{definition}

\begin{proposition}\label{proposition-compatibilite-bg}
Soit $S\in\Sm[k]$. Soit $U$ un $G$-torseur étale sur $S$. À la classe du
$G$-torseur $U$ correspond un morphisme $[U]\colon S\to \Bet G$ dans
$\Hosimp[k]$. Le diagramme suivant est commutatif :
\[\xymatrix{K_r\ar[r]^-{P^{\univ}_G} \ar[rd]_-{P_{G,U}}&
\R\SheafHom(\Bet G,K_{rn})\ar[d]^{[U]^\star}\\
& \R\SheafHom(S,K_{rn})}\]
\end{proposition}

La propriété énoncée est vraie par définition de $P^\univ_G$ dans le cas
des $G$-torseurs $U_d$ sur $S_d$ pour tout $d\geq 1$. Si $U$ est l'image
inverse de $U_d$ par un morphisme $f\colon S\to S_d$ dans $\Sm[k]$, on
déduit la compatibilité pour $S$ de celle pour $S_d$ en utilisant la
commutativité des deux triangles dans le diagramme ci-dessous :
\[\xymatrix@C=1.5cm{K_r\ar[r]^-{P^\univ_G}\ar[rd]^-{P_{G,U_d}}
\ar[rdd]_-{P_{G,U}} &
\R\SheafHom(\Bet G,K_{rn}) \ar[d]^{[U_d]^\star} \\
&\R\SheafHom(S_d,K_{rn})\ar[d]^{f^\star} \\
&\R\SheafHom(S,K_{rn}) }\]
Le triangle du haut commute parce que c'est le cas du
torseur $U_d$ sur $S_d$ et celui du bas d'après la
proposition~\ref{proposition-fonctorialite-u}.

Pour traiter le cas d'un $G$-torseur arbitraire $U$ sur $S\in\Sm[k]$, on
peut commencer par supposer que $S$ est affine en utilisant l'astuce de
Jouanolou. Le lemme suivant montre que l'on est alors
dans le cas traité ci-dessus :

\begin{lemme}\label{lemme-morphisme-g-torseur-dans-u-d}
Soit $S$ un $k$-schéma affine de type fini. Soit $U$ un $G$-torseur étale
sur $S$. Alors, il existe un entier $d\geq 1$ et un morphisme
$G$-équivariant $U\vers{\tilde{f}} U_d$. Ainsi, le $G$-torseur $U$
s'identifie à l'image inverse de $U_d$ par le morphisme $f\colon S\to S_d$
induit par $\tilde{f}$ par passage au quotient par $G$. (On peut en outre
s'arranger pour que $f$ et $\tilde{f}$ soient des immersions fermées.)
\end{lemme}

Notons $A$ l'anneau du schéma affine $U$. Comme $A$ est une $k$-algèbre de
type fini, il existe un sous-$k$-espace vectoriel de dimension finie $W$ de
$A$ tel que le morphisme d'algèbres $\AlgSym^\star W\to A$ associé soit
surjectif. L'action à gauche de $G$ sur $U$ se traduit en une action à
droite de $G$ sur $A$. Quitte à remplacer $W$ par $\sum_{g\in G}g^\star W$,
on peut supposer que $W\subset A$ est stable par l'action de $G$. Comme
module à droite sur $k[G]$, $W$ est évidemment de type fini, il existe donc
un morphisme surjectif $k[G]^d\to W$ de modules à droites sur $k[G]$. Ceci
induit un morphisme surjectif et $G$-équivariant de $k$-algèbres :
\[\AlgSym^\star(k[G]^d)\to A\;\text{.}\] Géométriquement, compte tenu de
l'autodualité de la représentation régulière, cela correspond à une
immersion fermée $G$-équivariante $U\to k[G]^{\oplus d}$. Compte tenu de
cela, le fait que l'action de $G$ sur $U$ soit libre implique que
l'immersion fermée $W\to k[G]^d$ se factorise par l'ouvert $U_d$. On a
ainsi construit le morphisme voulu $\tilde{f}\colon U\to U_d$.

\bigskip

On introduit ici une notion \emph{ad hoc} qui nous sera utile plus loin :

\begin{definition}
On suppose donnés deux objets $\mathcal X$ et $\mathcal Z$ de $\Ho[k]$ et
un objet $\mathcal Y$ de $\Hosimp[k]$ (qui induit donc aussi un objet de
$\Ho[k]$). On dira de deux morphismes $\mathcal X\to \R\SheafHom(\mathcal
Y,\mathcal Z)$ qu'ils sont presqu'égaux si pour tout $Y\in\Sm[k]$ et tout
morphisme $Y\to\mathcal Y$ dans $\Hosimp[k]$, les morphismes composés dans
$\Ho[k]$
\[\mathcal X\to \R\SheafHom(\mathcal Y,\mathcal Z)\to\R\SheafHom(Y,\mathcal
Z)\]
sont égaux.
\end{definition}

On s'intéresse tout particulièrement ici à cette définition dans le cas où
$\mathcal Y=\Bet G$. Dans ce cas, deux morphismes $\mathcal X\to
\R\SheafHom(\Bet G,\mathcal Z)$ sont presqu'égaux si pour tout $S\in
\Sm[k]$ et tout $G$-torseur $U$ sur $S$, les morphismes composés
\[\mathcal X\to \R\SheafHom(\Bet G,\mathcal Z)\vers{[U]^\star}
\R\SheafHom(S,\mathcal Z)\]
sont égaux.

La proposition~\ref{proposition-compatibilite-bg} montre que le choix du
morphisme $P^\univ_G\colon K_r\to \R\SheafHom(\Bet G,K_{rn})$ correspond
exactement au choix d'un représentant dans une classe de morphismes à
presqu'égalité près, les morphismes composés $K_r\to \R\SheafHom(S,K_{rn})$
étant en effet imposés pour chaque $G$-torseur $U$ sur $S$.

À partir de cette presqu'égalité de morphismes, on peut introduire une
notion évidente de diagramme \guil{presque commutatif} dont on fait usage
dans la proposition suivante :

\begin{proposition}
Soit $i\colon H\to G$ un morphisme injectif de groupes. Le $G$-ensemble $A$
est muni par restriction d'une structure de $H$-ensemble. Par la
définition~\ref{definition-p-univ}, on dispose ainsi non seulement d'un
morphisme $P^\univ_G$ mais aussi d'un morphisme $P^\univ_H$, et ils font
presque commuter le diagramme suivant :
\[\xymatrix{K_r\ar[r]^-{P^\univ_G}\ar[rd]_-{P^\univ_H} &
\R\SheafHom(\Bet G,K_{rn})\ar[d]^{(\Bet i)^\star}\\
& \R\SheafHom(\Bet H,K_{rn})}\]
\end{proposition}

Soit $U$ un $H$-torseur sur $S\in\Sm[k]$. 
La classe du $H$-torseur $U$ définit un morphisme $S\vers{[U]} \Bet H$. Par
la fonctorialité sur les torseurs, le morphisme composé $S\vers{[U]}\Bet
H\to \Bet G$ correspond au $G$-torseur défini comme le quotient de $G\times
U$ par l'action de $H$ sur $G\times U$ donnée par $h.(g,u)=(gh^{-1},hu)$, ce
quotient étant muni de l'action à gauche donnée par la multiplication à
gauche dans $G$.

Après composition avec $[U]^\star\colon \R\SheafHom(\Bet H,K_{rn})\to
\R\SheafHom(S,K_{rn})$, les deux morphismes dont on doit montrer la
presqu'égalité deviennent deux morphismes $K_r\to \R\SheafHom(S,K_{rn})$.
Ce sont deux cas particuliers de la construction $P$ pour des torseurs sur
$S$ : l'un pour le $H$-torseur $U$, l'autre pour le $G$-torseur $(G\times
U)/H$ défini ci-dessus. L'égalité de ces deux morphismes découle du
corollaire~\ref{corollaire-fonctorialite-u-g} appliqué au diagramme
$H$-équivariant évident :
\[\xymatrix@C=1.5cm{U \ar[d]_{H\text{-torseur}}\ar[r]^-{u\longmapsto (1,u)}
&\ar[d]^{G\text{-torseur}} (G\times U)/H \\
S\ar@{=}[r] & S }\]

\begin{remarque}
Il est vraisemblablement possible de se débarasser des \guil{presque}
apparaissant ci-dessus et de définir le morphisme
$P^\univ_G\colon K_r\to \R\SheafHom(\Bet G,K_{rn})$ sans équivoque. Pour
cela, il suffirait d'obtenir une version convenable de l'isomorphisme de
Thom (cf. proposition~\ref{proposition-isomorphisme-thom-relatif})
pour une famille compatible de fibrés vectoriels sur un diagramme de schéma
(dont la colimite serait $\Bgm G$). Ceci ne présente pas un intérêt décisif
puisque dans le cas le plus intéressant du groupe symétrique $\Sym \ell$ et
de $\Lambda=\mathbf{Z}/\ell\mathbf{Z}$, le $\lim^1$ mesurant l'obstruction
à l'unicité est nul (cf.~corollaire~\ref{corollaire-cohomologie-b-sym-ell}).
\end{remarque}

\subsection{Théorème de symétrie}

\begin{theoreme}\label{theoreme-symetrie}
Soit $G$ un groupe agissant sur un ensemble $A$ à $n$ éléments. Soit
$S\in\Sm[k]$ un schéma muni d'un $G$-torseur $U$.
Pour tout $\mathcal X\in\Hopt[k]$, pour tout $u\in\tH^{2r,r}(\mathcal X)$,
on peut considérer $P_G(u)\in\tH^{2rn,rn}(\mathcal X\wedge S_+)$ puis
$P_G(P_G(u))\in\tH^{2rn^2,rn^2}(\mathcal X\wedge (S\times S)_+)$. Si on
note $\tau\colon S\times S\to S\times S$ l'échange des deux facteurs et
$\tau^\star$ son action sur $\tH^{\star,\star}(\mathcal X\wedge
(S\times S)_+)$, on a :
\[\tau^{\star}(P_G(P_G(u))=P_G(P_G(u)))\;\text{.}\]
\end{theoreme}

La proposition suivante permet d'interpréter $P_G\circ P_G$ comme une seule
construction $P$, appliquée au groupe $G\times G$, et ainsi de réduire le
théorème~\ref{theoreme-symetrie} à une formule $\tau^\star \circ P_{G\times
G}=P_{G\times G}$.

\begin{proposition}\label{proposition-compatibilite-composition-sur-k}
Outre une action d'un groupe fini $G$ sur un ensemble $A$ à $n$ éléments,
on se donne un groupe fini $G'$ agissant sur un ensemble $A'$ à $n'$
éléments. Soit $S\in\Sm[k]$ un schéma muni d'un $G$-torseur $U$. Soit
$S'\in\Sm[k]$ un schéma muni d'un $G'$-torseur $U'$. Il en résulte une
manière évidente de considérer $U\times U'$ comme un $G\times G'$-torseur
sur $S\times S'$ et $A\times A'$ comme un $G\times G'$-ensemble.

Pour tout $r\geq 0$, le diagramme suivant est commutatif dans $\Hopt$ :
\[\xymatrix{K_r\ar[r]^{P_G}\ar[d]_{P_{G\times G'}} &
\R\SheafHom(S,K_{rn})\ar[d]^{\R\SheafHom(S,P_{G'})}\\
\R\SheafHom(S\times
S',K_{rnn'})\ar@{=}[r]&\R\SheafHom(S,\R\SheafHom(S',K_{rnn'})}\]
\end{proposition}

Cela résulte immédiatement de \ref{proposition-compatibilite-composition}
et de la compatibilité des constructions $P$ aux changements de base
associés aux morphismes de projection $S\times S'\to S$ et $S\times S'\to
S'$ (cf. corollaire~\ref{corollaire-compatibilite-cb}).

\medskip

\paragraph{Première démonstration du théorème~\ref{theoreme-symetrie}}
Observons que le diagramme
suivant est commutatif dans $\Hosimp[k]$, où $\tau$ désigne l'échange des
facteurs $S$ ou $G$ :
\[\xymatrix{S\times S\ar[d]_\tau \ar[r]^-{[U\times U]} & \Bet(G\times
G)\ar[d]^{\Bet\tau}\\
S\times S\ar[r]^-{[U\times U]}& \Bet(G\times G)}\]
En effet, si on note $(U\times U)^\tau$ le $G\times G$-torseur déduit du
$G\times G$-torseur $U\times U$ sur $S\times S$ par le changement de groupe
associé au morphisme $\tau\colon G\times G\to G\times G$ (c'est-à-dire que
le schéma sous-jacent à $(U\times U)^\tau$ est $U\times U$ et l'action est
définie par la formule $(g,h).(u,v)=(hu,gv)$), alors le $G\times G$-torseur
$(U\times U)^\tau$ s'identifie aussi à l'image inverse du $G\times
G$-torseur $U\times U$ par $\tau\colon S\times S\to S\times S$.
(L'isomorphisme entre les deux $G\times G$-torseurs est donné par l'échange
des deux facteurs $U$.)

On en déduit que pour démontrer le théorème, il suffit de montrer que le
diagramme suivant est presque commutatif :
\[\xymatrix@C=1.5cm{K_r\ar[rd]_-{P^\univ_{G\times G}}
\ar[r]^-{P^\univ_{G\times G}} & \R\SheafHom(\Bet(G\times
G),K_{rn^2}) \ar[d]^{\R\SheafHom(\Bet \tau,K_{rn^2})}\\
& \R\SheafHom(\Bet(G\times G),K_{rn^2})} \]
Notons $\widetilde{G\times G}$ le produit semi-direct $(G\times G)\rtimes
\mathbf{Z}/2\mathbf{Z}$ où l'élément non trivial de
$\mathbf{Z}/2\mathbf{Z}$ agit par l'échange
$\tau\colon G\times G\to G\times G$ des deux facteurs $G$. Notons
$\tilde{\tau}$ l'automorphisme intérieur de $\widetilde{G\times G}$ donné
par la conjugaison par l'élément non trivial de
$\mathbf{Z}/2\mathbf{Z}\subset \widetilde{G\times G}$.
Le diagramme de groupes suivant est commutatif :
\[\xymatrix{G\times G\ar[r]\ar[d]^\tau & \widetilde{G\times
G}\ar[d]^{\tilde{\tau}}\\
G\times G\ar[r] & \widetilde{G\times G}}\]
Comme $\tilde{\tau}$ est un automorphisme intérieur, il induit l'identité
sur $\Bet(\widetilde{G\times G})$, d'où un diagramme commutatif dans
$\Hosimp[k]$ :
\[\xymatrix{\Bet (G\times G)\ar[rd] \ar[d]^\tau& \\
\Bet (G\times G)\ar[r] & \Bet (\widetilde{G\times G}) }\]
Comme l'action de $G\times G$ sur $A\times A$ s'étend en une action de
$\widetilde{G\times G}$ (en faisant agit l'élément non trivial de
$\mathbf{Z}/2\mathbf{Z}\subset \widetilde{G\times G}$ par l'échange
$\tau\colon A\times A\to A\times A$), on en déduit un diagramme presque
commutatif :
\[\xymatrix{
& & & \R\SheafHom(\Bet (G\times G),K_{rn^2}) \ar[dd]^{\R\SheafHom(\Bet
\tau,K_{rn^2})}\\
K_r\ar@/^2pc/[rrru]^{P^{\univ}_{G\times G}}
\ar@/_2pc/[rrrd]_{P^{\univ}_{G\times G}}
\ar[rr]_-{P^\univ_{\widetilde{G\times G}}} & & \R\SheafHom(\Bet (\widetilde{G\times
G}),K_{rn^2})\ar[ur]\ar[dr] & \\
& & & \R\SheafHom(\Bet (G\times G),K_{rn^2}) \\
}\]
La presque commutativité du contour externe est ce qu'il fallait démontrer.

\begin{remarque}
La première démonstration est dans le même ordre d'idées que celle de
\cite[\S{}7]{voevodsky-reduced}, mais utilise de façon un peu plus
systématique les classifiants des groupes $G\times G$ et $\widetilde{G\times
G}$. La version \cite[Corollary~7.3]{voevodsky-reduced} de la
proposition~\ref{proposition-compatibilite-composition-sur-k} énonçait les
hypothèses de façon peu claire. La compatibilité n'était obtenue qu'après
multiplication par une certaine classe qui dans certains cas était
simplifiable. On a obtenu cette compatibilité plus directement ici, et
c'est là un des intérêts d'avoir introduit les opérations $P$ au niveau
d'espaces $K_E$ associés à des fibrés vectoriels $E$ arbitraires plutôt
que dans le seul cas des fibrés triviaux qui donnent les modèles usuels des
espaces d'Eilenberg-MacLane $K_r$.
\end{remarque}

\paragraph{Deuxième démonstration du théorème~\ref{theoreme-symetrie}}

Le schéma $U\times U$ est considéré comme un $G\times G$-torseur sur
$S\times S$ de la façon évidente. À chaque $G\times G$-ensemble $X$ est
associé un morphisme
\[P_X\colon K_r\to \R\SheafHom(S\times S,K_{r\cdot \#X})\]
qui ne dépend bien évidemment que de la classe d'isomorphisme du
$G\times G$-ensemble $X$.

Notons $A\times A$ le $G\times G$-ensemble évident et $(A\times A)^\tau$ le
même ensemble $A\times A$ muni de l'action de $G\times G$ donnée par la
formule $(g,h).(a,b)=(ha,gb)$.
D'après le corollaire~\ref{corollaire-fonctorialite-u-g} appliqué aux
morphismes d'échanges $\tau\colon G\times G\to G\times G$ et $\tau\colon
U\times U\to U\times U$ (induisant $\tau\colon S\times S\to S\times S$), on
obtient le diagramme suivant dans $\Hopt[k]$ :
\[\xymatrix{K_r\ar[r]^-{P_{A\times A}}\ar[rd]_-{P_{(A\times A)^\tau}} &
\R\SheafHom(S\times S,K_{rn^2})\ar[d]^{\tau^\star }\\
& \R\SheafHom(S\times S,K_{rn^2})}\]
Comme les $G\times G$-ensembles $A\times A$ et $(A\times A)^\tau$ sont
isomorphes (\emph{via} l'échange $\tau\colon A\times A\to A\times A$), on a
$P_{A\times A}=P_{(A\times A)^\tau}$, ce qui fournit la compatibilité
annoncée.

\subsection{Annulation du Bockstein}
\label{subsection-annulation-bockstein}

\begin{theoreme}\label{theoreme-annulation-bockstein}
Soit $\ell$ un nombre premier. Soit $S\in\Sm[k]$ un schéma muni d'un $\Sym
\ell$-torseur étale $U$. On note $P_\ell$ l'opération en cohomologie
à coefficients dans $\mathbf{Z}/\ell\mathbf{Z}$
associée à ce torseur et à l'action évidente
de $\Sym \ell$ sur $\{1,\dots,\ell\}$
(cf.~définition~\ref{definition-p-sur-k}). Alors, pour
tous $r\geq 0$ et $\mathcal X\in \Hopt[k]$, la composition suivante est
nulle :
\[\tH^{2r,r}(\mathcal X)\vers{P_\ell} \tH^{2r\ell,r\ell}(\mathcal
X\wedge S_+)\vers {\beta} \tH^{2r\ell+1,r\ell}(\mathcal X\wedge
S_+)\]
où $\beta$ est le Bockstein. Autrement dit, pour tout $u\in
\tH^{2r,r}(\mathcal X)$, la classe $P_\ell(u)$ est dans l'image du
morphisme évident
$\tH^{2r\ell,r\ell}(\mathcal X\wedge S_+,\mathbf{Z}/\ell^2\mathbf{Z})\to
\tH^{2r\ell,r\ell}(\mathcal X\wedge S_+,\mathbf{Z}/\ell\mathbf{Z})$
\end{theoreme}

La démonstration du théorème \cite[Theorem~8.4]{morel-voevodsky} contient
des erreurs : il y a par exemple
un problème d'exactitude à gauche dans
la suite \cite[(8.2), p.~30]{morel-voevodsky}. La démonstration
présentée ci-dessous suit peut-être plus fidèlement encore la démonstration
originale de \cite[\S{}4, Chapter VII]{steenrod} dans le cadre
classique. Pour ce faire, nous introduisons des complexes qui permettent de
représenter les classes de cohomologie motivique par des cocycles.

\begin{definition}
Soit $S\in \Sm[k]$. Soit $K$ un objet de $\CompNeg(\PST[S])$ (par exemple
un objet $\mathcal F$ de $\PST[S]$ identifié à un complexe
concentré en degré $0$).
Pour tout $\mathcal X\in \Delta^\opp \Sc$, on note $C(\mathcal X,K)$ le
complexe simple (défini en termes de produits) du bicomplexe 
$C^{p,q}(\mathcal X,K)=K_{-q}(\mathcal X_p)$ où les différentielles qui
incrémentent $p$ sont définies par des sommes alternées de flèches induites
par la structure simpliciale de $\mathcal X$ (et la structure de
préfaisceau sur $K$)
et où les différentielles qui incrémentent $q$ sont
induites par celles de $K$\;\footnote{On a commis ici un petit abus de
notations en faisant comme si on pouvait simplement évaluer les
préfaisceaux avec transferts $K_{-q}$ sur les objets $\mathcal X_p$ qui
sont des sommes directes de faisceaux représentables. Cette notation cache
en réalité un $\Hom$ dans la catégorie des préfaisceaux
d'ensembles sur $\Sm$.}.
Ce complexe satisfaisant des fonctorialités
évidentes par rapport à $\mathcal X$ et $K$, on peut noter $\tC(\mathcal
X,K)$ le noyau de l'épimorphisme scindé $C(\mathcal X,K)\to C(S,K)$ induit
par le point-base de $\mathcal X$.
\end{definition}

Le foncteur $\Ztr\colon\Sm[S]\to \PST[S]$ donne naissance à un foncteur
$\tZtr\colon \Sc\to\PST[S]$ qui pour $X\in \Sm[S]$ est tel que
$\tZtr(X_+)=\Ztr(X)$ et qui commute aux sommes directes. On étend ceci en
un foncteur $\tZtr\colon \Delta^\opp \Sc\to \Delta^\opp \PST[S]$. En
utilisant la correspondance de Dold-Kan
$\Delta^\opp\PST[S]\simeq\CompNeg(\PST[S])$
(cf.~\cite[\S2, Chapter~III]{goerss-jardine}), on en déduit un foncteur
$\tM\colon \Delta^\opp\Sc\to\CompNeg(\PST[S])$ : pour tout $\mathcal X\in
\Delta^\opp\Sc$, $\tM(\mathcal X)$ est le complexe
normalisé associé au préfaisceau avec transferts simplicial $\tZtr(\mathcal
X)$. On dispose d'une inclusion fonctorielle $\tM(\mathcal X)\subset
\tMnn(\mathcal X)$ où $\tMnn(\mathcal X)$ est le complexe dont les
différentielles sont des sommes alternées de cofaces :
\[\xymatrix@C=1.5cm{\tZtr(\mathcal X_0)&\ar[l]_-{d_0^\star-d_1^\star}
\tZtr(\mathcal X_1) &\ar[l]_-{d_0^\star-d_1^\star+d_2^\star}
\tZtr(\mathcal X_2)&\ar[l] \dots}\]

La notation $\tM$ est compatible avec celle de la
proposition~\ref{proposition-adjonction-non-a-1-localise}. On rappelle que
d'après \cite[Theorem~2.5~(3), Chapter~III]{goerss-jardine},
l'inclusion $\tM(\mathcal
X)\subset \tMnn(\mathcal X)$ est une équivalence d'homotopie fonctorielle.

Avec ces notations, le complexe $\tC(\mathcal X,K)$ n'est autre que le
complexe d'homomorphismes dans la catégorie abélienne $\PST[S]$ entre
$\tMnn(\mathcal X)$ et $K$. On dispose donc, pour tout $n\in\mathbf{Z}$,
d'une bijection canonique :
\[H^n(\tC(\mathcal X,K))\isomto \Hom_{K^-(\PST[S])}(\tMnn(\mathcal X),K[n])\]
où $K^-(\PST[S])$ est la catégorie à homotopie près des complexes
bornés supérieurement d'objets de $\PST[S]$.
Comme $\tMnn(\mathcal X)$ est formé d'objets projectifs, ce $\Hom$ dans
$K^-(\PST[S])$ s'identifie à un $\Hom$ dans la catégorie dérivée de
$\PST[S]$. Retenons que l'on a défini des bijections pour tout $n\geq 0$ :
\[H^n(\tC(\mathcal X,K)\isomto \Hom_{\DerNeg(\PST[S])}(\tMnn(\mathcal
X),K[n])\]

En utilisant le foncteur évident $\DerNeg(\PST[S])\to\DMNeg[S]$ et
l'adjonction entre $\DMNeg[S]$ et $\Hopt$
(cf.~\S\ref{subsection-adjonction}), on peut définir des applications
\[H^n(\tC(\mathcal X,K))\to \Hom_{\DMNeg[S]}(\tM(\mathcal X),K[n])\simeq
\Hom_{\Hopt}(\mathcal X, K(K[n]))\;\text{.}\]

Dans le cas particulier où $K=\tM(\Th_S E)\otimes_\mathbf{Z} \Lambda$
(placé en degré $0$) pour
un fibré vectoriel $E$ de rang $r$ sur $S$, le complexe $\tC(\mathcal X,K)$
sera noté $\tC(\mathcal X,K_{E,\Lambda})$, auquel cas, compte tenu de
l'isomorphisme de Thom relatif, la construction ci-dessus fournit une
application de la forme suivante pour tout $n\geq 0$ :
\[H^n(\tC(\mathcal X,K_{E,\Lambda}))\to \tH^{2r+n,r}(\mathcal
X,\Lambda)\;\text{.}\]
Pour $u$ un $n$-cocycle du complexe
$C(\mathcal X,K_{E,\Lambda})$, on notera $[u]$
la classe de cohomologie associée dans $\tH^{2r+n,r}(\mathcal X,\Lambda)$.

Le lemme suivant montre que toute classe de cohomologie motivique peut-être
ainsi représentée par un cocycle :

\begin{lemme}\label{lemme-representation-zero-cocycle}
Pour tout $\mathcal X\in\Hopt$ et $u\in \tH^{2r,r}(\mathcal X,\Lambda)$, il
existe $\mathcal X'\in\Delta^\opp\Sc$, un isomorphisme $\mathcal X\simeq
\mathcal X'$ dans $\Hopt$ et un $0$-cocycle $u'\in \tC^0(\mathcal
X,K_{r,\Lambda})$
tel que la classe $[u']\in \tH^{2r,r}(\mathcal X',\Lambda)$ corresponde à
$u$ \emph{via} l'isomorphisme $\mathcal X\simeq \mathcal X'$.
\end{lemme}

Tout d'abord, remarquons une petite évidence. Supposons que $u$ soit donné
par un morphisme de préfaisceaux simpliciaux pointés $v\colon \mathcal X\to
K_r$, alors l'action de $v$ en dimension simpliciale $0$ est décrite par un
élément $v_0\in K_r(\mathcal X_0)$ qui s'identifie à un $0$-cocycle de
$\tC(\mathcal X,K_r)$ qui est évidemment tel que $u=[v_0]$. Il s'agit ici de
se ramener à ce cas.

On sait déjà (cf. corollaire~\ref{corollaire-calcul-hom-hopt}) que quitte à
changer $\mathcal X$, on peut représenter le morphisme $\mathcal X\to K_r$
dans $\Hopt$ par un morphisme de préfaisceaux simpliciaux pointés $\mathcal
X\to \Sing K_r$. En utilisant une variante de l'adjonction
\cite[Proposition~3.14, p.~91]{morel-voevodsky}, ce morphisme correspond à
un morphisme de préfaiscaux $\left|\mathcal
X\right|_{\mathbf{\Delta}^\bullet}\to K_r$. Il suffit alors de choisir pour
$\mathcal X'$ une équivalence faible (simpliciale) $\mathcal X'\to
\left|\mathcal X\right|_{\mathbf{\Delta}^\bullet}$ avec $\mathcal X'\in
\Delta^\opp\Sc$.

\bigskip

Le lemme suivant permet de comprendre le calcul du Bockstein :

\begin{lemme}\label{lemme-calcul-bockstein}
Soit $0\to \mathcal F'\vers i \mathcal F\vers p
\mathcal F''\to 0$ une suite exacte
dans $\PST[S]$. Soit $\mathcal X\in\Delta^\opp\Sc$. Alors, on a une suite
exacte de complexes :
\[0\to \tC(\mathcal X,\mathcal F')\vers i \tC(\mathcal X,\mathcal F)\vers p
\tC(\mathcal X,\mathcal F'')\to 0\;\text{.}\]
Soit $u\in \tC^n(\mathcal X,\mathcal F'')$ un $n$-cocycle. Il existe
$\tilde{u}\in \tC^n(\mathcal X,\mathcal F)$ un élément tel que
$p(\tilde{u})=u$. L'élément $d\tilde{u}$ de $\tC^{n+1}(\mathcal X,\mathcal
F)$ est l'image d'un unique élément $v\in \tC^{n+1}(\mathcal X,\mathcal F')$
qui vérifie $dv=0$.
Alors, la classe de $v$ dans $H^{n+1}(\mathcal X,\mathcal F')$ ne dépend
que de la classe de $u$ dans $H^{n}(\mathcal X,\mathcal F'')$ et elle
s'identifie, au signe près, à la classe induite par la
composition suivante dans
$\DerNeg(\PST[S])$ :
\[\tMnn(\mathcal X)\vers{u} \mathcal F''[n] \vers {\beta[n]} \mathcal F'[n+1]\]
où $\beta\colon \mathcal F''\to \mathcal F'[1]$
correspond à l'extension $0\to \mathcal F'\to\mathcal
F\to\mathcal F''\to 0$.
\end{lemme}

L'exactitude vient de ce que $\tMnn(\mathcal X)$ soit formé d'objets
projectifs de $\PST[S]$. Le reste n'est qu'un rappel d'algèbre homologique.

\bigskip

Nous allons montrer le lemme suivant :

\begin{lemme}\label{lemme-image-transfert}
Soit $\mathcal X\in \Hopt$. Soit $u\in \tH^{2r,r}(\mathcal
X,\mathbf{Z}/\ell\mathbf{Z})$. Alors, $\beta (Pu)\in
\tH^{2r\ell+1,r\ell}(\mathcal X,\mathbf{Z}/\ell\mathbf{Z})$
appartient à l'image du
transfert associé au revêtement étale $U\to S$ :
\[\Tr\colon
\tH^{2r\ell+1,r\ell}(\mathcal X\wedge U_+,\mathbf{Z}/\ell\mathbf{Z})\to
\tH^{2r\ell+1,r\ell}(\mathcal X,\mathbf{Z}/\ell\mathbf{Z})\;\text{.}\]
\end{lemme}

D'après le lemme~\ref{lemme-representation-zero-cocycle}, on peut supposer
que $\mathcal X\in\Delta^\opp\Sc$ et que la classe de cohomologie est
donnée par un $0$-cocycle, c'est-à-dire un élément $u\in \tC^0(\mathcal
X,K_{r,\mathbf{Z}/\ell\mathbf{Z}})$ tel que $du=0$. L'élément $u$
s'interprète comme un élément de $K_{r,\mathbf{Z}/\ell\mathbf{Z}}(\mathcal
X_0)$. La classe $Pu\in \Hom_{\Hopt}(\mathcal
X,K_{\xi^r,\mathbf{Z}/\ell\mathbf{Z}})$ est représentée par le $0$-cocycle
$Pu\in \tC^0(\mathcal
X,K_{\xi^r,\mathbf{Z}/\ell\mathbf{Z}})=K_{\xi^r,\mathbf{Z}/\ell\mathbf{Z}}(\mathcal
X)$ obtenu en utilisant le morphisme de préfaisceaux pointés $P\colon
K_{r,\mathbf{Z}/\ell\mathbf{Z}}\to K_{\xi^r,\mathbf{Z}/\ell\mathbf{Z}}$. Le
lemme~\ref{lemme-calcul-bockstein} permet de représenter $\beta Pu$ par un
cocycle. Pour ce faire, on choisit $\tilde{u}\in \tC^0(\mathcal
X,K_{r,\mathbf{Z}/\ell^2\mathbf{Z}})$ relevant $u$. On note $P\tilde{u}\in
\tC^0(\mathcal X,K_{\xi^r,\mathbf{Z}/\ell^2\mathbf{Z}})$ l'image de
$\tilde{u}$ par le morphisme $P\colon K_{r,\mathbf{Z}/\ell^2\mathbf{Z}}\to
K_{\xi^r,\mathbf{Z}/\ell^2\mathbf{Z}}$ associé à l'anneau de coefficients
$\mathbf{Z}/\ell^2\mathbf{Z}$.

Bien sûr, $P\tilde{u}\in \tC^0(\mathcal
X,K_{\xi^r,\mathbf{Z}/\ell^2\mathbf{Z}})$ est un relèvement de $Pu\in
\tC^0(\mathcal X,K_{\xi^r,\mathbf{Z}/\ell\mathbf{Z}})$, mais il n'y a pas de
raison pour que ce soit un $0$-cocycle. Comme $Pu$ est un $0$-cocycle, il
existe un unique $v\in \tC^1(\mathcal X,K_{\xi^r,\mathbf{Z}/\ell\mathbf{Z}})$
tel que $dP\tilde{u}=\ell v$ dans $\tC^1(\mathcal
X,K_{\xi^r,\mathbf{Z}/\ell^2\mathbf{Z}})$. Alors, $v$ est un $1$-cocycle
tel que $[v]=\beta Pu$ dans $\tH^{2r\ell+1,r\ell}(\mathcal X,
\mathbf{Z}/\ell\mathbf{Z})$. Nous allons montrer que
ce $1$-cocycle est l'image par le transfert d'un $1$-cocycle représentant
une classe dans $\tH^{2r\ell+1,r\ell}(\mathcal X\wedge
U_+,\mathbf{Z}/\ell\mathbf{Z})$.

Soit $\mathcal F\in \PST[S]$. Pour tout $X\in\Sm$, on a une application
\[\pi_\star\colon \mathcal F(X\times_S U)\to \mathcal F(X)\]
induite par l'action de la correspondance finie évidente de $X$ dans
$X\times_S U$ donnée par le revêtement étale $X\times_S U\to X$ déduit de
$U\to S$ par changement de base.
Ceci définit un morphisme $\pi_\star\colon \pi_\star \pi^\star \mathcal
F\to\mathcal F$ dans $\PST[S]$. On dispose aussi d'un morphisme plus
évident
\[\pi^\star\colon \mathcal F\to \pi_\star \pi^\star \mathcal F\]
(c'est un morphisme d'adjonction).
Dans le cas où $\mathcal F=K_{\xi^r,\mathbf{Z}/\ell\mathbf{Z}}$ (ou
$K_{r\ell,\mathbf{Z}/\ell\mathbf{Z}}$ si on préfère), le morphisme
$\pi_\star\colon \pi_\star \pi^\star \mathcal F\to\mathcal F$ induit le
transfert dont il est question dans l'énoncé du lemme :
\[\tH^{2r\ell+\star,r\ell}(\mathcal X\wedge U_+,\mathbf{Z}/\ell\mathbf{Z})
\to \tH^{2r\ell+\star,r\ell}(\mathcal X,\mathbf{Z}/\ell\mathbf{Z})\;\text{.}\]
Pour établir le lemme, il va donc être suffisant de montrer que le
$1$-cocycle $v$ de $\tC(\mathcal X,K_{\xi^r,\mathbf{Z}/\ell\mathbf{Z}})$
représentant $\beta Pu$ est l'image d'un $1$-cocycle par le morphisme de
complexes induit par $\pi_\star$ :
\[\Tr\colon \tC(\mathcal X,\pi_\star K_{\pi^\star
\xi^r,\mathbf{Z}/\ell\mathbf{Z}})\to
\tC(\mathcal X,K_{\xi^r,\mathbf{Z}/\ell\mathbf{Z}})\;\text{.}\]
Observons que l'on dispose d'un isomorphisme canonique de fibrés vectoriels
$\pi^\star E^\ell\simeq \pi^\star \xi^r$ où $E$ est le fibré trivial de rang $r$
sur $S$. Nous noterons $\tw$ cet isomorphisme. Nous allons
en réalité montrer que $v$ est l'image d'un $1$-cocycle de $\tC(\mathcal
X,K_{E^\ell,\mathbf{Z}/\ell\mathbf{Z}})$ \emph{via} la composition suivante :
\[\xymatrix{K_{E^\ell,\mathbf{Z}/\ell\mathbf{Z}}\ar[r]^-{\pi^\star}&
\pi_\star K_{\pi^\star E^\ell,\mathbf{Z}/\ell\mathbf{Z}}\ar[r]^-{tw}_\sim&
\pi_\star K_{\pi^\star \xi^r,\mathbf{Z}/\ell\mathbf{Z}}\ar[r]^-{\pi_\star}&
K_{\xi^r,\mathbf{Z}/\ell\mathbf{Z}}}\;\text{.}\]

Notons $w$ l'unique élément de $K_{r,\mathbf{Z}/\ell\mathbf{Z}}(\mathcal
X_1)$ tel que $d\tilde{u}=d_0^\star \tilde{u}-d_1^\star\tilde{u}=\ell w\in
\tC^1(\mathcal X,K_{r,\mathbf{Z}/\ell^2\mathbf{Z}})$ (ceci est possible
puisque $d\tilde{u}$ est congru modulo $\ell$ à $du$ qui est nul).
Comme $dd\tilde{u}=0$, il vient que $\ell dw=0$ dans $\tC^2(\mathcal
X,K_{r,\mathbf{Z}/\ell^2\mathbf{Z}})$, d'où $dw=0$ dans $\tC^2(\mathcal
X,K_{r,\mathbf{Z}/\ell\mathbf{Z}})$. On a ainsi un $1$-cocycle $w$ du
complexe $\tC(\mathcal X,K_{r,\mathbf{Z}/\ell\mathbf{Z}})$.

Si $a_1,a_2,\dots,a_\ell$ sont des éléments de
$K_{E,\mathbf{Z}/\ell\mathbf{Z}}(\mathcal X_i)$ pour un certain $i\geq 0$,
on notera $a_1a_2\dots a_\ell\in
K_{E^\ell,\mathbf{Z}/\ell\mathbf{Z}}(\mathcal X_i)$ l'image du $\ell$-uplet
qu'ils forment par le morphisme de
préfaisceaux \guil{produit} $K_{E,\mathbf{Z}/\ell\mathbf{Z}}\times \dots
\times K_{E,\mathbf{Z}/\ell\mathbf{Z}} \to
K_{E^\ell,\mathbf{Z}/\ell\mathbf{Z}}$.

Notons $u'=d_0^\star u=d_1^\star u\in
K_{E,\mathbf{Z}/\ell\mathbf{Z}}(\mathcal X_1)$. On peut former le produit
$\gamma=u'^{\ell-1}w\in K_{E^\ell,\mathbf{Z}/\ell\mathbf{Z}}(\mathcal
X_1)$. J'affirme que $\gamma$ est un $1$-cocycle de $\tC(\mathcal
X,K_{E^\ell,\mathbf{Z}/\ell\mathbf{Z}})$. En effet, si on note
$u''=d_0^\star u'=d_1^\star u'=d_2^\star u'\in
K_{E^\ell,\mathbf{Z}/\ell\mathbf{Z}}(\mathcal X_2)$, on a :
\begin{eqnarray*}
d\gamma &=& d_0^\star \gamma-d_1^\star \gamma+d_2^\star\gamma\\
&=& u''^{\ell-1}(d_0^\star w-d_1^\star w+d_2^\star w) \\
&=& u''^{\ell-1}dw =0
\end{eqnarray*}

Par ailleurs, dans $\tC^1(\mathcal X,K_{\xi^r,\mathbf{Z}/\ell^2\mathbf{Z}})$,
on a 
$\ell v=dP\tilde{u}=d_0^\star P\tilde{u}-d_1^\star
P\tilde{u}=Pd_0^\star \tilde{u}-Pd_1^\star \tilde u$. Comme $d_0^\star
\tilde{u}=d_1^\star \tilde{u}+\ell w$, le lemme suivant permet d'affirmer
que :
\[\ell v = Pd_0^\star\tilde{u} -Pd_1^\star\tilde{u}
= P(d_1^\star \tilde{u}+\ell w)-P(d_1^\star \tilde{u})=
\ell \cdot \left(\frac{1}{(\ell-1)!}\pi_\star \tw \pi^\star
\gamma\right)\in \tC^1(\mathcal X,K_{\xi^r,\mathbf{Z}/\ell^2\mathbf{Z}})\;\text{.}\]
Ceci permet de conclure que $v=
\frac{1}{(\ell-1)!}\pi_\star \tw \pi^\star \gamma\in
\tC^1(\mathcal X,K_{\xi^r,\mathbf{Z}/\ell\mathbf{Z}})$. Pour achever la
démonstration du lemme~\ref{lemme-image-transfert}, il ne reste plus qu'à
établir le lemme suivant :

\begin{lemme}
Soit $X\in\Sm$. Soit $\tilde{a}\in K_{E,\mathbf{Z}/\ell^2\mathbf{Z}}(X)$.
Soit $b\in K_{E,\mathbf{Z}/\ell\mathbf{Z}}$. On note $a$ l'image de
$\tilde{a}$ dans $K_{E,\mathbf{Z}/\ell\mathbf{Z}}(X)$. Alors, dans
$K_{E,\mathbf{Z}/\ell^2\mathbf{Z}}(X)$, on a l'égalité :
\[P(\tilde{a}+\ell b)-P(\tilde{a})=\ell\cdot
\left(\frac{1}{(\ell-1)!}\pi_\star \tw \pi^\star (a^{\ell-1}b)\right)\]
\end{lemme}

On note encore $\tilde{a}$ un représentant de $\tilde{a}$ dans
$\cequi(X\times_S E/X,0)\otimes \mathbf{Z}/\ell^2\mathbf{Z}$ et de même pour
$b\in \cequi(X\times_S E/X,0)\otimes \mathbf{Z}/\ell\mathbf{Z}$.
Par construction, $P(\tilde{a}+\ell b)\in \cequi(X\times_S (E\otimes
\xi)/X,0)\otimes \mathbf{Z}/\ell^2\mathbf{Z}$ est l'unique élément qui
vérifie
\[\pi_\xi^\star P(\tilde{a}+\ell b)=\tw \pi^\star ((\tilde{a}+\ell b)^\ell)=
\pi_\xi^\star P(\tilde{a})+\ell\tw \pi^\star\left(\sum_{i=0}^{\ell-1}a^i b
a^{\ell-1-i}\right) \]
Pour lever toute ambiguité, dans cette identité interviennent deux
morphismes évidents 
$\pi\colon X\times_S E^\ell \times_S U\to X\times_S
E^\ell$ et $\pi_\xi\colon X\times_S \xi^r\times_S U\to X\times_S \xi^r$
qui étaient jusqu'ici tous les deux notés $\pi$.

Pour conclure, on utilise l'identité suivante
dans $\cequi(X\times_S (E\otimes \xi)/X,0)\otimes
\mathbf{Z}/\ell\mathbf{Z}$ qui résulte du calcul standard de la composition
$\pi^\star_\xi\pi_{\xi\star}$ au niveau des cycles pour une projection
$\pi_\xi$ d'un $\Sym \ell$-torseur étale :
\[\pi^\star_\xi \pi_{\xi\star} (\tw \pi^\star (a^{\ell-1}b))=
(\ell-1)!\cdot \tw \pi^\star\left(
\sum_{i=0}^{\ell-1} a^iba^{l-1-i}\right)\;\text{.}\]

\bigskip

Nous sommes maintenant en mesure de démontrer le
théorème~\ref{theoreme-annulation-bockstein}. D'après les considérations du
\S\ref{subsection-bet-g} (notamment le
lemme~\ref{lemme-morphisme-g-torseur-dans-u-d}), on peut supposer que $S$
est l'un des schémas $S_d$ définis au début du \S\ref{subsection-bet-g}
(ils forment un système inductif dont la colimite s'identifie à $\Bet
\Sym\ell$) et que $U$ est le $\Sym \ell$-torseur $U_d$ sur $S_d$ défini au
même endroit.

Notons $P_{\ell,d}\colon \tH^{2r,r}(\mathcal X)\to
\tH^{2r\ell,r\ell}(\mathcal X\wedge S_{d+})$ les différentes opérations
$P_\ell$ associées à ces $\Sym \ell$-torseurs $U_d$ sur $S_d$ pour $d\geq
1$. Il s'agit de montrer que $\beta P_{\ell,d}(x)=0$ pour tout $d\geq 1$.
On sait d'une part que les classes $\beta P_{\ell,d}(x)$ forment un système
compatible pour $d\geq 1$ et d'autre part que pour tout $e\geq 1$, il
existe $\tilde{x}_e\in \tH^{2r\ell,r\ell}(\mathcal X\wedge U_{e+})$ tel que
$\beta P_{\ell,e}(x)$ soit l'image de $\tilde{x}_e$ par le transfert (c'est
le lemme~\ref{lemme-image-transfert}).
Ainsi, pour $e\geq d\geq 1$, il vient que la classe $\beta P_{\ell,d}(x)$
est dans l'image de
\[\tH^{\star,\star}(\mathcal X\wedge U_{e+})\to \tH^{\star,\star}(\mathcal
X\wedge U_{d+})\vers{\Tr} \tH^{\star,\star}(\mathcal X\wedge
S_{d+})\;\text{.}\]
Fixons $d\geq 1$. J'affirme que cette composition est nulle pour $e$ assez
grand. Cela résulte formellement de l'annulation de la composition suivante
dans $\DMmoins[k,\Fl]$ pour $e$ assez grand :
\[M(S_d)\vers{Tr}M(U_d)\to M(U_e)\;\text{.}\]
Justifions ceci. $M(S_d)$ étant un objet \guil{compact} de
$\DMmoins[k,\Fl]$, pour montrer l'annulation de la composition pour $e$
assez grand, il suffit de montrer que la composition s'annule si on
remplace $M(U_e)$ par la colimite de ce système pour $e\geq d$, qui est le
motif unité $\Fl$ (cf. proposition~\ref{proposition-bgm}). Que la composition
\[M(S_d)\vers{\Tr}M(U_d)\to \Fl\]
soit nulle provient du fait que la classe $1\in
H^{0,0}(U_d)$ s'envoie dans $H^{0,0}(S_d)$ par le transfert sur
la classe modulo $\ell$ du degré $U_d\to S_d$ qui vaut $l!$.

\subsection{Comparaison avec la construction de Voevodsky}
\label{subsection-comparaison-voevodsky}

Nous allons expliquer brièvement ici pourquoi la construction de Voevodsky
\cite[Construction~5.3]{voevodsky-reduced} donne la même opération
cohomologique que celle de la définition~\ref{definition-p-sur-k}.

Les données de la définitions~\ref{definition-p-sur-k} sont celles d'un
groupe fini agissant sur un ensemble fini $A$ de cardinal $n$ et d'un
schéma $S\in \Sm[k]$ muni d'un $G$-torseur étale $U$. Voevodsky se donne en
outre un fibré vectoriel $L$ et un isomorphisme de fibrés vectoriels sur
$S$ entre $\xi\oplus L$ et un fibré trivial de rang $N$ où $\xi$ est le
fibré vectoriel de la définition~\ref{definition-xi}. Pour tout $r\geq 0$,
Voevodsky définit un morphisme
\[\tilde{P}\colon K_i\wedge \Th_{G\backslash U} L^i\to K_{iN}\;\text{.}\]
Par un jeu d'adjonctions qu'il serait peu intéressant de détailler, ce
morphisme dans la catégorie des préfaisceaux d'ensembles pointés sur
$\Sm[k]$ correspond à un morphisme dans la catégorie des préfaisceaux
d'ensembles pointés sur $\Sm[S]$ :
\[\tilde{P}\colon K_i\to \SheafHom_\bullet(\Th_S L^i,K_{iN})\]
où $K_i$, $K_{iN}$ et $\Th_S L^i$ sont considérés ici étant comme
étant \guil{au-dessus de $S$}. Il résulte aussitôt des définitions que ce
morphisme est égal au composé :
\[K_i\vers{P_{\mathbf{A}^i}} K_{\xi^i}\vers{\cdot \tau_{L^i}}
\SheafHom_\bullet(\Th_S L^i,K_{\xi^i\oplus L^i})\simeq
\SheafHom_\bullet(\Th_S L^i,K_{iN})\]
où $P_{\mathbf{A}^i}$ est le morphisme de la
définition~\ref{definition-operation-totale} appliqué au fibré trivial de
rang $i$ sur $S$ et $\cdot \tau_{L^i}$ le morphisme de multiplication par
la classe tautologique (cf.
proposition~\ref{proposition-classes-tautologiques}).

Il est dès lors évident que la transformation
\[\tilde{P}\colon \tH^{2i,i}(\mathcal X)\to \tH^{in,2in}(\mathcal X\wedge
S_+)\]
que Voevodsky en a déduit pour tout $\mathcal X\in\Hopt[k]$ est la même que
celle de la définition~\ref{definition-p-sur-k}. Ceci montre \emph{a
posteriori} que bien que la démonstration de ce fait
fût incorrecte, l'énoncé
d'indépendance en le choix du supplémentaire $L$ de $\xi$
\cite[Lemmas~5.2 \& 5.3]{voevodsky-reduced} était vrai.

\section{Opérations de Steenrod}
\label{section-operations}

Soit $\ell$ un nombre premier. Soit $k$ un corps parfait de caractéristique
différente de $\ell$. Nous allons enfin voir comment les constructions
précédentes s'appliquent au cas de l'action évidente du groupe symétrique
$\Sym \ell$ agissant sur $\{1,\dots,\ell\}$ et donnent lieu aux opérations de
Steenrod motiviques sur la cohomologie modulo $\ell$. Il faut pour cela
commencer par décomposer le motif à coefficients
$\mathbf{Z}/\ell\mathbf{Z}$ du classifiant $\Bet\Sym \ell$.
Dans cette section, s'il n'est pas précisé, l'anneau des coefficients des
groupes de cohomologie sera toujours $\Fl$.

\subsection{Le motif de $\Bet \Sym\ell$}
\label{subsection-motif-b-sym-ell}

\begin{proposition}\label{proposition-decomposition-motif-de-b-mu-ell}
Le motif $M(\Bet \mu_\ell;\Fl)\in \DMNeg[k,\Fl]$ se décompose de la façon
suivante :
\[M(\Bet \mu_\ell;\Fl)\simeq \oplus_{i\geq 0}\Fl(i)[2i]\bigoplus\oplus_{i\geq
0}\Fl(i+1)[2i+1]\;\text{.}\]
Les projections sur les différents facteurs correspondent aux classes de
cohomologie $v^i$ et $uv^i$ dans $H^{\star,\star}(\Bet \mu_\ell)$ où $v\in
H^{2,1}(\Bet \mu_\ell)$ est l'image inverse par le morphisme évident
$\Bet\mu_\ell\to \Bet\Gm\simeq \mathbf{P}^\infty$ de $c_1(\OO(1))\in
H^{2,1}(\mathbf{P}^\infty)\simeq \lim_n H^{2,1}(\mathbf{P}^n)$ et où $u\in
H^{1,1}(\Bet\mu_\ell)$ est l'unique élément tel que $\beta u=v$ et dont la
restriction au point-base de $\Bet \mu_\ell$ soit nulle dans $H^{1,1}(k)$.
\end{proposition}

On définit $\Bgm \mu_\ell$ en utilisant le caractère évident
$\mu_\ell\subset \Gm$ et les ouverts $U_d=\mathbf{A}^d-\{0\}$ les plus
grands possibles (cf.~proposition~\ref{proposition-bgm}) pour obtenir
l'identification $\Bgm\mu_\ell\simeq \Bet \mu_\ell$. La proposition résulte
alors de la décomposition, pour tout $d\geq 1$, du motif du schéma quotient
$\mu_\ell\backslash (\mathbf{A}^d-\{0\})$ faisant intervenir ceux des
termes apparaissant ci-dessus qui ont un indice $i$ vérifiant $0\leq i\leq
d-1$.

Cette décomposition provient d'une formule plus générale dans la situation
d'un $\Gm$-torseur $T$ au-dessus de $X\in\Sm[k]$ (comme ici
$\mathbf{A}^d-\{0\}$ est un $\Gm$-torseur au-dessus de $\mathbf{P}^{d-1}$).
Il existe alors à isomorphisme unique près
un unique fibré en droites $L$ sur $X$
muni d'un isomorphisme de $\Gm$-torseurs entre $T$ et le complémentaire de
la section nulle de $L$. À partir de la suite de Kummer $0\to \mu_\ell \to
\Gm \vers {x\mapsto x^\ell} \Gm\to 0$ (exacte pour la topologie étale), il
vient immédiatement que le quotient $\mu_\ell\backslash T$ s'identifie au
complémentaire de la section nulle dans le fibré en droites $L^{\otimes
\ell}$. Le $\Gm$-torseur $\mathbf{A}^d-\{0\}$ étant le complémentaire de la
section nulle dans $\OO(-1)$, il vient ainsi que $\mu_\ell\backslash
(\mathbf{A}^d-\{0\})$ s'identifie au complémentaire de la section nulle
dans le fibré en droites $\OO(-\ell)$ sur $\mathbf{P}^{d-1}$.

Dans la situation générale d'un $\Gm$-torseur $T$ sur $X\in \Sm[k]$ qui
soit le complémentaire de la section nulle d'un fibré en droites $L$ sur
$X$, en utilisant l'isomorphisme $M(\Th_X L)\simeq M(X)(1)[2]$ donné par
l'isomorphisme de Thom (cf.~définition~\ref{definition-classe-thom}),
on voit que le motif de $T$ s'insère dans un triangle distingué dans
$\DMmoins[k,\Fl]$ : \[M(X)(1)[1]\vers\delta M(T)\to M(X)\to M(X)(1)[2]\] où
le morphisme $M(X)\to M(X)(1)[2]$ est induit par la classe $c_1(L)\in
H^{2,1}(X)$\;\footnote{Si on dispose implicitement d'un morphisme $p\colon
Y\to X$ dans $\Sm[k]$, une classe $x\in H^{p,q}(Y)$ (correspondant à un
morphisme $M(Y)\to \Fl(q)[p]$) induit un morphisme $M(Y)\to M(X)(q)[p]$ qui
est la composition $M(Y)\vers{(\Id,p)} M(Y\times X)\simeq M(Y)\otimes
M(X)\vers{x\otimes \Id_{M(X)}} M(X)(q)[p]$.}. Si cette classe est nulle, il
existe une rétraction $M(T)\to M(X)(1)[1]$ de $\delta$ et on obtient une
décomposition $M(T)\isomto M(X)\oplus M(X)(1)[1]$. Plus précisément, si
$w\in H^{1,1}(T)$ vérifie $\delta^\star w=1\in H^{0,0}(X)$, la classe $w$
et la projection $T\to X$ induisent un morphisme $M(T)\to M(X)(1)[1]$.
C'est une rétraction de $\delta$. En effet, il s'agit de montrer que le
morphisme composé
\[M(X)(1)[1]\vers{\delta} M(T)\vers w M(X)(1)[1]\]
est l'identité. Par construction, cette composition fixe la classe $1\in
H^{0,0}(X)$ et on peut conclure en utilisant le fait évident que $\delta\colon
M(X)(1)[1]\to M(T)$ (et donc cette composition) est un morphisme de
$M(X)$-comodules où $M(X)(1)[1]$ et $M(T)$ sont munis des structures
évidentes de $M(X)$-comodules.
Ainsi, les classes $1$ et $w$ induisent une décomposition 
$M(T)\isomto M(X)\oplus M(X)(1)[1]$ et on peut observer que la classe $w\in
H^{1,1}(T)$ est bien définie modulo l'image de $H^{1,1}(X)$.

Comme $c_1(\OO(-\ell))=\ell c_1(\OO(-1))=0\in H^{2,1}(\mathbf{P}^{d-1})$,
on obtient l'existence d'une classe $w\in H^{1,1}(\mu_\ell\backslash
(\mathbf{A}^d-\{0\}))$ telle que $1$ et $w$ induisent un isomorphisme
$M(\mu_\ell\backslash (\mathbf{A}^d-\{0\}))\isomto
M(\mathbf{P}^{d-1})\oplus M(\mathbf{P}^{d-1})(1)[1]$. On peut assurer
l'unicité de $w$ (vérifiant $\delta^\star w=1$ comme plus haut) en
demandant en outre que la restriction de $w$ au point-base de
$\mu_\ell\backslash (\mathbf{A}^d-\{0\})$ provenant d'un élément quelconque
de $k^d-\{0\}$ soit nulle (on laisse en exercice au lecteur de montrer que
les morphismes de restriction $H^{1,1}(\mu_\ell\backslash
\mathbf{A}^d-\{0\})\to H^{1,1}(k)$
associés aux différents éléments de $k^d-\{0\}$ sont tous égaux).

Si on avait $\beta w=0$, la classe $w$ se relèverait en une classe de
cohomologie $\tilde{w}\in H^{1,1}(\mu_\ell\backslash (\mathbf{A}^d-\{0\}))$
qui vérifierait $\delta^\star \tilde{w}\equiv 1 \mod \ell$.
L'exactitude de la suite
\[\xymatrix@C=1.5cm{H^{1,1}(\mu_\ell\backslash (\mathbf{A}^d-\{0\}),
\mathbf{Z}/\ell^2\mathbf{Z})\ar[r]^-{\delta^\star} &
H^{0,0}(\mathbf{P}^{d-1},\mathbf{Z}/\ell^2\mathbf{Z})\ar[r]^-{\cdot
c_1(\OO(-\ell))} & H^{2,1}(\mathbf{P}^{d-1},\mathbf{Z}/\ell^2\mathbf{Z})}\]
impliquerait alors que $c_1(\OO(-\ell))=0\in
H^{2,1}(\mathbf{P}^{d-1},\mathbf{Z}/\ell^2\mathbf{Z})$ ce
qui est faux pour $d\geq 2$.

Pour $d\geq 2$, on a donc $\beta w\neq 0$, ce qui implique l'existence de
$\lambda\in\Fl^\times$ tel que $\lambda\cdot \beta w=v\in
H^{2,1}(\mu_\ell\backslash (\mathbf{A}^d-\{0\}))$. En notant $u=\lambda
w\in H^{1,1}(\mu_\ell\backslash (\mathbf{A}^d-\{0\}))$ qui vérifie $\beta
u=v$, on obtient la décomposition sous la forme énoncée dans la proposition.

\bigskip

\begin{proposition}\label{proposition-decomposition-motif-coinvariants-bmu}
Pour tout $\lambda\in\Fl^\times$, notons $m_\lambda\colon
\mu_\ell\to\mu_\ell$ l'élévation à la puissance $\lambda$ et le morphisme
$\Bet\mu_\ell\to\Bet\mu_\ell$ qu'elle induit. Alors, dans
$H^{\star,\star}(\Bet\mu_\ell)$, on a : $m_\lambda^\star
u=\lambda u$ et $m_\lambda^\star v=\lambda v$.
Ainsi, \emph{via} la décomposition de la
proposition~\ref{proposition-decomposition-motif-de-b-mu-ell}, $m^\star$
agit par la multiplication par $\lambda^i$ sur le terme $\Fl(i)[2i]$ et par
la multiplication par $\lambda^{i+1}$ sur le terme $\Fl(i+1)[2i+1]$.

Enfin, la projection $M(\Bet\mu_\ell)\to M(\Bet \mu_\ell)_{\Fl^\times}$ sur
le facteur direct
des co-invariants sous l'action de $\Fl^\times$ est la projection sur la
sous-somme suivante :
\[M(\Bet \mu_\ell)_{\Fl^\times}\simeq \oplus_{j\geq
0}\Fl((\ell-1)j)[2(\ell-1)j]\bigoplus
\oplus_{j\geq 0} \Fl((\ell-1)(j+1))[2(\ell-1)(j+1)-1]\;\text{.}\]

Les classes $d=-v^{\ell-1}\in H^{2\ell-2,\ell-1}(\Bet \mu_\ell)$ et
$c=-uv^{\ell-2}\in H^{2\ell-3,\ell-1}(\Bet\mu_\ell)$ sont invariantes sous
l'action de $\Fl^\times$. Les projections de $M(\Bet
\mu_\ell)_{\Fl^\times}$ sur les différents termes de la décomposition
ci-dessus correspondent aux classes $(-1)^jd^j$ et $(-1)^{j+1}cd^j$.
\end{proposition}

La formule concernant $m_\lambda^\star v=\lambda v$ est évidente. Compte
tenu des propriétés de $u$, l'identité $m_\lambda^\star u=\lambda u$ en
résulte. Les autres énoncés en découlent immédiatement.

\begin{proposition}\label{proposition-isomorphisme-coinvariants-mu-ell-z-sur-ell}
Il existe un isomorphisme canonique dans $\DMNeg[k,\Fl]$ compatible à
l'extension du corps parfait $k$ :
\[M(\Bet \mathbf{Z}/\ell\mathbf{Z})_{\Fl^\times}\simeq
M(\Bet \mu_\ell)_{\Fl^\times}\;\text{.}\]
\end{proposition}

S'il existe $\zeta\in k^\times$ une racine primitive $\ell$-ième de
l'unité, l'isomorphisme est induit par l'isomorphisme
$\mathbf{Z}/\ell\mathbf{Z}\isomto \mu_\ell$ envoyant $1$ sur $\zeta$.
L'isomorphisme induit
au niveau des co-invariants est évidemment indépendant du
choix de $\zeta$. En général, il existe une extension finie galoisienne
$k'/k$ de degré divisant $\ell-1$ et donc premier à $\ell$ tel que $k'$
contienne une racine primitive $\ell$-ième de l'unité. Par un argument de
transfert, on obtient l'existence et l'unicité d'un isomorphisme de la
forme cherchée qui après extension des scalaires à $k'$ soit celui
construit précédemment dans le cas d'un corps possédant une racine
primitive $\ell$-ième de l'unité.

\medskip

On peut ainsi procéder à une identification $H^{\star,\star}(\Bet
\mathbf{Z}/\ell\mathbf{Z})^{\Fl^\times}\simeq H^{\star,\star}(\Bet
\mu_\ell)^{\Fl^\times}$ des sous-algèbres fixées par l'action de
$\Fl^\times$, ce qui fournit des classes $c$ et $d$ dans la
cohomologie de $\Bet \mathbf{Z}/\ell\mathbf{Z}$.

\begin{proposition}\label{proposition-motif-sym-ell}
Le morphisme canonique 
\[M(\Bet\mathbf{Z}/\ell\mathbf{Z})_{\Fl^\times}\to M(\Bet \Sym \ell)\]
est un isomorphisme dans $\DMNeg[k,\Fl]$.
\end{proposition}

On identifie ici implicitement $\mathbf{Z}/\ell\mathbf{Z}$ au sous-groupe
de $\Sym \ell$ engendré par le $\ell$-cycle de $\Sym \ell$ correspondant à
la permutation $x\longmapsto x+1$ de $\mathbf{Z}/\ell\mathbf{Z}$ que l'on
identifie ensemblistement à $\{1,\dots,\ell\}$ (par l'inverse de la
réduction modulo $\ell$). \emph{Via} ces identifications, le normalisateur
$N$ de $\mathbf{Z}/\ell\mathbf{Z}$ dans $\Sym \ell$ est le groupe affine
des bijections $x\longmapsto ax+b$ sur $\mathbf{Z}\ell\mathbf{Z}$ avec
$a\in\Fl^\times $ et $b\in \Fl^\times$. L'action de $N$ sur
$\mathbf{Z}/\ell\mathbf{Z}$ par \guil{conjugaison dans $\Sym \ell$} se
factorise par le quotient évident $N\to \Fl^\times$. Le morphisme considéré
dans la proposition n'est donc autre que l'épimorphisme scindé donné par le
corollaire~\ref{corollaire-transfert-motif-classifiants}.

Pour obtenir le résultat, il reste à montrer que la composition \[M(\Bet
\mathbf{Z}/\ell\mathbf{Z})_{\Fl^\times}\to M(\Bet \Sym\ell)\vers{\Tr}
M(\Bet \mathbf{Z}/\ell)_{\Fl^\times}\] est un isomorphisme. Notons
$\varphi\in\End(M(\Bet\mathbf{Z}/\ell)_{\Fl^\times})$ cette composition
divisée par $(\ell-2)!$. C'est un projecteur dont l'image s'identifie à
$M(\Bet \Sym \ell)$. Pour conclure, il faut montrer que $\varphi$ est
l'identité de $M(\Bet \mathbf{Z}/\ell)_{\Fl^\times}$. Par ailleurs, on peut
identifier $H^{\star,\star}(\Bet \Sym \ell)$ à une sous-algèbre de
$H^{\star,\star}(\Bet \mathbf{Z}/\ell)^{\Fl^\times}$. Cette sous-algèbre
est précisément l'image du projecteur que $\varphi$ induit. Avant de
montrer que $\varphi$ est l'identité, on va montrer que $\varphi$ induit
l'identité de $H^{\star,\star}(\Bet \mathbf{Z}/\ell)^{\Fl^\times}$, ce qui
revient à montrer le lemme suivant :

\begin{lemme}
L'inclusion $H^{\star,\star}(\Bet \Sym \ell)\subset 
H^{\star,\star}(\Bet \mathbf{Z}/\ell)^{\Fl^\times}$ est une égalité.
\end{lemme}

Il s'agit de montrer que les classes $c$ et $d$ appartiennent à la
sous-algèbre $H^{\star,\star}(\Bet \Sym \ell)$.

Le classifiant $\Bet \Sym \ell$ s'identifie à la colimite $\Bgm \Sym
\ell=\colim_i \mu_\ell\backslash U_i$ d'un système de faisceaux
représentables (cf. proposition~\ref{proposition-bgm}). Sur chacun de ces
quotients $\mu_\ell\backslash U_i$, on dispose d'après la
définition~\ref{definition-xi} d'un fibré vectoriel $\xi$ de rang $\ell$
associé à l'action tautologique de $\Sym \ell$ sur $\{1,\dots,\ell\}$. Ces
fibrés étant compatibles pour les différentes valeurs de $i$, on peut
donner un sens à $c_{\ell-1}(\xi)\in H^{2\ell-2,\ell-1}(\Bet \Sym \ell)$
(ceci sera rendu plus précis dans l'isomorphisme du
corollaire~\ref{corollaire-cohomologie-b-sym-ell}).
On peut vérifier que $d=c_{\ell-1}(\xi)$. En effet, on peut supposer qu'il
existe une racine primitive $\ell$-ième de l'unité $\zeta\in k$, ce qui
permet de définir une inclusion $\mu_\ell\to \Sym \ell$. On vérifie alors
facilement que l'image inverse de $\xi$ sur $\mu_\ell\backslash
(\mathbf{A}^d-\{0\})$ est $\oplus_{i=0}^{\ell-1}L^{\otimes i}$ où $L$ est
l'image inverse de $\OO(1)$ sur $\mathbf{P}^{d-1}$, ce qui permet d'obtenir
la formule : \[c_{\ell-1}(\xi)_{|\mu_\ell\backslash
\mathbf{A}^{d-0}}=\prod_{i=1}^{\ell-1} (iv) = -v^{\ell-1}=d\;\text{.}\]

Il reste à montrer que $c$ appartient aussi à $H^{\star,\star}(\Bet \Sym
\ell)$. Pour tout $\mathcal X\in\Ho[k]$, l'image du Bockstein $\beta$ sur
$H^{\star,\star}(\mathcal X;\mathbf{Z}/\ell\mathbf{Z})$ est le noyau de la
\guil{multiplication par $\ell$}
$H^{\star,\star}(\mathcal{X};\mathbf{Z}/\ell\mathbf{Z})\to
H^{\star,\star}(\mathcal{X};\mathbf{Z}/\ell^2\mathbf{Z})$.
Comme le morphisme évident $H^{\star,\star}(\Bet \Sym \ell,\Lambda)\to \
H^{\star,\star}(\Bet \mathbf{Z}/\ell,\Lambda)$ est injectif non seulement pour
$\Lambda=\mathbf{Z}/\ell\mathbf{Z}$ mais aussi pour
$\Lambda=\mathbf{Z}/\ell^2\mathbf{Z}$ (cf.
corollaire~\ref{corollaire-transfert-motif-classifiants}), le fait
que $d$ soit dans l'image du Bockstein sur $H^{\star,\star}(\Bet
\mathbf{Z}/\ell\mathbf{Z})$ implique que $d$ peut aussi s'écrire $d=\beta
(c')$ avec $c' \in H^{2\ell-3,\ell-1}(\Bet \Sym \ell)$. La
décomposition du motif $M(\Bet \mathbf{Z}/\ell\mathbf{Z})_{\Fl^\times}$
montre qu'il existe $a\in H^{2\ell-3,\ell-1}(k)$ et $b\in \Fl$ tels que
dans $H^{\star,\star}(\Bet \mathbf{Z}/\ell\mathbf{Z})^{\Fl^\times}$, on ait
$c' = a + b c$. L'identité $\beta (c')=d$ impose que $b=1$, d'où $c'=a+c$.
Ceci montre que $c=c'-a$ appartient bien à la sous-algèbre
$H^{\star,\star}(\Bet\Sym \ell)$ de
$H^{\star,\star}(\mathbf{Z}/\ell\mathbf{Z})^{\Fl^\times}$.

En utilisant que les classes $d^j$ et $cd^j$ dans $H^{\star,\star}(\Bet
\mathbf{Z}/\ell)^{\Fl^\times}$ pour tous $j\geq 0$ appartiennent à la
sous-algèbre $H^{\star,\star}(\Bet \Sym \ell)$, on obtient que si $p$ est
le morphisme de projection de $M(\Bet \mathbf{Z}/\ell)_{\Fl^\times}$ vers
un des termes $\Fl(i)[2i]$ ou $\Fl(i)[2i-1]$ apparaissant dans sa
décomposition (cf.
propositions~\ref{proposition-decomposition-motif-coinvariants-bmu} et
\ref{proposition-isomorphisme-coinvariants-mu-ell-z-sur-ell}), alors
$p\circ \varphi=\varphi$. Le lemme suivant appliqué à $\psi=\varphi-\Id$
permet de conclure que $\varphi$ est l'identité de $M(\Bet
\mathbf{Z}/\ell)_{\Fl^\times}$ et donc de terminer la démonstration de la
proposition~\ref{proposition-motif-sym-ell} :

\begin{lemme}
Soit $M=\oplus_{i\in I} M_i$ une décomposition en somme directe dans
$\DMNeg[k]$ telle que les objets $M_i$ appartiennent à $\DMeffgm[k]$.
Pour tout $i\in I$, on note $p_i\colon M\to M_i$ le morphisme de
projection.
Soit $\psi$ un endomorphisme de $M$ tel que pour tout $i\in I$, $p_i\circ
\psi=0$. Alors, $\psi=0$.
\end{lemme}

Pour tout $i\in I$, on note $j_i\colon M_i\to M$ le morphisme d'inclusion.
Par définition de la somme directe, pour conclure, il s'agit de montrer que
pour tout $i\in I$, $\psi\circ j_i=0$. Le morphisme $\psi\circ j_i$ étant
de la forme $M_i\to \oplus_{j\in I} M_j$, comme $M_i$ est compact, ce
morphisme $\psi\circ j_i$ se factorise par un facteur direct évident
$\oplus_{j\in J_i} M_j$ avec $J_i$ une partie finie de $I$. On est alors
ramené à montrer que pour tous $i\in I$ et $j\in J_i$, on a $p_j \circ
\psi\circ j_i=0$, ce qui résulte bien de l'hypothèse $p_j\circ \psi=0$.

\begin{corollaire}\label{corollaire-cohomologie-b-sym-ell}
Soit $V$ un $k$-espace vectoriel de dimension finie et $\rho\colon
\Sym \ell\to \GL(V)$ une représentation linéaire fidèle de $\Sym \ell$.
Soit $(U_i)_{i\geq 1}$ une suite d'ouverts des espaces affines $V^{\oplus
i}$ satisfaisant les hypothèses de la proposition~\ref{proposition-bgm}.
Alors, pour tout $\mathcal X\in \Hopt$, l'application canonique est un
isomorphisme
\[\tH^{\star,\star}(\mathcal X\wedge (\Bet \Sym \ell)_+)\isomto
\lim_{i\geq 1}\tH^{\star,\star}(\mathcal X\wedge (\Sym \ell\backslash
U_i)_+)\;\text{.}\]
\end{corollaire}

Comme la suite exacte de Milnor contient la surjectivité de cette
application, il suffit de montrer l'injectivité. Soit $x\in
\tH^{p,q}(\mathcal X\wedge (\Bet\Sym \ell)_+)$ une classe dont les images
dans les groupes $\tH^{p,q}(\mathcal X\wedge (\Sym \ell\backslash U_i)+)$
soient toutes nulles. Il s'agit de montrer que $x$ est nulle. Pour cela, on
peut supposer que $k$ contient $\zeta$, une racine primitive $\ell$-ième
de l'unité.
D'après la proposition~\ref{proposition-motif-sym-ell},
il suffit de montrer l'annulation de l'image inverse de $x$ par $\mathcal
X\wedge f_+$ où $f$ est le morphisme canoniquement associé à $\zeta$ :
$\Bet \mu_\ell\simeq \Bet \mathbf{Z}/\ell\mathbf{Z}\to \Bet \Sym \ell$.
La proposition~\ref{proposition-decomposition-motif-de-b-mu-ell}
fournit un isomorphisme
\[\tH^{p,q}(\mathcal X \wedge (\Bet \mu_\ell)_+)\isomto
\lim_{d\geq 1} \tH^{p,q}(\mathcal X\wedge (\mu_\ell\backslash
(\mathbf{A}^d-\{0\}))_+)\;\text{.}\]
Ainsi, pour montrer que $x=0$, il suffit de montrer
l'annulation de l'image inverse de $x$ par $\mathcal X\wedge {g_d}_+$ pour
tout $d\geq 1$ où $g_d$ est la composition $\mu_\ell\backslash
(\mathbf{A}^d-\{0\})\to \Bet \mu_\ell\simeq \Bet
\mathbf{Z}/\ell\mathbf{Z}\to \Bet \Sym \ell$. Pour obtenir cela pour
un certain $d$, vu l'hypothèse sur $x$, il suffit de montrer que
pour $i$ assez grand,
$g_d$ se factorise dans $\Ho[k]$ par le morphisme canonique
$\Sym\ell\backslash U_i\to \Bet \Sym \ell$, ce qui est évident.

\bigskip

La proposition suivante permet de comprendre la structure multiplicative
sur les algèbres de cohomologie motivique des classifiants de $\mu_\ell$ et
de $\Sym \ell$ :

\begin{proposition}[{\cite[Theorem~6.10]{voevodsky-reduced}}]
\label{proposition-formule-u-carre}
Dans $H^{\star,\star}(\Bet \mu_\ell)$, on a l'identité $u^2=0$ si $\ell\neq
2$ et si $\ell=2$, $u^2=\tau v + \rho u$ où $\tau\in
H^{0,1}(k)=\mu_2(k)=\{\pm 1\}$ est l'élément non nul et $\rho\in
H^{1,1}(k)\simeq k^\times/k^{\times 2}$ est la classe de $-1$.
Dans $H^{\star,\star}(\Bet \Sym \ell)$, on a aussi $c^2=0$ si $\ell\neq 2$
et $c^2=\tau d+\rho c$ si $\ell=2$.
\end{proposition}

\subsection{Construction des opérations $P^i$ et $B^i$ sur $\tH^{2\star,\star}$}
\label{subsection-construction-operation-p-i-b-i}

Pour tout $r\geq 0$, on note 
\[P_\ell\colon \tH^{2r,r}(\mathcal X)\to \tH^{2r\ell,r\ell}(\mathcal X
\wedge (\Bet \Sym \ell)_+)\]
la transformation naturelle, pour $\mathcal X\in\Hopt[k]$, correspondant à
la construction $P^{\univ}_{\Sym \ell}$ (cf.
définition~\ref{definition-p-univ})
pour l'action de $\Sym \ell$ sur
$\{1,\dots,\ell\}$ (et ce bien sûr dans le cas de l'anneau de coefficients
$\mathbf{F}_\ell$).
\emph{A priori}, cette opération $P^{\univ}_{\Sym
\ell}$ n'est définie qu'à \guil{presqu'égalité près}. Le
corollaire~\ref{corollaire-cohomologie-b-sym-ell} montre que dans ce
cas particulier, elle est définie sans ambiguité aucune. 

D'après les
proposition~\ref{proposition-decomposition-motif-coinvariants-bmu},
\ref{proposition-isomorphisme-coinvariants-mu-ell-z-sur-ell} et
\ref{proposition-motif-sym-ell}, on peut définir de manique unique des
transformations naturelles \[P^i\colon \tH^{2r,r}(\mathcal X)\to
\tH^{2r+2i(\ell-1),r+i(\ell-1)}(\mathcal X)\] et
\[B^i\colon \tH^{2r,r}(\mathcal X)\to
\tH^{2r+2i(\ell-1)+1,r+i(\ell-1)}(\mathcal X)\] pour tout $i\in\mathbf{Z}$
telles que pour tout $x\in \tH^{2r,r}(\mathcal X)$, on ait :
\[P_\ell(x)=\sum_{j\geq 0}P^{r-j}(x)d^j+\sum_{j\geq 0}B^{r-1-j}(x)cd^j\]
et que $P^i=0$ pour $i\geq r+1$ et $B^i=0$ pour $i\geq r$.

Il résulte de \cite[Proposition~3.6]{voevodsky-reduced} que $P^i=0$ et
$B^i=0$ si $i<0$. On a ainsi :
\[P_\ell(x)=\sum_{j=0}^rP^{r-j}(x)d^j+
\sum_{j=0}^{r-1}B^{r-1-j}(x)cd^j\;\text{.}\]
Que $\beta P_\ell(x)=0$ (cf. théorème~\ref{theoreme-annulation-bockstein})
revient à dire que $B^i(x)=\beta P^i(x)$ et $\beta B^i(x)=0$.

\subsection{Les opérations stables $P^i$ et $B^i$}

\begin{definition}\label{definition-operation-stable}
Un nombre premier $\ell$ et un corps parfait $k$ ayant été fixés, une
opération cohomologique stable $F$ de bidegré $(i,j)\in\mathbf{Z}^2$ est
la donnée, pour tous $(p,q)\in\mathbf{Z}^2$ de transformations naturelles
(pour $\mathcal X\in\Hopt[k]$)
\[F\colon \tH^{p,q}(\mathcal X)\to \tH^{p+i,q+j}(\mathcal X)\]
commutant aux isomorphismes de suspension $\tH^{\star,\star}(\mathcal
X)\isomto \tH^{\star+2,\star+1}(\mathcal X\wedge
(\mathbf{A}^1/\mathbf{A}^1-\{0\}))$ et $\tH^{\star,\star}(\mathcal
X)\isomto \tH^{\star+1,\star}(\mathcal X\wedge S^1)$ donnés par la
multiplication à droite par les classes tautologiques
dans $\tH^{2,1}(\mathbf{A}^1/(\mathbf{A}^1-\{0\}))$ et $\tH^{1,0}(S^1)$.
On dira que $F$ est de degré $i$ et de poids $j$.
\end{definition}

On considère ici chaque foncteur $\tH^{p,q}$ comme un préfaisceau
d'\emph{ensembles} sur $\Hopt[k]$. Cependant, les transformations
naturelles qui composent une opération cohomologique stable sont
automatiquement additives.

Une opération cohomologique stable importante est le Bockstein $\beta$.
Elle est de bidegré $(1,0)$.

Une opération cohomologique stable est déterminée par son action sur
les classes dans $\tH^{2r,r}(\mathcal X)$ pour pour $r\geq 0$ et $\mathcal
X\in\Hopt$ et une famille de transformations naturelles 
$\tH^{2r,r}(\mathcal X)\to \tH^{2r+i,r+j}(\mathcal X)$ pour $r\geq 0$
s'étend en une opération cohomologique stable si et seulement si elle commute
aux isomorphismes $\tH^{\star,\star}(\mathcal
X)\isomto \tH^{\star+2,\star+1}(\mathcal X\wedge
(\mathbf{A}^1/\mathbf{A}^1-\{0\}))$ donnés par la multiplication par $t\in
\tH^{2,1}(\mathbf{A}^1/(\mathbf{A}^1-\{0\}))$ ($t$ est la classe de Thom du
fibré trivial de rang $1$ sur $\Spec k$).

On définit de façon évidente une composition des opérations cohomologiques
stables. Pour tout $\lambda \in H^{i,j}(k)$, la multiplication \emph{à
gauche} par $\lambda$ définit une opération cohomologique stable de
bidegré $(i,j)$. On dispose ainsi d'un plongement de l'algèbre
$H^{\star,\star}(k)$ dans l'algèbre des opérations cohomologiques stables.
On considèrera principalement la structure de $H^{\star,\star}(k)$-module
(à gauche) sur l'algèbre des opérations cohomologiques stables provenant de
cette inclusion. Ainsi, pour $F$ une opération cohomologique stable,
$\lambda F$ est définie par $(\lambda F)(x)=\lambda \cdot F(x)$. On dispose
également d'une structure de module-$H^{\star,\star}(k)$\;\footnote{C'est
ainsi que l'on notera les structures de modules à droite.} définie par
$(F\lambda)(x)=F(\lambda x)$.

\begin{proposition}
Les opérations $P^i$ et $B^i$ définies précédemment s'étendent en des
opérations cohomologiques stables de bidegrés respectifs
$(2i(\ell-1),i(\ell-1))$ et $(2i(\ell-1)+1,i(\ell-1))$ qui vérifient
$B^i=\beta P^i$.
\end{proposition}

Si on note $-\cup t$ la multiplication à droite par $t$ :
$\tH^{\star,\star}(\mathcal X)\to \tH^{\star+2,\star+1}(\mathcal X\wedge
\mathbf{A}^1/(\mathbf{A}^1-\{0\}))$, comme on sait déjà que $B^i=\beta
P^i$, il suffit de vérifier que pour tout $x\in \tH^{2r,r}(\mathcal X)$, on
a $P^i(x\cup t)=P^i(x)\cup t$. D'après la
proposition~\ref{proposition-compatibilite-multiplication}, on a
$P_\ell(x\cup t)=P_\ell(x)\cup P_\ell(t)\in \tH^{\star,\star}(\mathcal
X\wedge (\mathbf{A}^1/\mathbf{A}^1-\{0\})\wedge (\Bet \Sym\ell)_+)$.
D'après la proposition~\ref{proposition-action-classes-de-thom}, on a
$P_\ell(t)=td$, ce qui montre que $P_\ell(x\cup t)=(P_\ell(x)\cup t)d$.
L'identification des coefficients montre que l'on a bien $P^i(x\cup
t)=P^i(x)\cup t$.

\subsection{Premières propriétés}
\label{subsection-premieres-proprietes}

\begin{proposition}\label{proposition-p-0-b-0}
$P^0$ est l'identité et $B^0=\beta$.
\end{proposition}

D'après \cite[Proposition~3.7]{voevodsky-reduced}, pour tout $r\geq 0$,
l'action de $P^0$ sur $\tH^{2r,r}$ est la multiplication par un nombre
$\lambda_r\in\mathbf{Z}/\ell\mathbf{Z}$. La condition de stabilité se
traduit en disant que $\lambda_r=\lambda_{r+1}$ pour tout $r\geq 0$. Ainsi,
il existe $\lambda\in\mathbf{F}_\ell$ tel que $P^0=\lambda\Id$. Pour
montrer que $\lambda=1$, on peut par exemple vérifier que pour la classe
$1\in H^{0,0}(\Spec k)\simeq \tH^{0,0}(S^0)$, on a $P_\ell(1)=1\in
H^{0,0}(\Bet \Sym \ell)$ ce qui implique que $P^0(1)=1$.

\begin{proposition}\label{proposition-p-d-un-produit}
Soit $x\in\tH^{p,q}(\mathcal X)$. Soit $x'\in \tH^{p',q'}(\mathcal X')$.

Pour tout $n\geq 0$, on a une égalité dans $\tH^{\star,\star}(\mathcal
X\wedge \mathcal X')$ :
\[P^n(xy)=\left\{\begin{array}{ll}\sum_{i+j=n}P^i(x)P^j(y) &
\text{si }\ell\neq 2\\
\sum_{i+j=n}P^i(x)P^j(y)+
\tau \sum_{i+j=n-1}B^i(x)B^j(y) & \text{si }\ell=2\;\footnotemark
\end{array}\right.\]
Par ailleurs, $\beta(xx')=\beta(x)\cdot x'+(-1)^p x\cdot \beta(x')$.
\end{proposition}
\footnotetext{Pour la définition de $\tau\in H^{0,1}(k)$, voir
proposition~\ref{proposition-formule-u-carre}. L'autre élément $\rho\in
H^{1,1}(k)$ qui y est défini satisfait aussi la relation $\rho=\beta(\tau)$.}
La formule pour $\beta$ est indiquée pour mémoire. Pour établir la première
série de formules, en prenant garde aux signes, on peut se ramener au cas
où $(p,q)=(2r,r)$ et $(p',q')=(2r',r')$ et alors cela résulte de l'identité
$P_\ell(xy)=P_\ell(x)P_\ell(y)\in\tH^{\star,\star}(\mathcal X\wedge
\mathcal X'\wedge (\Bet \Sym\ell)_+)$ (cf.
proposition~\ref{proposition-compatibilite-multiplication}).

\begin{proposition}\label{proposition-p-r-sur-h-2r-r}
Pour tout $x\in\tH^{2r,r}(\mathcal X)$, on a $P^r(x)=x^\ell\in
\tH^{2r\ell,r\ell}(\mathcal X)$.
\end{proposition}

Il s'agit de déterminer le \guil{terme constant} du développement de
$P_\ell(x)$. Comme les restrictions de $c$ et $d$ au point-base de
$\Bet\Sym \ell$ sont nulles, il s'agit en fait de calculer l'image de
$P_\ell(x)$ par la flèche de restriction \[\tH^{2r\ell, r\ell}(\mathcal
X\wedge (\Bet \Sym \ell)_+)\to \tH^{2r\ell, r\ell}(\mathcal X)\] induite
par le point-base $\Spec k\to\Bet \Sym\ell$. Celui-ci est induit par la
classe du $\Sym \ell$-torseur trivial $U$ sur $S=\Spec k$. D'après la
proposition~\ref{proposition-compatibilite-bg}, $P^r(x)$ est l'image de $x$
par la construction $P$ pour ce torseur trivial et l'action de $\Sym \ell$
sur $A=\{1,\dots,\ell\}$. En revenant à la
définition~\ref{definition-operation-totale}, il est alors évident que $P$
n'est alors que l'élévation à la puissance $\#A=\ell$.

\begin{corollaire}
Soit $x\in \tH^{p,q}(\mathcal X)$. Soit $n\geq 0$. Si on a $p-q<n$ et
$q\leq n$, alors $P^n(x)=0$.
\end{corollaire}

En effet, quitte à remplacer $\mathcal X$ par $\mathcal X\wedge
S^{n-p+q-1}\wedge \Gm^{\wedge n-q}$, on peut supposer que $p=2n-1$ et
$q=n$. Si on note $e\in H^{1,0}(S^1)$ la classe tautologique, il faut
montrer l'annulation de $P^n(x)e=P^n(xe)=(xe)^\ell\in
H^{2n\ell,n\ell}(\mathcal X\wedge S^1)$ (la dernière égalité est donnée par
la proposition~\ref{proposition-p-r-sur-h-2r-r}). On conclut en remarquant que
$(xe)^\ell$ peut être identifié (au signe près) au produit de $x^\ell\in
H^{2n\ell-\ell,n\ell}(\mathcal X)$ et de $e^\ell=0\in H^{\ell,0}(S^1)$.

\begin{proposition}\label{proposition-p-c-1}
Soit $X\in \Sm[k]$. Soit $L$ un fibré en droites sur $X$. On a alors :
\[P(c_1(L))=c_1(L)^\ell + c_1(L)d\in \tH^{2\ell, \ell}(X\times
\Bet\Sym\ell)\;\text{.}\]
Autrement dit, $P^0(c_1(L))=c_1(L)$, $P^1(c_1(L))=c_1(L)^\ell$ et les
autres opérations $P^i$ et $B^i$ s'annulent sur $c_1(L)$.
\end{proposition}

On sait déjà \emph{a priori} que
$P(c_1(L))=P^1(c_1(L))+P^0(c_1(L))d+B^0(c_1(L))c$. La formule résulte
simplement du fait que $P^0$ est l'identité, $B^0$ le Bockstein
(cf. proposition~\ref{proposition-p-0-b-0})
et que sur
$H^{2,1}$, $P^1$ soit l'élévation à la puissance $\ell$ (cf.
proposition~\ref{proposition-p-r-sur-h-2r-r}).

Il était également possible de partir de la
proposition~\ref{proposition-action-classes-de-thom}, puis de montrer
l'annulation des classes de Chern $c_i(\xi)$ pour $1\leq i\leq \ell-1$ afin
d'obtenir une formule pour $P(t_L)$ de laquelle on peut déduire une formule
pour $P(c_1(L))$ en observant que $c_1(L)$ est l'image inverse de $t_L$ par
le morphisme $X\to \Th_X L$ donné par la section nulle.

\begin{proposition}
Dans $H^{\star,\star}(\Bet\mu_\ell)$, on a les relations suivantes pour
$i\geq 0$ et $j\geq 0$ :
\[P^i(v^j)=\binom j i v^{j+(\ell-1)i}\;\text{,}\quad
\beta P^i(v^j)=0\;\text{,}\]
\[P^i(uv^j)=\binom j i uv^{j+(\ell-1)i}\;\text{,}\quad
\beta P^i(uv^j)= \binom j i v^{j+(\ell-1)i+1}\;\text{.}\]
Dans $H^{\star,\star}(\Bet\Sym \ell)$, on a les relations suivantes pour
$i\geq 0$ et $j\geq 0$ :
\[P^i(d^j)=(-1)^i\binom {(\ell-1)j} i d^{i+j}\;\text{,}\quad
\beta P^i(d^j)=0\;\text{,}\]
\[P^i(cd^j)=(-1)^{i}\binom {(\ell-1)(j+1)-1} i cd^{i+j}\;\text{,}\quad
\beta P^i(cd^j)=(-1)^i\binom {(\ell-1)(j+1)-1} i d^{i+j+1}\;\text{.}\]
\end{proposition}

Dans $H^{\star,\star}(\Bet\mu_\ell\times \Bet \Sym \ell)\simeq
H^{\star,\star}(\Bet\mu_\ell) \otimes_{H^{\star,\star(k)}}
H^{\star,\star}(\Bet \Sym \ell)$, on a $P(v)=v^\ell\otimes 1+v\otimes d$
d'après la proposition~\ref{proposition-p-c-1}. On en déduit le
développement de $P(v^j)=P(v)^j$ en utilisant la formule du binôme. Ceci
permet de calculer les images de $v^j$ par $P^i$ et $B^i$. Les formules
pour $uv^j$ s'en déduisent en utilisant que les seules opérations $P^i$ et
$B^i$ ne s'annulant pas sur $u$ sont $P^0=\Id$ et $B^0=\beta$. Le résultat
pour $d^j$ et $cd^j$ se déduit des autres formules en utilisant le fait que
le plongement canonique de $H^{\star,\star}(\Bet\Sym \ell)$ dans
$H^{\star,\star}(\Bet\mu_\ell)$ soit compatible avec l'action des
opérations $P^i$ et $B^i$ (se ramener au cas où $k$ possède une racine
$\ell$-ième de l'unité).

\subsection{Relations d'Adem}
\label{subsection-adem}

On énonce ci-dessous les relations d'Adem. Pour le cas $\ell=2$, on a noté
$\Sq^{2n}=P^n$ et $\Sq^{2n+1}=B^n$. En général, l'opération cohomologique
stable $\Sq^i$ est de bidegré $(i,\floor{\frac i 2})$.

\begin{theoreme}[Relations d'Adem pour $\ell=2$]
\label{theoreme-adem-modulo-2}
On suppose que $a$ et $b$
sont des entiers naturels tels que $0<a<2b$.
\begin{enumerate} \item[(1)]
Si $a$ est pair et $b$ impair : \[\Sq^a\Sq^b=\sum_{j=0}^{\floor{\frac a 2}}
\binom{b-1-j}{a-2j}\Sq^{a+b-j}\Sq^j+
\sum_{\underset{\text{impair}}{j=1}}^{\floor{\frac a
2}}\binom{b-1-j}{a-2j}\rho\Sq^{a+b-j-1}\Sq^j\]
\item[(1')] Si $a$ et $b$
sont impairs :
\[\Sq^a\Sq^b=\sum_{\underset{\text{impair}}{j=0}}^{\floor{\frac a 2}}
\binom{b-1-j}{a-2j}\Sq^{a+b-j}\Sq^j\]
\item[(2)] Si $a$ et $b$ sont pairs :
\[\Sq^a\Sq^b=\sum_{j=0}^{\floor{\frac a 2}} \tau^{j\;\mathrm{mod}\;2}
\binom{b-1-j}{a-2j}\Sq^{a+b-j}\Sq^j\]
\item[(2')] Si $a$ est impair et $b$
pair : \[\quad\Sq^a\Sq^b=\sum_{\underset{\text{pair}}{j=0}}^{\floor{\frac a
2}} \binom{b-1-j}{a-2j}\Sq^{a+b-j}\Sq^j+
\sum_{\underset{\text{impair}}{j=1}}^{\floor{\frac a 2}}
\binom{b-1-j}{a-1-2j}\rho\Sq^{a+b-j-1}\Sq^j\]
\end{enumerate}
\end{theoreme}

Les identités (1') et (2') se déduisent respectivement des identités (1) et
(2) par application du Bockstein. Les formules apparaissant dans
\cite[Theorem~2.3]{voevodsky-reduced}
souffraient de la présence de coquilles dans les
termes faisant intervenir $\rho$. On a vérifié par ordinateur que les
coefficients des formules ci-dessus étaient les bons dans le cas où
$a+b\leq 100$\;\footnote{Le code source C du programme ayant permis cette
vérification est disponible à l'adresse
\url{http://www.math.u-psud.fr/~riou/doc/adem.tar.gz}.}. Le principe de la
démonstration rédigée dans \cite{voevodsky-reduced} est cependant correct.
(Cette démonstration est légèrement plus compliquée que dans le cas $\ell\neq
2$ qui apparaît plus bas.)

\begin{theoreme}[Relations d'Adem pour $\ell\neq 2$]
Si des entiers $a$ et
$b$ vérifient $0<a<\ell b$, on a : \[P^aP^b=\sum_{t=0}^{\floor{\frac a
\ell}}(-1)^{a+t}\binom{(\ell-1)(b-t)-1}{a-\ell t}P^{a+b-t}P^t\]
S'ils vérifient $0<a\leq \ell b$, on a :
\[P^a\beta P^b=\sum_{t=0}^{\floor{\frac
a \ell}}(-1)^{a+t}\binom{(\ell-1)(b-t)}{a-\ell t}\beta P^{a+b-t}P^t
+\sum_{t=0}^{\floor{\frac {a-1}
\ell}}(-1)^{a+t+1}\binom{(\ell-1)(b-t)-1}{a-\ell t-1}P^{a+b-t}\beta P^t\]
\end{theoreme}

Comme indiqué dans
\cite[Theorem~10.3]{voevodsky-reduced}, le principe de la démonstration est
le même que dans \cite[Chapter~VIII]{steenrod}.
Pour en faciliter la compréhension et
la reconstitution des moindres détails pour le lecteur intéressé, nous
allons préciser certains points. Traitons par exemple le cas de la deuxième
identité. Il s'agit d'établir une égalité entre deux opérations
cohomologiques stables. Pour vérifier une telle égalité,
il suffit de la tester sur des classes
$x\in \tH^{2r,r}(\mathcal X)$ pour $\mathcal X\in\Hopt[k]$ et $r\geq 0$.
Nous ne le ferons que pour une infinité d'entiers $r$ ayant une forme très
particulière, mais ce sera suffisant puisque si l'identité est vraie pour
un certain $r$, elle l'est automatiquement pour tous les entiers $r'\leq
r$.

La démonstration repose sur le théorème de symétrie
(théorème~\ref{theoreme-symetrie}) qui énonce que
$P_\ell(P_\ell(x))\in \tH^{2r\ell^2,r\ell^2}(\mathcal X\wedge
(\Bet\Sym\ell\times \Bet\Sym\ell)_+)$ est invariant par l'échange des deux
facteurs $\Bet\Sym\ell$. En utilisant les formules du
\S\ref{subsection-premieres-proprietes}, il est possible de développer
$P_\ell(P_\ell(x))$ pour l'écrire sous la forme d'une somme de termes
$\lambda\otimes c^id^j\otimes c^{i'}d^{j'}$ avec $\lambda\in
\tH^{\star,\star}(\mathcal X)$ et $i,i',j,j'$ des entiers naturels, $i$ et
$i'$ devant valoir $0$ ou $1$. L'échange des deux facteurs envoie 
$\lambda\otimes c^id^j\otimes c^{i'}d^{j'}$ sur 
$(-1)^{ii'}\lambda\otimes c^{i'}d^{j'}\otimes c^id^j$. Le théorème de
symétrie implique ainsi un certain nombre de relations entre les
coefficients.

Précisément, pour $r>a+b$, on pose $p=r-b$ et $q=r\ell-a$ et on s'intéresse
au coefficient devant $cd^{p-1}\otimes d^q$ dans $P_\ell(P_\ell(x))$. Le
calcul montre que ce coefficient vaut :
\[\sum_{\delta=0}^a (-1)^\delta \binom{(\ell-1)(p-\delta)-1}{\delta}
P^{a-\delta}\beta P^{b+\delta}(x)\;\text{.}\]
Le terme correspondant à $\delta=0$ est précisément $P^a\beta P^b(x)$ que
l'on cherche à réexprimer autrement. Il convient ainsi de s'arranger pour
que ce soit le seul terme non nul.
J'affirme que si $N$ est assez grand, c'est le cas pour
$r=\ell^N+b$\;\footnote{Ainsi, $p=r-b=\ell^N$. C'est le même choix qui est
fait dans \cite[p. 121]{steenrod}, mais des coquilles peuvent faire croire
à un choix différent ($p=\frac{\ell^N-1}{\ell-1}$), ce qui entraînerait,
qu'outre le terme voulu $P^a\beta P^b$ apparaîtrait
le terme $P^{a-1}\beta P^{b+1}$ ainsi que quelques autres.}.
Il s'agit de montrer que le coefficient binomial
$\binom{(\ell-1)(p-\delta)-1}{\delta}$ est divisible par $\ell$ pour $1\leq
\delta \leq a$ et $N$ assez grand.
On a :
\[\binom{(\ell-1)(p-\delta)-1}{\delta}=\binom{(\ell-1)\ell^N-1-(\ell-1)\delta}{\delta}\;\text{.}\]
Pour $x$ et $y$ des entiers naturels, par un raisonnement sur les
valuations $\ell$-adiques des factorielles et des coefficients binomiaux,
on obtient que $\binom x y\not\equiv 0 \mod \ell$ si et seulement
s'il n'y a pas de
retenue lors du calcul de la soustraction $x-y$ en termes des
développements $\ell$-adiques.

Il est évident que le développement de $(\ell-1)\ell^N-1$ se termine par $N$
chiffres égaux à $\ell-1$. Si $N$ est assez grand, aucune retenue
n'apparaît donc quand on soustrait $(\ell-1)\delta$ à $(\ell-1)\ell^N-1$.
S'il n'y avait pas de retenue dans la soustraction de $\delta$ à
$(\ell-1)\ell^N-1-(\ell-1)\delta$, il viendrait donc qu'il n'y aurait pas de
retenue dans l'addition de $(\ell-1)\delta$ et de $\delta$, ce qui est
manifestement faux si $\delta\neq 0$. Par conséquent, pour $1\leq
\delta\leq a$, on a bien $\binom{(\ell-1)(p-\delta)-1}{\delta}\equiv 0\mod
\ell$ pour $N$ assez grand.

Sous l'hypothèse que $r=\ell^N+b$ soit assez grand, on a ainsi obtenu que
le coefficient de $cd^{p-1}\otimes d^q$ (pour $p=r-b=\ell^N$ et
$q=r\ell-a$) dans $P_\ell(P_\ell(x))$ est précisément $P^a\beta P^b(x)$.
Par le théorème de symétrie, il est égal au coefficient de $d^q\otimes
cd^{p-1}$, qui vaut :
\[\sum_{t=0}^{a+b}(-1)^{a+t}\binom{(\ell-1)(r-t)}{(\ell-1)r-a+t}\beta
P^{a+b-t}P^t(x)+
\sum_{t=0}^{a+b}(-1)^{a+t+1}\binom{(\ell-1)(r-t)-1}{(\ell-1)r-a+t}P^{a+b-t}
\beta P^t(x) \;\text{.}\]
En utilisant l'identité $\binom x y=\binom x {x-y}$, on peut le récrire :
\[\sum_{t=0}^{\floor{\frac a \ell}}
(-1)^{a+t}\binom{(\ell-1)(r-t)}{a-\ell t}\beta
P^{a+b-t}P^t(x)+
\sum_{t=0}^{\floor{\frac{a-1}\ell}}
(-1)^{a+t+1}\binom{(\ell-1)(r-t)-1}{a-\ell t-1}P^{a+b-t}
\beta P^t(x) \;\text{.}\]
On a $\binom{(\ell-1)(r-t)}{a-\ell t}=\binom{(\ell-1)(b-t)+(\ell-1)\ell^N}
{a-\ell t}$. Pour $N$ assez grand, ce coefficient binomial est congru
modulo $\ell$ à $\binom{(\ell-1)(b-t)}{a-\ell t}$. En effet, si $y$ est un
entier naturel donné, la réduction modulo $\ell$ de $\binom x y$ pour
$x\geq 0$ ne dépend que de la classe de congruence de $x$ modulo $\ell^N$
pour $N$ assez grand. C'est à cet endroit que l'on utilise $a\leq \ell
b$ puisque cela permet de s'assurer que $(\ell-1)(b-t)\geq 0$ sous
l'hypothèse que $0\leq t\leq \frac a \ell$.
De même, pour chacun des termes de la deuxième somme, on a une congruence
modulo $\ell$ entre $\binom{(\ell-1)(r-t)-1}{a-\ell t-1}$ et 
$\binom{(\ell-1)(b-t)-1}{a-\ell t-1}$,
ce qui achève la démonstration de la deuxième série
des relations d'Adem énoncées pour $\ell\neq 2$. (La première série
d'identités donnant une formule pour $P^aP^b$ pour $0<a<\ell b$
s'obtient de façon similaire en considérant les coefficients
devant $d^p\otimes d^q$ et $d^q\otimes d^p$ avec $p=r-b$, $q=r\ell-a$
comme ci-dessus, mais avec $r=\frac{\ell^N-1}{\ell-1}+b$ pour
$N$ assez grand.)

\section{L'algèbre de Steenrod et sa duale}
\label{section-algebre-steenrod}

On fixe un corps parfait $k$ et un nombre premier $\ell$ inversible dans
$k$. On notera simplement $H^{\star,\star}$ l'algèbre de cohomologie
motivique $H^{\star,\star}(k)$ de $k$ à coefficients dans $\Fl$.

\subsection{L'algèbre $A^{\star,\star}$ et la comultiplication $\Psi^\star$}

\begin{definition}
Un monôme en $\beta$ et les $P^i$ est une expression de la forme
\[\beta^{\varepsilon_0}P^{s_1}\beta^{s_1}\dots P^{s_k}\beta^{\varepsilon_k}\]
avec $k\geq 0$, $\varepsilon_i\in\{0,1\}$ et $s_i>0$. On note $k$ la
longueur du monôme.
Le monôme est dit admissible si pour $1\leq i\leq k-1$, on a $s_i\geq \ell
s_{i+1}+\varepsilon_i$.
\end{definition}

Par composition des opérations $P^i$ et $\beta$ définies dans la section
\ref{section-operations}, à chaque monôme est associée une opération
cohomologique stable dont le bidegré est $(n,w)$ avec
$w=(\ell-1)\sum_{i=1}^k s_i$ et $n=2w+\sum_{i=0}^k\varepsilon_i$. L'entier
$w$ sera appelé le poids du monôme, et $n$ son degré.

\begin{definition}
On note $A^{\star,\star}$ le $H^{\star,\star}$-module libre formellement
défini comme ayant pour base l'ensemble des monômes admissibles.
\end{definition}

\begin{proposition}\label{proposition-plongement-a-star-star}
Le morphisme évident de $H^{\star,\star}$-modules de $A^{\star,\star}$ vers
l'algèbre des opérations cohomologiques stables est injectif et identifie
$A^{\star,\star}$ à une sous-algèbre de l'algèbre des opérations
cohomologiques stables.
\end{proposition}

Nous allons commencer par montrer que l'image de $A^{\star,\star}$ par ce
morphisme est une sous-algèbre de l'algèbre des opérations cohomologiques
stables. Compte tenu des formules de la
proposition~\ref{proposition-p-d-un-produit}, il suffit pour cela d'établir
le lemme suivant :

\begin{lemme}
Pour tout monôme $M$ en $\beta$ et les $P^i$, l'opération cohomologique
stable associée à $M$ est une $H^{\star,\star}$-combinaison linéaire
d'opérations cohomologiques stables associées à des monômes admissibles.
\end{lemme}

Nous allons démontrer ce lemme par récurrence sur le poids $w$ du monôme
(le résultat est évident si $w=0$). Supposons donc le résultat obtenu pour
les monômes de poids $<w$. Au cours de cette démonstration, nous dirons
qu'une opération cohomologique stable est admissible si elle est
$H^{\star,\star}$-combinaison linéaire d'opérations associée à des monômes
admissibles.

Fixons donc $w>0$. On voudrait montrer que l'opération cohomologique
associée à un monôme $M=\beta^\varepsilon P^a M'$ de poids $w$ est
admissible. Commençons pour cela par montrer que l'on peut supposer que le
monôme $M'$ est admissible. Par hypothèse de récurrence sur le poids,
l'opération cohomologique associée à $M'$ est une somme d'opérations de la
forme $\lambda M''$ où $M''$ est un monôme admissible et $\lambda\in
H^{\star,\star}$. Il s'agit donc de montrer que les opérations de la forme
$\beta^\varepsilon P^a \lambda M''$ avec $M''$ monôme admissible sont
admissibles. Si $\lambda\not\in H^{0,0}$, en utilisant les formules de la
proposition~\ref{proposition-p-d-un-produit}, il est possible de réexprimer
l'opération $\beta^\varepsilon P^a \lambda M''$ sous la forme d'une
combinaison linéaire d'opérations associées à des monômes de poids $<w$.
Grâce à l'hypothèse de récurrence sur $w$, on peut donc supposer que
$\lambda\in H^{0,0}$, bref on se ramène au cas où $\lambda=1$. On s'est
donc bien ramené au cas où $M'$ est un monôme admissible.

Ainsi, il s'agit de montrer que si $M=\beta^\varepsilon P^a M'$ est un
monôme de poids $w$ avec $M'$ un monôme admissible, alors l'opération
cohomologique associée à $M$ est admissible. Nous allons le montrer par
récurrence descendante sur $\ell a+\varepsilon$ (cette grandeur est
majorée en raison de l'hypothèse sur le poids). Si $M$ n'est pas déjà
admissible, on peut écrire $M'=\beta^{\varepsilon'} P^b M''$ et
$M=\beta^\varepsilon P^a\beta^{\varepsilon'} P^b M''$. En appliquant les
relations d'Adem (cf. \S\ref{subsection-adem}) au monôme non admissible
$P^a\beta^{\varepsilon'}P^b$, on obtient une expression de l'opération
associée à $M$ comme une somme d'opérations de la forme $\mu
\beta^{\eta'}P^{a'} \beta^{\eta''} M'''$ (avec $\mu\in H^{p,q}$ et $M'''$ un
monôme) et il s'agit de montrer que chacun de ces termes induit une
opération admissible. Si le coefficient $\mu$ est de poids $q>0$, c'est le
cas puisque le monôme $\beta^{\eta'}P^{a'} \beta^{\eta''} M'''$ étant de
poids $<w$, l'hypothèse de récurrence s'applique. Ceci permet de ne se
soucier que des termes apparaissant dans les relations d'Adem dont les
coefficients appartiennent à $\Fl$. On remarque qu'on a alors $\ell
a'+\eta'>\ell a+\varepsilon$, ce qui permet de conclure.

\medskip

Pour finir la démonstration de la
proposition~\ref{proposition-plongement-a-star-star}, il reste à démontrer
le fait intéressant, à savoir que les opérations cohomologiques stables
associées aux monômes admissibles sont linéairement indépendantes sur
$H^{\star,\star}$. Pour $w\geq 0$, notons $A^{\star,\star\leq w}$ le
sous-$H^{\star,\star}$-module libre de $A^{\star,\star}$ de base les
monômes admissibles de poids $\leq w$. Il suffit d'établir le lemme
suivant :

\begin{lemme}\label{lemme-injection-evaluation-u-v}
Soit $w\geq 0$. Alors, pour $N$ et $n$ assez grands (par exemple $N>w+1$
et $n>w$),
le morphisme
$H^{\star,\star}$-linéaire d'évaluation sur la classe
$(u\otimes v)^{\otimes n}$
\[A^{\star,\star\leq w}\to H^{\star,\star}((\mu_\ell
\backslash(\mathbf{A}^N-\{0\}))^{2n})\]
est l'inclusion d'un facteur direct.
\end{lemme}

Le $H^{\star,\star}$-module d'arrivée est muni de la base évidente obtenue
par produit tensoriel des éléments
$1,v,v^2,\dots,v^{N-1},u,uv,uv^2,\dots,uv^{N-1}$ (voir la démonstration de
la proposition~\ref{proposition-decomposition-motif-de-b-mu-ell}). On peut
considérer la matrice $B$ de l'évaluation en $u\otimes v \otimes
\dots\otimes u\otimes v$. Essentiellement pour des raisons de poids, les
matrices obtenues pour différents entiers $N>w+1$ (mais à $n$ fixé)
ne diffèrent que par
l'ajout ou la suppression de lignes de zéros.

Il est en principe possible de calculer cette matrice $B$ en utilisant les
formules du \S\ref{subsection-premieres-proprietes}. Pour $\ell\neq 2$, on
peut remarquer que la matrice $B$, \emph{a priori} à coefficients dans
$H^{\star,\star}$ (ou plus précisément de son anneau opposé) est en fait à
coefficients dans $\Fl=H^{0,0}\subset H^{\star,\star}$. Il est évident
aussi que cette matrice est la même que celle correspondant à une étude
homologue pour les opérations de Steenrod dans le cadre topologique
classique. Une variante des arguments conduisant à
\cite[Corollary~2.6, Chapter~VI]{steenrod} montre que cette matrice est la
matrice d'une application $\Fl$-linéaire injective. Par extension des
scalaires à $H^{\star,\star}$, on obtient que le morphisme considéré dans
le lemme est l'inclusion d'un facteur direct.

Si $\ell=2$, on note $I$ l'idéal bilatère de $H^{\star,\star}$ défini par
$I=\oplus_{p\in\mathbf{Z},q>0} H^{p,q}$. Le quotient
$H^{\star,\star}/I$ s'identifie à $\Fl$. L'idéal $I$ étant stable par les
opérations $P^i$ et $\beta$, il vient que les opérations cohomologiques
stables associées aux monômes induisent des transformations
naturelles
\[H^{\star,\star}(\mathcal X)/IH^{\star,\star}(\mathcal X)
\to H^{\star,\star}(\mathcal X)/IH^{\star,\star}(\mathcal X)\]
pour tout $\mathcal X\in\Hopt$.
Comme $\tau\in I$, on peut remarquer que la formule pour $P^n(xy)$ de la
proposition~\ref{proposition-p-d-un-produit} prend modulo $I$ exactement la
forme qu'elle a dans le cas $\ell\neq 2$. Pourvu que l'on raisonne
modulo $I$, il n'y a donc plus lieu de distinguer les cas $\ell\neq 2$ et
$\ell=2$. Ainsi, les calculs faits dans la démonstration de
\cite[Corollary~2.6, Chapter~VI]{steenrod} peuvent être appliqués au cas
$\ell=2$ afin de montrer que si on note $\overline{B}$ la matrice 
à coefficients dans $\Fl$ obtenue en réduisant modulo $I$ les coefficients
de $B$, alors $\overline{B}$ est la matrice d'une application
$\Fl$-linéaire injective. Que le morphisme considéré dans le
lemme~\ref{lemme-injection-evaluation-u-v} soit
l'inclusion d'un facteur direct résulte alors du lemme suivant :

\begin{lemme}
Soit $H^{\star,\star}$ un anneau bigradué \guil{commutatif} (au sens où
$x\cdot y=(-1)^{pq}y\cdot x$ si $x\in H^{p,\star}$ et $y\in H^{q,\star}$).
On suppose que $H^{\star,w}=0$ si $w<0$, que
$H^{\star,0}=H^{0,0}$ et que $H^{0,0}$
est un corps $k$. On dispose ainsi d'un morphisme
d'anneaux canonique $H^{\star,\star}\to k$.

Soit $\varphi\colon M^{\star,\star}\to N^{\star,\star}$ un morphisme (de
bidegré $(0,0)$) entre
$H^{\star,\star}$-modules bigradués libres de rang fini. On suppose que le
morphisme $\overline{\varphi}\colon
k\otimes_{H^{\star,\star}}M^{\star,\star}\to
k\otimes_{H^{\star,\star}}N^{\star,\star}$ est injectif.
Alors, $\varphi$ est l'inclusion d'un facteur direct.
\end{lemme}

Choisissons des bases $(e_j)_{j\in J}$ et $(f_i)_{i\in I}$ formées
d'éléments homogènes de $M^{\star,\star}$ et $N^{\star,\star}$
respectivement. Notons $B=(b_{i,j})_{i,j}$ la matrice de $\varphi$ dans ces
bases. Autrement dit, on a $\varphi(e_j)=\sum_{i} b_{ij}f_i$. Notons
$\overline{B}=(\overline{b_{i,j}})_{i,j}$ la matrice à coefficients dans
$k$ obtenue par réduction. Par hypothèse, c'est une matrice de rang égal au
cardinal de $J$. Quitte à enlever des lignes à cette matrice, pour
obtenir la conclusion voulue, on peut supposer que $I$ et $J$ ont le même
cardinal. Sous cette hypothèse supplémentaire, le but est de montrer que
$\varphi$ est un isomorphisme.

Comme le déterminant de $\overline{B}$ (pour des choix d'ordres totaux sur
$I$ et $J$) est non nul, il existe une bijection $\sigma\colon J\to I$
telle que $\overline{b}_{\sigma(j),j}\neq 0$. On peut évidemment
supposer que $I=J$ et
que $\sigma$ est l'identité. Alors, comme pour tout $i\in I$,
$\overline{b}_{i,i}\neq 0$, le bidegré de cet élément ne peut être que
$(0,0)$, ce qui montre que $e_i$ et $f_i$ sont de même bidegré. Ceci permet
de supposer que $M^{\star,\star}=N^{\star,\star}$, autrement dit que
$\varphi$ est un endomorphisme de $M^{\star,\star}$. Pour tout $w\in
\mathbf{Z}$, notons $W_{\leq w} M^{\star,\star}$ le
sous-$H^{\star,\star}$-module de $M^{\star,\star}$ engendré par ceux des
éléments $(e_i)_{i\in I}$ dont le bidegré $(p,q)$ vérifie $q\leq w$.

Les hypothèses sur $H^{\star,\star}$ font que $\varphi(W_{\leq
w}M^{\star,\star})\subset W_{\leq w}M^{\star,\star}$. Les éléments $e_i$
induisent une base évidente du gradué $\Gr_W M^{\star,\star}$. La matrice
de $\Gr_W(\varphi)$ dans cette base est à coefficients dans $H^{0,0}=k$ et
n'est autre que $\overline{B}$ qui est inversible. Ainsi, $\varphi$ induit
un isomorphisme sur $\Gr_W M^{\star,\star}$. Il résulte d'une utilisation
itérée du lemme des cinq que $\varphi$ est un isomorphisme.

\bigskip

Avant de procéder aux définitions suivantes, précisons que sauf mention du
contraire, on considérera toujours la structure de $H^{\star,\star}$-module
évidente sur $A^{\star,\star}$. Si $M$ et $N$ sont deux
$H^{\star,\star}$-modules bigradués, on peut former le produit tensoriel
$M\otimes_{H^{\star,\star}} N$ en munissant $M$ de la structure de
$H^{\star,\star}$-bimodule-$H^{\star,\star}$ obtenue en combinant la
structure de module à gauche donnée et la structure de module à droite
définie par $m\lambda=(-1)^{\abs{\lambda}\abs{m}}\lambda m$ où
$\lambda\in H^{\star,\star}$ et $m\in M^{\star,\star}$ où $\abs{\lambda}$
(resp. $\abs{m}$) désignent le premier entier constituant le bidegré de
$\lambda$ (resp. $m$).

\bigskip

Si $x\in \tH^{\star,\star}(\mathcal X)$ et $y\in \tH^{\star,\star}(\mathcal
Y)$ sont deux classes de cohomologie d'espaces pointés
arbitraires, on définit une
application $H^{\star,\star}$-linéaire :
\[\left\{\begin{array}{ccc}A^{\star,\star}\otimes_{H^{\star,\star}}A^{\star,\star}&\longrightarrow&
\tH^{\star,\star}(\mathcal X\wedge \mathcal Y)\\
F&\longmapsto& F(x,y) \end{array}\right.\]
définie par $(C\otimes D)(x,y)=(-1)^{\abs{D}\abs{x}}Cx\cdot Dy$ pour $C$ et
$D$ des éléments de $A^{\star,\star}$.

\begin{proposition}\label{proposition-psi}
On définit une application $H^{\star,\star}$-linéaire 
$\Psi^\star\colon A^{\star,\star}\to
A^{\star,\star}\otimes_{H^{\star,\star}} A^{\star,\star}$ telle que pour
tout $F\in A^{\star,\star}$, $\Psi^{\star}(F)$ soit l'unique élément tel
que 
\[(\Psi^{\star}(F))(x,y)=F(xy)\in \tH^{\star,\star}(\mathcal
X\wedge\mathcal Y)\]
pour tous $x\in H^{\star,\star}(\mathcal X)$ et $y\in
H^{\star,\star}(\mathcal Y)$, $\mathcal X\in \Hopt$, $\mathcal
Y\in\Hopt$\;\footnote{Si on souhaite se dispenser d'écrire certains signes,
on peut se contenter du cas où $\abs{x}$ et $\abs{y}$ sont pairs, et même
que les bidegrés respectifs soient de la forme $(2r,r)$ et $(2s,2s)$ pour
$r,s\geq 0$.}.

\begin{itemize}
\item $\Psi^\star(\Id)=\Id\otimes \Id$ ;
\item si $\Psi^\star(F)=\sum_i C_i\otimes D_i$ et $\Psi^\star(G)=\sum_j
C'_j\otimes D'_j$, on a
$\Psi^\star(GF)=\sum_{i,j}(-1)^{\abs{C_i}\abs{D'_j}} C'_jC_i\otimes D'_j
D_i$.
\item cette comultiplication $\Psi^\star$ est cocommutative et
coässociative (au sens bigradué) et la coünité est le morphisme
$\xi_0\colon A^{\star,\star}\to H^{\star,\star}$ d'évaluation sur la classe
$1\in H^{\star,\star}$.
\end{itemize}

\end{proposition}

L'unicité de $\Psi^{\star}(F)$ pour tout $F\in A^{\star,\star}$ résulte du
lemme suivant qui est une conséquence immédiate du
lemme~\ref{lemme-injection-evaluation-u-v} :

\begin{lemme}\label{lemme-injection-evaluation-u-v-double}
Soit $w\geq 0$. Alors, pour $N$ et $n$ assez grands,
le morphisme
$H^{\star,\star}$-linéaire d'évaluation sur $((u\otimes v)^{\otimes
n},(u\otimes v)^{\otimes n}))$
\[A^{\star,\star\leq w}\otimes_{H^{\star,\star}} A^{\star,\star\leq w}
\to H^{\star,\star}((\mu_\ell \backslash(\mathbf{A}^N-\{0\}))^{4n})\]
est l'inclusion d'un facteur direct.
\end{lemme}

Si $\Psi^\star(F)$ et $\Psi^\star(G)$ existent, alors on vérifie
immédiatement que $\Psi^{\star}(GF)$ existe et est donné par la formule
donnée ci-dessus. Pour obtenir l'existence de $\Psi^\star(F)$, on se ramène
donc au cas où $F=P^i$ ou $F=\beta$, ce pour quoi la
proposition~\ref{proposition-p-d-un-produit} donne des formules explicites. Par
exemple, $\Psi^\star(\beta)=\beta\otimes \Id+\Id\otimes\beta$ et si
$\ell\neq 2$, $\Psi^\star(P^n)=\sum_{i+j=n}P^i\otimes P^j$.

Enfin, à propos des propriétés de la comultiplication $\Psi^\star$,
l'énoncé concernant la coünité vient de ce que $1\in H^{\star,\star}$ soit
l'unité, la cocommutativité vient de la commutativité au sens gradué de la
multiplication sur la cohomologie motivique et la coässociavitié résulte
formellement de l'associativité de cette multiplication. Plus précisément,
grâce à une version en trois variables du
lemme~\ref{lemme-injection-evaluation-u-v-double} qui en utilise deux, on
peut montrer que pour tout $F\in A^{\star,\star}$,
$(\Id\otimes \Psi^{\star})(\Psi^{\star}(F))$ et
$(\Psi^{\star}\otimes \Id)(\Psi^{\star}(F))$ sont deux éléments du produit
tensoriel triple
$A^{\star,\star}\otimes_{H^{\star,\star}}A^{\star,\star}\otimes_{H^{\star,\star}}A^{\star,\star}$
qui ont la même vertu caractéristique que si on les écrit comme une somme
$\sum_i A_i\otimes B_i\otimes C_i$, alors quelles que soient les classes de
cohomologie motivique $x$, $y$ et $z$ (que l'on suppose de degrés pairs pour
ne pas avoir à introduire de signe), on ait :
\[F(xyz)=\sum_i A_i(x)B_i(y)C_i(z)\;\text{.}\]

\subsection{L'algèbre duale $A_{\star,\star}$}

\begin{definition}
Pour tout $(i,j)\in\mathbf{Z}^2$, on note $A_{i,j}$ l'ensemble des familles
d'homomorphismes de groupes $\varphi\colon A^{p,q}\to H^{p-i,q-j}$ pour
tout $(p,q)\in\mathbf{Z}^2$, telles que pour tout $\lambda\in
H^{\star,\star}$ et $F\in A^{\star,\star}$, on ait $\varphi(\lambda
F)=\lambda \varphi(F)$.
\end{definition}

L'accouplement canonique $A^{p,q}\times A_{i,j}\to H^{p-i,q-j}$ qui à
$(F,\varphi)$ associe $\varphi(F)$ sera noté $\acc{F}{\varphi}$. Il vérifie
$\acc{\lambda F}{\varphi}=\lambda\acc{F}{\varphi}$ pour $\lambda\in
H^{\star,\star}$, $F\in A^{\star,\star}$ et $\varphi\in A_{\star,\star}$.
On munit $A_{\star,\star}$ d'une structure de module-$H^{\star,\star}$ de
façon à ce que $\acc{F}{\varphi\mu}=\acc{F}{\varphi}\mu$ pour
$F\in A^{\star,\star}$, $\varphi\in A_{\star,\star}$ et $\mu\in
H^{\star,\star}$.

\begin{definition}
On généralise l'accouplement $\acc{\cdot}{\cdot}\colon
A^{\star,\star}\times A_{\star,\star}\to H^{\star,\star}$ en un
accouplement $\acc{\cdot}{\cdot}\colon A^{\star,\star}\times
(A_{\star,\star}\otimes_{H^{\star,\star}} \tH^{\star,\star}(\mathcal X))\to
\tH^{\star,\star}(\mathcal X)$ pour tout $\mathcal X\in\Hopt$ en posant :
\[\acc{F}{\alpha\otimes x}=\acc{F}{\alpha}x\]
pour $F\in A^{\star,\star}$, $\alpha\in A_{\star,\star}$ et $x\in
\tH^{\star,\star}(\mathcal X)$.
\end{definition}

\begin{lemme}\label{lemme-base-duale}
Soit $\mathcal X\in\Hopt$ un espace tel qu'il existe $d\geq 0$ tel que
$\tH^{p,q}(\mathcal X)=0$ pour $p>q+d$ (on dira que $\mathcal X$ est un
espace
cohomologiquement de dimension finie $\leq d$\;\footnote{Les espaces $X_+$
pour $X\in\Sm[k]$ en sont.}). Si $\varphi\colon A^{\star,\star}\to
\tH^{\star,\star}(\mathcal X)$ est une application
$H^{\star,\star}$-linéaire d'un certain bidegré $(-i,-j)$,
il existe un unique élément $\alpha\in
A_{\star,\star}\otimes_{H^{\star,\star}} \tH^{\star,\star}(X)$ tel que pour
tout $F\in A^{\star,\star}$, on ait :
\[\varphi(F)=\acc{F}{\alpha}\;\text{.}\]
Plus précisément, $\alpha=\sum_I \theta(I)^\star \otimes \varphi(P^I)$
où les éléments $\theta(I)^\star\in A_{\star,\star}$ sont définis par
$\acc{P^J}{\theta(I)^\star}=\delta_{IJ}$ (pour tous les monômes admissibles
$P^I$ et $P^J$) et tous sauf un nombre fini des coefficients
$\varphi(P^I)$ sont nuls.

Dans le cas particulier $\mathcal X=S^0$, on obtient que
$A_{\star,\star}$ est un module-$H^{\star,\star}$ libre de base les
éléments $\theta(I)^\star$.
\end{lemme}

Le point-clef à justifier est qu'il n'y ait qu'un nombre fini de monômes
admissibles $P^I$ tels que $\varphi(P^I)\neq 0$. Si un monôme $P^I$ est de
poids $w$, il est de bidegré $(2w+e,w)$ avec $e\geq 0$ et donc
$\varphi(P^I)\in\tH^{2w+e-i,w-j}(\mathcal X)$. Si $\varphi(P^I)\neq 0$, on
a donc $2w+e-i\leq w-j+d$, d'où $w\leq w+e\leq d+i-j$, ce qui donne une
borne sur le poids des monômes $P^I$ tels que $\varphi(P^I)\neq 0$. Il n'y
en a donc qu'un nombre fini.
Ceci étant acquis, le cas particulier où $\mathcal X=S^0$ (et $d=0$) montre
que $A_{\star,\star}$ est un module-$H^{\star,\star}$ libre dont les
éléments $\theta(I)^\star$ forment une base. Dès lors, dans le cas général,
l'élément $\alpha$ cherché doit être une somme finie de termes
$\theta(I)^\star\otimes \beta_I$ avec $\beta_I\in
\tH^{\star,\star}(\mathcal X)$ et cela ne peut être que $\alpha=\sum_I
\theta(I)^\star\otimes \varphi(P^I)$ qui a bien un sens d'après le
point-clef établi ci-dessus.

\medskip

Ce lemme permet de procéder à la définition suivante :

\begin{definition}
Soit $\mathcal X\in \Hopt$ un espace cohomologiquement de dimension finie.
Soit $x\in \tH^{p,q}(\mathcal X)$. On note $\lambda^\star(x)$ l'unique
élément de
$A_{\star,\star}\otimes_{H^{\star,\star}}{\tH^{\star,\star}}(\mathcal X)$
tel que pour tout $F\in A^{\star,\star}$, on ait :
\[F(x)=\acc{F}{\lambda^\star(x)}\;\text{.}\]
\end{definition}

\begin{definition}\label{definition-xi-k-tau-k}
Pour tout $k\geq 0$, on note $M_k=P^{\ell^{k-1}}\dots P^\ell P^1\in
A^{2w,w}$ (avec $w=\ell^k-1$). Ce sont des monômes admissibles. On note
$\xi_k\in A_{2w,w}$ les éléments qui leurs correspondent dans la base duale
formée des éléments $\theta(I)^\star$. Les monômes $M_k\beta\in
A^{2w+1,2w}$ sont également admissibles, on note aussi
$\tau_k\in A_{2w+1,w}$
les éléments correspondants dans la base duale de celle des monômes
admissibles.
\end{definition}

\begin{proposition}\label{proposition-action-monomes-sur-u-v}
Soit $v=c_1(\OO(1))\in H^{2,1}(\mathbf{P}^\infty)\subset H^{2,1}(\Bet
\mu_\ell)$ (cf. proposition~\ref{proposition-decomposition-motif-de-b-mu-ell}).
Si $M$ est un monôme (admissible ou non) tel que $M(v)\neq 0$,
alors il existe $k\geq 0$ tel que $M=M_k$ et alors
$M(v)=M_k(v)=v^{\ell^k}$.

Si $M$ est un monôme (admissible ou non) tel que $M(u)\neq 0$ dans
$H^{\star,\star}(\Bet\mu_\ell)$, alors soit $M=\Id$, soit $M$ est de la
forme $M=M_k\beta$ pour $k\geq 0$ et $M(u)=M_k(v)=v^{\ell^k}$.
\end{proposition}

Grâce aux formules du \S\ref{subsection-premieres-proprietes}, ceci
s'établit par récurrence sur la longueur d'un monôme $M$ tel que $M(v)\neq
0$.

\begin{corollaire}\label{corollaire-lambda-u-et-v}
Pour tout $n\geq 1$, on a :
\begin{eqnarray*}
\lambda^\star(v)&=&\sum_{k=0}^\infty \xi_k\otimes v^{\ell^k}\in
A_{\star,\star}\otimes_{H^{\star,\star}}
H^{\star,\star}(\mathbf{P}^{n-1})\;\text{,}\\
\lambda^\star(u)
&=&\xi_0\otimes u + \sum_{k=0}^\infty \tau_k\otimes v^{\ell^k}\in
A_{\star,\star}\otimes_{H^{\star,\star}}
H^{\star,\star}(\mu_\ell\backslash (\mathbf{A}^n-\{0\}))\;\text{.}
\end{eqnarray*}
\end{corollaire}

Que ceci ait un sens vient du fait que $v$ soit nilpotent dans
$H^{\star,\star}(\mathbf{P}^{n-1})$ pour tout $n\geq 1$.

\begin{definition}\label{definition-accouplement-double-trivial}
On définit un accouplement
\[(A^{\star,\star}\otimes_{H^{\star,\star}} A^{\star,\star})\times
(A_{\star,\star}\otimes_{H^{\star,\star}} A_{\star,\star})\to
H^{\star,\star}\]
par la formule
\[\acc{C\otimes D}{\alpha\otimes \beta}=(-1)^{\abs{D}\abs{\alpha}}
\acc{C}{\alpha}\cdot \acc{D}{\beta}\;\text{.}\]
\end{definition}

\begin{definition}
Si $\alpha$ et $\beta$ sont deux éléments de $A_{\star,\star}$, on note
$\alpha\beta\in A_{\star,\star}$ l'élément défini par l'identité
\[\acc{F}{\alpha\beta}=\acc{\Psi^{\star}F}{\alpha\otimes \beta}\]
pour tout $F\in A^{\star,\star}$. (Ceci définit bien un élément de
$A_{\star,\star}$ puisque $\acc{\Psi^{\star}(\lambda
F)}{\alpha\otimes\beta}=
\acc{\lambda
\Psi^{\star}(F)}{\alpha\otimes\beta}=\lambda\acc{\Psi^{\star}(F)}{\alpha\otimes\beta}$
pour tous $F\in A^{\star,\star}$ et $\lambda\in H^{\star,\star}$.)

Ceci munit $A_{\star,\star}$ d'une structure d'anneau commutatif au sens
bigradué. L'unité est l'élément $\xi_0$ (cf. définitions~\ref{proposition-psi}
et \ref{definition-xi-k-tau-k}).
\end{definition}

Ceci ne fait que refléter les propriétés de la comultiplication
$\Psi^\star$ sur $A^{\star,\star}$ (cf. proposition~\ref{proposition-psi}).
(Afin d'obtenir l'associativité, pour
bien faire, il faudrait en toute rigueur introduire un accouplement entre
les produits tensoriels triples de $A^{\star,\star}$ et $A_{\star,\star}$,
ce pour quoi la seule difficulté réside dans les choix de signes.)

\medskip

Cette multiplication sur $A_{\star,\star}$ vérifie une compatibilité avec
la structure de module-$H^{\star,\star}$ sur $A_{\star,\star}$ :
\[\alpha(\beta\lambda)=(\alpha\beta)\lambda\]
pour $\alpha,\beta\in A_{\star,\star}$ et $\lambda\in H^{\star,\star}$.

En particulier, en faisant $\beta=\xi_0$, on obtient que la multiplication
à droite par $\lambda\in H^{\star,\star}$ pour la structure de
module-$H^{\star,\star}$ sur $A_{\star,\star}$ est donnée par la
multiplication à droite par $\xi_0\lambda$ pour la structure d'anneau sur
$A_{\star,\star}$. Ceci permet aussi d'observer que l'application
$H^{\star,\star}\to A^{\star,\star}$ qui à $\lambda$ associe $\xi_0\lambda$
est un morphisme d'anneaux et on peut réinterpréter ce qui précède en
disant que la structure de module-$H^{\star,\star}$ sur $A_{\star,\star}$
est induite par ce morphisme d'anneaux $H^{\star,\star}\to
A_{\star,\star}$.

Ces considérations permettent de donner un sens à l'anneau (commutatif au
sens bigradué) $A_{\star,\star}\otimes_{H^{\star,\star}}
H^{\star,\star}(\mathcal X)$ intervenant dans la proposition suivante :

\begin{proposition}\label{proposition-lambda-morphisme-d-anneaux}
Soit $\mathcal X\in \Ho$ un espace cohomologiquement de dimension finie.
Alors, l'application $\lambda^\star\colon H^{\star,\star}(\mathcal X)\to
A_{\star,\star}\otimes_{H^{\star,\star}} H^{\star,\star}(X)$ est un
morphisme d'anneaux.
\end{proposition}

Ceci ne présente aucune difficulté.
Dans le cas particulier $\mathcal X=\bullet$, on obtient :

\begin{corollaire}
On dispose d'un morphisme d'anneaux $\lambda^\star\colon H^{\star,\star}\to
A_{\star,\star}$ défini par le fait que pour tous $x\in H^{\star,\star}$ et
$F\in A^{\star,\star}$, on ait l'identité suivante dans $H^{\star,\star}$ :
\[\acc{F}{\lambda^\star(x)}=F(x)\;\text{.}\]
\end{corollaire}

\begin{remarque}\label{remarque-deux-structures-de-module-sur-la-duale}
On n'a pour le moment considéré $A^{\star,\star}$ que comme un
$H^{\star,\star}$-module pour la composition à gauche par la
multiplication (à gauche) par des éléments de $H^{\star,\star}$. On peut
aussi obtenir 
une structure de module-$H^{\star,\star}$ en utilisant la composition à
droite : $F\lambda = F\circ \lambda$ pour tous $F\in A^{\star,\star}$ et
$\lambda\in H^{\star,\star}$.

Ceci fait de $A^{\star,\star}$ un
$H^{\star,\star}$-bimodule-$H^{\star,\star}$. Comme $A_{\star,\star}$ est,
au sens gradué, le groupe des morphismes $A^{\star,\star}\to
H^{\star,\star}$ de $H^{\star,\star}$-modules, la structure de
module-$H^{\star,\star}$ sur $A^{\star,\star}$ induit une structure de
$H^{\star,\star}$-module sur $A_{\star,\star}$. (Précisément, si $\alpha\in
A_{\star,\star}$, $\lambda\in H^{\star,\star}$ et $F\in A^{\star,\star}$,
on a : $\acc{F}{\lambda \alpha}=\acc{F\lambda}{\alpha}$.)
Cette structure de $H^{\star,\star}$-module sur $A_{\star,\star}$ est celle
induite par le morphisme d'anneaux $\lambda^\star\colon H^{\star,\star}\to
A^{\star,\star}$. Plus précisément, on a :
\[\acc{F}{x\alpha}=\acc{F\circ x}{\alpha}=\acc{F}{\lambda^\star(x)\alpha}\]
pour $F\in A^{\star,\star}$, $x\in H^{\star,\star}$ et $\alpha\in
A_{\star,\star}$.
\end{remarque}

La démonstration de cette compatiblité n'étant pas inintéressante, nous
allons la détailler. Choisissons une décomposition $\Psi^\star F=\sum_i
A_i\otimes B_i$ avec $A_i$ et $B_i$ des éléments homogènes de
$A^{\star,\star}$. Observons que pour tous $\mathcal X\in \Hopt$ et
$y\in \tH^{\star,\star}(\mathcal X)$, on a :
\[(F\circ x)(y)=F(xy)=\sum_i (-1)^{\abs{x}\abs{B_i}} A_i(x)B_i(y)\;\text{.}\]
Ainsi, $F\circ x=\sum_i (-1)^{\abs{x}\abs{B_i}} A_i(x)B_i\in A^{\star,\star}$.
On en déduit :
\begin{eqnarray*}
\acc{F}{x\alpha}&=& \acc{F\circ x}{\alpha}\\
&=& \sum_i (-1)^{\abs{x}\abs{B_i}} \acc{A_i(x)B_i}{\alpha} \\
&=& \sum_i (-1)^{\abs{x}\abs{B_i}} A_i(x)\acc{B_i}{\alpha} \\
&=& \sum_i (-1)^{\abs{x}\abs{B_i}} \acc{A_i}{\lambda^\star(x)}\acc{B_i}{\alpha} \\
&=& \sum_i \acc{A_i\otimes B_i}{\lambda^\star(x)\otimes \alpha} \\
&=& \acc{\Psi^\star F}{\lambda^\star(x)\otimes \alpha} \\
&=& \acc{F}{\lambda^\star(x)\alpha}\;\text{.}
\end{eqnarray*}
Ceci étant vrai pour tout $F\in A^{\star,\star}$, on a bien
$x\alpha=\lambda^\star(x)\alpha$.

On a vu que les deux morphismes d'anneaux $H^{\star,\star}\to
A^{\star,\star}$ définis par $x\longmapsto \lambda^\star(x)$ et
$x\longmapsto \xi_0 x$ induisaient respectivement les structures de modules à
gauche et à droite sur $A_{\star,\star}$. Ils permettent de considérer
de deux manières différentes $A_{\star,\star}$ comme une algèbre (comutative
au sens bigradué) sur $H^{\star,\star}$. La structure provenant de
l'application $x\longmapsto \xi_0 x$ est étudiée dans le théorème qui suit
et elle sera appellée \guil{algèbre-$H^{\star,\star}$} :

\begin{theoreme}\label{theoreme-presentation-duale}
En tant qu'algèbre-$H^{\star,\star}$ commutative au sens
bigradué,
$A_{\star,\star}$ admet une présentation par générateurs et relations, les
générateurs étant les éléments $\xi_k\in A_{2(\ell^k-1),\ell^k-1}$ pour
$k\geq 0$ et $\tau_k\in A_{2(\ell^k-1)+1,\ell^k-1}$ pour $k\geq 0$, soumis
à la relation $\xi_0=1$ et à la famille de relations suivantes,
pour tous $k\geq 0$ :
\[\left\{\begin{array}{ll} \tau_k^2 = 0 & \text{si }\ell\neq 2\\
\tau_k^2=\xi_{k+1}\tau+\tau_0\xi_{k+1}\rho+\tau_{k+1}\rho& \text{si }\ell=2
\end{array}\right.\]
\end{theoreme}

Il convient de commencer par vérifier que ces relations sont satisfaites
dans $A_{\star,\star}$.
Dans le cas $\ell\neq 2$, cela résulte simplement de la commutativité au
sens gradué de $A_{\star,\star}$.
La formule pour $\tau_k^2$ dans le cas $\ell=2$
provient de l'examen du coefficient de $v^{2^{k+1}}$
dans le développement de $\lambda^\star(u)^2=\lambda^\star(u^2)$, ce qui
s'obtient en utilisant que $u^2=\tau v+\rho u$
(cf.~proposition~\ref{proposition-formule-u-carre})
et donc que
$\lambda^\star(u^2)=\lambda^\star(\tau)\lambda^\star(v)+
\lambda^\star(\rho)\lambda^\star(u)$, ainsi que les formules évidentes
$\lambda^\star(\tau)=\xi_0 \tau+\tau_0\rho$ et
$\lambda^\star(\rho)=\xi_0\rho$. La démonstration du théorème sera
terminée un peu plus loin, après le corollaire~\ref{corollaire-base-duale}.

\begin{definition}\label{definition-omega}
Pour tout monôme admissible
$\beta^{\varepsilon_0}P^{s_1}\beta^{\varepsilon_1}\dots
P^{s_k}\beta^{\varepsilon_k}$, on peut prolonger la suite des entiers
$(\varepsilon_i)_{i\geq 0}$ et $(s_i)_{i\geq 1}$ par $\varepsilon_i=s_i=0$
pour $i>k$ et on note $r_i=s_i-(\ell s_{i+1}+\varepsilon_i)\geq 0$ pour
tout $i\geq 1$. On obtient ainsi toutes les suites
$I=(\varepsilon_0,r_1,\varepsilon_1,r_2,\varepsilon_2,\dots,r_k,\varepsilon_k,0,\dots)$ d'entiers
nuls à partir d'un certain rang telles que la sous-suite
$(\varepsilon_i)_{i\geq 0}$ soit faite de $0$ et de $1$. À une telle
suite sont associés les éléments de même bidegré
$\omega(I)=\tau_0^{\varepsilon_0}\dots
\tau_k^{\varepsilon_k}\xi_1^{r_1}\dots \xi_k^{r_k}\in A_{\star,\star}$ et
$\theta(I)=\beta^{\varepsilon_0}P^{s_1}\beta^{\varepsilon_1}\dots
P^{s_k}\beta^{\varepsilon_k}\in A^{\star,\star}$
avec $s_i=\sum_{j\geq
i}(r_j+\varepsilon_j)\ell^{j-i}$. On ordonne ces suites
$I=(\varepsilon_0,r_1,\varepsilon_1,\dots)$ en utilisant l'ordre
lexicographique correspondant à un ordre de lecture des lettres de la
droite vers la gauche.
\end{definition}

\begin{proposition}\label{proposition-matrice-triangulaire}
\begin{itemize}
\item Si $I=J$, $\acc{\theta(I)}{\omega(J)}=\pm 1$.
\item Si $I<J$, $\acc{\theta(I)}{\omega(J)}=0$.
\end{itemize}
\end{proposition}

Si $\ell\neq 2$, on peut observer que quels que soient $I$ et $J$, le
coefficient $\acc{\theta(I)}{\omega(J)}\in H^{\star,\star}$ appartient à
$\Fl=H^{0,0}$ (en particulier, si les bidegrés de $\theta(I)$ et
$\omega(J)$ ne sont pas les mêmes, ce coefficient est nul). En effet, il
résulte des définitions, du corollaire~\ref{corollaire-lambda-u-et-v} et de
la proposition~\ref{proposition-lambda-morphisme-d-anneaux} que si on écrit
$\omega(J)=\tau_{i_1}\dots\tau_{i_m}\xi_{j_1}\dots \xi_{j_n}$, alors
$\acc{\theta(I)}{\omega(J)}$ est le coefficient devant
$v^{\ell^{i_1}}\otimes \dots \otimes v^{\ell^{i_m}}\otimes
v^{\ell^{j_1}}\otimes\dots \otimes v^{\ell^{j_m}}$ dans la décomposition
de la classe $\theta(I)(u^{\otimes m}\otimes
v^{\otimes n})\in H^{\star,\star}((\Bet\mu_\ell)^{m+n})$ sur la base
déduite de la décomposition de la
proposition~\ref{proposition-decomposition-motif-de-b-mu-ell}.

Or, à partir des formules du \S\ref{subsection-premieres-proprietes}, il
est évident que la sous-$\Fl$-algèbre de
$H^{\star,\star}((\Bet\mu_\ell)^{m+n})$ engendrée par les éléments
$1^{\otimes i-1}\otimes u\otimes 1^{\otimes m+n-i}$ et 
$1^{\otimes i-1}\otimes v\otimes 1^{\otimes m+n-i}$ pour $1\leq i\leq m+n$
est stable par les opérations $P^k$ et $\beta$ et donc aussi par les
opérations $\theta(I)$. Les coefficients $\acc{\theta(I)}{\omega(J)}$
appartiennent donc bien à $\Fl$.

Ce principe de calcul vaut aussi évidemment pour la version topologique
usuelle de l'algèbre de Steenrod et de sa duale. Les coefficients
$\acc{\theta(I)}{\omega(J)}$ sont donc les mêmes qu'en topologie. La
propriété voulue pour $\ell\neq 2$ se déduit donc de l'énoncé topologique
\cite[Lemma~8]{milnor}.

Pour $\ell=2$, on ne peut pas utiliser le même raisonnement, mais il est
possible d'adapter les arguments de la démonstration de \cite{milnor}. Notons
$J=(\varepsilon_0,r_1,\varepsilon_1,\dots,r_k,\varepsilon_k,0,\dots)$. On
raisonne par récurrence sur $\sum_{i\geq 1}r_i+\sum_{i\geq
0}\varepsilon_i$. Si cette quantité est nulle, le résultat est évident.
Sinon, on note $k$ le plus grand entier tel que $r_k\neq 0$ ou
$\varepsilon_k\neq 0$.

Dans le premier cas, on suppose que $\varepsilon_k=0$ et donc $r_k\geq 1$.
Notons
$J'=(\varepsilon_0,r_1,\varepsilon_1,\dots,r_{k-1},\varepsilon_{k-1},r_k-1,0,\dots)$
de sorte que $\omega(J)=\omega(J')\xi_k$. Pour tout $I$, on a :
\[\acc{\theta(I)}{\omega(J)}=\acc{\Psi^{\star}\theta(I)}{\omega(J')\otimes
\xi_k}\;\text{.}\]
Si $I\leq J$, $I$ est de la forme $I=
(\tilde{\varepsilon}_0,\tilde{r}_1,\tilde{\varepsilon}_1,\dots,
\tilde{r}_{k-1},\tilde{\varepsilon}_{k-1},\tilde{r}_k,0,\dots)$.
À partir des formules $\Psi^\star \beta=\beta\otimes\Id+\Id\otimes \beta$
et $\Psi^\star P^n = \sum_{a+b=n}P^a\otimes P^b+\sum_{a+b=n-1}\tau \beta
P^a\otimes \beta P^b$ et de la formule de la
proposition~\ref{proposition-psi} énonçant une compatibilité de
$\Psi^\star$ à la composition, on peut développer $\Psi^\star
\theta(I)$ en une somme de termes de la forme $F\otimes M$ où $F\in
A^{\star,\star}$ et où $M$ est un monôme (ou $0$ si deux $\beta$ se suivent
dans la composition). Ces termes ne peuvent contribuer à l'accouplement
avec $\omega(J')\otimes \xi_k$ que si $\acc{M}{\xi_k}\neq 0$, c'est-à-dire
si le monôme est $M=M_k$ (cf.
proposition~\ref{proposition-action-monomes-sur-u-v} et
corollaire~\ref{corollaire-lambda-u-et-v}). On observe que l'on ne peut
obtenir un tel terme que si $\tilde{s}_k=\tilde{r}_k\geq 1$. Ainsi, si
$\tilde{r}_k=0$, on a $I<J$ et $\acc{\theta(I)}{\omega(J)}=0$. Sinon,
$\tilde{r}_k\geq 1$ et on peut introduire $I'=(\tilde{\varepsilon}_0,\tilde{r}_1,\tilde{\varepsilon}_1,\dots,
\tilde{r}_{k-1},\tilde{\varepsilon}_{k-1},\tilde{r}_k-1,0,\dots)$ et alors,
dans le développement de $\Psi^\star\theta(I)$, le terme de la forme
$F\otimes M_k$ cherché est $\theta(I')\otimes M_k$, de sorte que
$\acc{\theta(I)}{\omega(J)}=\acc{\theta(I')}{\omega(J')}$. Comme on a
évidemment $I'\leq J'$ et $I'=J'$ si et seulement si $I=J$, on peut
conclure en appliquant l'hypothèse de récurrence à $J'$.

On procède de même si
$J=(\varepsilon_0,r_1,\varepsilon_1,\dots,r_k,\varepsilon_k,0,\dots)$
avec $\varepsilon_k=1$. On note
$J'=(\varepsilon_0,r_1,\varepsilon_1,\dots,r_k,0,0,\dots)$ de sorte que
$\omega(J)=\omega(J')\tau_k$, d'où :
\[\acc{\theta(I)}{\omega(J)}=\acc{\Psi^\star\theta(I)}{\omega(J')\otimes
\tau_k}\;\text{.}\]
Dans ce cas, on cherche quelque terme de la forme $F\otimes M_k\beta$
dans le développement de $\Psi^\star\theta(I)$. Si
$I=(\tilde{\varepsilon}_0,\tilde{r}_1,\tilde{\varepsilon}_1,\dots,
\tilde{r}_{k-1},\tilde{\varepsilon}_{k-1},\tilde{r}_k,
\tilde{\varepsilon}_k,0,\dots)\leq J$, pour
qu'apparaisse un tel terme, il faut que $\tilde{\varepsilon}_k=1$. On peut
alors introduire $I'=(\tilde{\varepsilon}_0,\tilde{r}_1,\tilde{\varepsilon}_1,\dots,
\tilde{r}_{k-1},\tilde{\varepsilon}_{k-1},\tilde{r}_k,0,0,\dots)$ et
obtenir
\[\acc{\omega(I)}{\omega(J)}=\acc{\theta(I')\otimes
M_k\beta}{\omega(J')\otimes \tau_k}=
\acc{\omega(I')}{\omega(J')}\;\text{,}\]
ce qui permet de conclure par récurrence.

\begin{corollaire}\label{corollaire-base-duale}
Les éléments $\omega(I)$ constituent une base de $A_{\star,\star}$ comme
module-$H^{\star,\star}$.
\end{corollaire}

Le fait que la famille soit libre ne pose pas de difficulté : si $\sum_J
\omega(J)\lambda_J=0$ était une relation non triviale, en notant
$I=\min\{J, \lambda_J\neq 0\}$ (le minimum étant pris pour l'ordre total
défini plus haut), on aurait $0\neq \lambda_I=\pm\acc{\theta(I)}{\sum_J
\omega(J)\lambda_J}=0$, d'où une contradiction. La démonstration du fait
que les $\omega(J)$ engendrent
$A_{\star,\star}$ est plus délicate qu'on pourrait l'imaginer \emph{a
priori}. Soit $(i,j)\in\mathbf{Z}^2$. Montrons que tout élément de
$A_{i,j}$ est combinaison linéaire des $\omega(J)$. Notons $w=i-j$ et
$\mathcal B$ l'ensemble \emph{fini} des $I$ tels que $\theta(I)$ soit de
poids $\leq w$. Il résulte de la démonstration du
lemme~\ref{lemme-base-duale} que tout élément $\alpha\in A_{i,j}$ s'écrit
de manière unique $\alpha=\sum_{I\in\mathcal B} \theta(I)^\star \lambda_I$.

Comme $\mathcal B$ est fini, pour l'ordre induit, c'est un honnête ensemble
fini totalement ordonné. Cela a donc un sens de montrer par récurrence
descendante sur $I_0=\min\left\{I,\lambda_I\neq 0\right\}$
que tout élément (non nul)
$\alpha\in A_{i,j}$ comme ci-dessus est combinaison linéaire des
$\omega(J)$ pour $J\in\mathcal B$. On peut ainsi supposer que
$\alpha=\sum_{I\geq I_0,I\in\mathcal B} \theta(I)^\star \lambda_I\in
A_{i,j}$ et que pour tout $I\in\mathcal B$ tel que $I>I_0$, si $\lambda\in
H^{\star,\star}$ est un élément tel que $\theta(I)^\star\lambda\in
A_{i,j}$, alors $\theta(I)^\star \lambda$ est combinaison linéaire des
$\omega(J)$ pour $J\in\mathcal B$. Ceci permet bien sûr de supposer qu'en
fait $\alpha=\theta(I_0)^\star \lambda_{I_0}$. La
proposition~\ref{proposition-matrice-triangulaire} implique qu'il existe
d'uniques coefficients $\mu_I\in H^{\star,\star}$ pour $I>I_0$ (mais
\emph{a priori} $I$
n'est pas forcément dans $\mathcal B$ !) tels que
\[\pm \omega(I_0)=\theta(I_0)^\star+\sum_{I>I_0}
\theta(I)^\star\mu_I\;\text{,}\]
d'où
\[\alpha=\theta(I_0)^\star \lambda_{I_0}=\pm
\omega(I_0)\lambda_{I_0}-\sum_{I>I_0,I\in\mathcal B}\theta(I)^\star
\mu_I\lambda_{I_0}\;\text{.}\]
En effet, le choix de l'ensemble $\mathcal B$ fait que les termes
$\theta(I)^\star \mu_I\lambda_{I_0}$ sont nuls pour 
$I\not\in \mathcal B$. De là, on peut conclure que
$\alpha$ est combinaison linéaire des
$\omega(J)$, $J\in\mathcal B$ en utilisant l'hypothèse de récurrence.

\medskip

Nous pouvons maintenant démontrer le
théorème~\ref{theoreme-presentation-duale}. Notons
$\tilde{A}_{\star,\star}$ l'algèbre définie par la présentation
apparaissant dans l'énoncé. Les relations vérifiées plus haut montrent
qu'on dispose d'un morphisme évident
$\tilde{A}_{\star,\star}\to A_{\star,\star}$. Il s'agit de montrer que
c'est un isomorphisme. Si $\ell\neq 2$, il est évident que
$\tilde{A}_{\star,\star}$ est un module-$H^{\star,\star}$ libre ayant pour
base les monômes $\omega(I)$ (définis plus haut dans $A_{\star,\star}$,
mais ils ont bien entendu aussi un sens dans $\tilde{A}_{\star,\star}$).
Cette base de $\tilde{A}_{\star,\star}$ s'envoie
sur une base de $A_{\star,\star}$
d'après le corollaire~\ref{corollaire-base-duale}, donc
$\tilde{A}_{\star,\star}\to A_{\star,\star}$ est un isomorphisme.

Si $\ell=2$, la présentation a une forme telle que le fait que les éléments
$\omega(I)\in \tilde{A}_{\star,\star}$ constituent une famille libre ou
génératrice de $\tilde{A}_{\star,\star}$ comme module-$H^{\star,\star}$ n'a
rien de flagrant, algébriquement parlant. Notons $M\subset
\tilde{A}_{\star,\star}$ le sous-module-$H^{\star,\star}$ engendré par les
éléments $\omega(I)$. Le morphisme composé $M\to \tilde{A}_{\star,\star}\to
A_{\star,\star}$ est un isomorphisme en vertu du
corollaire~\ref{corollaire-base-duale}. Pour conclure, il suffit de montrer
que $M=\tilde{A}_{\star,\star}$, autrement dit que $M$ est stable par
produits. Comme $M$ est évidemment stable par multiplication par les
éléments de $H^{\star,\star}$ et les $\xi_k$, on se ramène à montrer que
tout produit $\tau_{i_1}\tau_{i_2}\dots \tau_{i_n}$ appartient à $M$. C'est
le cas par définition si les indices $i_k$ sont distincts. Dans le cas
général, le résultat s'obtient par récurrence sur $n$ en utilisant les
relations $\tau_k^2=\xi_{k+1}\tau + \tau_0\xi_{k+1}\rho +\tau_{k+1}\rho$
pour $k\geq 0$.

\subsection{La comultiplication $\Psi_\star$ sur $A_{\star,\star}$}

\begin{definition}
On note $A^{\star,\star}\otimes_{\rho,H^{\star,\star}} A^{\star,\star}$ le
produit tensoriel de $A^{\star,\star}$ vu comme module-$H^{\star,\star}$
par la multiplication à droite par $H^{\star,\star}\subset A^{\star,\star}$
(d'où la notation $\rho$ en indice) et de $A^{\star,\star}$ et de
$A^{\star,\star}$ considéré comme $H^{\star,\star}$-module comme on l'a
fait jusqu'à présent.
\end{definition}

En vérité, les deux copies de $A^{\star,\star}$ apparaissant dans ce
produit tensoriel sont toutes les deux munies de la structure de
$H^{\star,\star}$-bimodule-$H^{\star,\star}$ provenant de la multiplication
à gauche et à droite par $H^{\star,\star}\subset A^{\star,\star}$. Le
produit tensoriel
$A^{\star,\star}\otimes_{\rho,H^{\star,\star}} A^{\star,\star}$ hérite donc
ainsi d'une structure de $H^{\star,\star}$-bimodule-$H^{\star,\star}$.

\begin{definition}
On note $A_{\star,\star}\otimes_{H^{\star,\star},\lambda} A_{\star,\star}$
le produit tensoriel de deux copies de $A_{\star,\star}$ vu comme
$H^{\star,\star}$-bimodule-$H^{\star,\star}$ \emph{via} la multiplication à
gauche par $\lambda^\star(x)$ et à droite par $\xi_0x$ pour tout $x\in
H^{\star,\star}$ (cf.
remarque~\ref{remarque-deux-structures-de-module-sur-la-duale}).
\end{definition}

Ainsi, $A_{\star,\star}\otimes_{H^{\star,\star},\lambda} A_{\star,\star}$
est également un $H^{\star,\star}$-bimodule-$H^{\star,\star}$. (En tant que
groupe abélien, il ne dépend que de la structure de
module-$H^{\star,\star}$ sur le facteur de gauche et de la structure de
$H^{\star,\star}$-module sur le facteur de droite induite par
$\lambda^\star\colon H^{\star,\star}\to A^{\star,\star}$, ce qui est la raison
pour laquelle on a mis $\lambda$ en indice.)

\begin{definition}\label{definition-accouplement-double}
On définit un accouplement
$\left(A^{\star,\star}\otimes_{\rho,H^{\star,\star}} A^{\star,\star}\right)
\times
\left(A_{\star,\star}\otimes_{H^{\star,\star},\lambda}
A_{\star,\star}\right)\to H^{\star,\star}$ par les formules équivalentes
suivantes :
\[\acccro{C\otimes D}{\alpha\otimes
\beta}=\acc{C\acc{D}{\alpha}}{\beta}=\acc{C}{\acc{D}{\alpha}\beta}=
\acc{C}{\lambda^\star(\acc{D}{\alpha})\beta}\]
pour $C,D\in A^{\star,\star}$ et $\alpha,\beta\in A_{\star,\star}$.
\end{definition}

\begin{lemme}\label{lemme-dualite-double-tenseurs}
Soit $\varphi\colon
A^{\star,\star}\otimes_{\rho,H^{\star,\star}}A^{\star,\star}\to
H^{\star,\star}$ une application $H^{\star,\star}$-linéaire faisant
décroître le bidegré de $(i,j)\in \mathbf{Z}^2$. Alors, il existe un unique
élément $\psi\in
A_{\star,\star}\otimes_{H^{\star,\star},\lambda}A_{\star,\star}$ de bidegré
total $(i,j)$ tel que
$\varphi=\acccro{\cdot}{\psi}$.
\end{lemme}

Il est évident que les éléments $P^I\otimes P^J$ où $P^I$ et $P^J$
parcourent les monômes admissibles forment une base de
$A^{\star,\star}\otimes_{\rho,H^{\star,\star}} A^{\star,\star}$ comme
$H^{\star,\star}$-module et que $\acccro{P^I\otimes P^J}{\theta(J')^\star
\otimes \theta(I')^\star}=\delta_{II'}\delta_{JJ'}$. Un raisonnement
semblable à celui du lemme~\ref{lemme-base-duale} montre que $\varphi$
s'annule sur tous les $P^I\otimes P^J$ sauf un nombre fini, ce qui permet
de donner une formule pour $\psi$ comme combinaison linéaire d'un nombre
fini des $\theta(J')^\star \otimes \theta(I')^\star$ (qui forment une base
de $A_{\star,\star}\otimes_{H^{\star,\star},\lambda}A_{\star,\star}$ comme
module-$H^{\star,\star}$).

\begin{definition}
On définit une application $A_{\star,\star}\vers{\Psi_\star}
A_{\star,\star}\otimes_{H^{\star,\star},\lambda}A_{\star,\star}$ de façon à
ce que la formule suivante soit satisfaite pour $C,D\in A^{\star,\star}$ et
$\alpha\in A_{\star,\star}$ :
\[\acccro{C\otimes D}{\Psi_\star \alpha}=\acc{CD}{\alpha}\;\text{.}\]
\end{definition}

La composition dans $A^{\star,\star}$ induit une application
$H^{\star,\star}$-linéaire
$A^{\star,\star}\otimes_{\rho,H^{\star,\star}}A^{\star,\star}\to
A^{\star,\star}$ qui à $C\otimes D$ associe $CD$. L'existence et l'unicité
de $\Psi_\star \alpha$ pour tout $\alpha\in A_{\star,\star}$ résulte de
cette observation et du lemme~\ref{lemme-dualite-double-tenseurs}. La
construction $\Psi_\star$ reflète ainsi en quelque sorte la composition
dans $A^{\star,\star}$. Pour réduire le nombre de signes dans la
définition, nous avons utilisé une convention différente de celle de
\cite{voevodsky-reduced}. On ne s'étonnera donc pas que les tenseurs
apparaissant dans la proposition~\ref{proposition-calcul-psi-lower-star}
soient inversés par rapport à \cite[Lemma~12.11]{voevodsky-reduced}.

\begin{proposition}
L'application $\Psi_\star\colon
A_{\star,\star}\to 
A_{\star,\star}\otimes_{H^{\star,\star},\lambda}A_{\star,\star}$
est un morphisme d'anneaux, où
$A_{\star,\star}\otimes_{H^{\star,\star},\lambda}A_{\star,\star}$ est muni
de la structure évidente d'anneau commutatif au sens bigradué.
\end{proposition}

Outre l'identité évidente $\Psi_\star(\xi_0)=\xi_0\otimes \xi_0$,
il s'agit de montrer que pour tous $C,D\in A^{\star,\star}$ et
$\alpha,\beta\in A_{\star,\star}$, on a $\acccro{C\otimes D}{\Psi_\star
\alpha\cdot \Psi_\star \beta}=\acc{CD}{\alpha\beta}$. Ceci s'obtient sans
difficulté autre que celle de vérifier les signes en partant d'expressions
de $\Psi_\star\alpha$, $\Psi_\star \beta$, $\Psi^\star C$, $\Psi^\star D$
comme sommes de tenseurs et en déroulant les définitions des deux membres
jusqu'à obtenir une égalité.

\begin{remarque}
On montre facilement que $\Psi_\star\colon A_{\star,\star}\to 
A_{\star,\star}\otimes_{H^{\star,\star},\lambda}A_{\star,\star}$ est un
morphisme de $H^{\star,\star}$-bimodule-$H^{\star,\star}$. Compte tenu du
fait que $\Psi_\star$ soit un morphisme d'anneaux, cela revient à énoncer
que pour tout $\mu\in H^{\star,\star}$, on a d'une part
$\Psi_\star(\mu\xi_0)=\Psi^\star(\lambda^\star(\mu))
=\lambda^{\star}(\mu) \otimes \xi_0$ et d'autre part
$\Psi_\star(\xi_0\mu)=\xi_0\otimes \xi_0\mu$.
\end{remarque}

\begin{proposition}\label{proposition-calcul-psi-lower-star}
Pour tout $k\geq 0$, on a les égalités suivantes dans
$A_{\star,\star}\otimes_{H^{\star,\star},\lambda}A_{\star,\star}$ :
\[\Psi_{\star}\xi_k = \sum_{i=0}^k \xi_i\otimes \xi_{k-i}^{\ell^i}\qquad
\Psi_{\star}\tau_k = \xi_0\otimes \tau_k +
\sum_{i=0}^k \tau_i\otimes \xi_{k-i}^{\ell^i}\]
\end{proposition}

On se donne $C$ et $D$ des éléments de $A^{\star,\star}$. Calculons
$CD(v)\in H^{\star,\star}(\mathbf{P}^{d-1})$ pour tout $d\geq 1$.
La formule $\lambda^\star(v)=\sum_{i=0}^\infty\xi_i\otimes v^{\ell^i}$ du
corollaire~\ref{corollaire-lambda-u-et-v} implique que
$D(v)=\sum_{i=0}^\infty \acc{D}{\xi_i}v^{\ell^i}$. On a aussi :
\[\lambda^\star(v^{\ell^i})=\left(\lambda^\star(v)\right)^{\ell^i}=
\sum_{i=0}^\infty \xi_j^{\ell^i}\otimes v^{\ell^{i+j}}\;\text{,}\]
ce qui signifie que pour tout $F\in A^{\star,\star}$, on a :
$F(v^{\ell^i})=\sum_{j=0}^\infty
\acc{F}{\xi_j^{\ell^i}}v^{\ell^{i+j}}$. Ainsi :
\[CD(v)=\sum_{i=0}^\infty (C\acc{D}{\xi_i})(v^{\ell^i}) 
= \sum_{i,j\geq 0} \acc{C\acc{D}{\xi_i}}{\xi_j^{\ell^i}}v^{\ell^{i+j}}
= \sum_{i,j\geq 0}\acccro{C\otimes D}{\xi_i\otimes
\xi_j^{\ell^i}}v^{\ell^{i+j}}\;\text{.}\]
Par ailleurs, on a :
\[CD(v)=\sum_{k=0}^\infty \acc{CD}{\xi_k}v^{\ell^k}\;\text{.}\]
Pour $d$ assez grand, on peut identifier les coefficients devant
$v^{\ell^k}$ dans les deux expressions précédentes, ce qui donne l'identité
suivante dans $H^{\star,\star}$ :
\[\acc{CD}{\xi_k}=\sum_{i=0}^k\acccro{C\otimes D}{\xi_i\otimes
\xi_{k-i}^{\ell^i}}\;\text{.}\]
Ceci établit la première identité. La deuxième s'obtient de la même
façon en calculant $CD(u)$.

\subsection{La base de Milnor}

\begin{lemme}\label{lemme-bidualite}
Soit $\varphi\colon A_{\star,\star}\to H^{\star,\star}$ une application
linéaire-$H^{\star,\star}$ appliquant $A_{i,j}$ dans $H^{p-i,q-j}$ pour un
certain couple $(p,q)\in\mathbf{Z}^2$. Il existe un unique élément $F\in
A^{p,q}$ tel que pour tout $\alpha\in A_{\star,\star}$, on ait
$\varphi(\alpha)=\acc{F}{\alpha}$.
\end{lemme}

Ceci s'obtient par un simple raisonnement utilisant le poids, la base des
monômes admissibles de $A^{\star,\star}$ comme $H^{\star,\star}$-module et
la base duale de $A_{\star,\star}$ comme module-$H^{\star,\star}$ formée
des éléments $\theta(I)^\star$ (cf. lemme~\ref{lemme-base-duale}).

\medskip

Pour $(\varepsilon_\bullet)=(\varepsilon_0,\varepsilon_1,\dots)$ et
$(r_\bullet)=(r_1,r_2,\dots)$ des suites d'entiers comme dans la
définition~\ref{definition-omega}, on note
$\tau_\bullet^{\varepsilon_\bullet}=\tau_0^{\varepsilon_0}\tau_1^{\varepsilon_1}\dots$
et $\xi_\bullet^{r_\bullet}=\xi_1^{r_1}\xi_2^{r_2}\dots$, ce qui définit un
élément $\tau_\bullet^{\varepsilon_\bullet}\xi_\bullet^{r_\bullet}\in
A_{\star,\star}$. C'est la façon dont on notera ici les éléments de
$A_{\star,\star}$ qui étaient notés $\omega(I)$ dans la
définition~\ref{definition-omega} et qui en constituent une base comme
module-$H^{\star,\star}$ d'après le corollaire~\ref{corollaire-base-duale}.
On les désignera parfois sous le simple nom de \guil{monôme}.

\begin{definition}[Base de Milnor]
Pour $(\varepsilon_\bullet)$ et $(r_\bullet)$ comme ci-dessus, on note
$\rho(\varepsilon_\bullet,r_\bullet)$ l'élément de $A^{\star,\star}$
caractérisé grâce au lemme~\ref{lemme-bidualite} par le fait que
\[\acc{\rho(\varepsilon_\bullet,r_\bullet)}{\tau_\bullet^{\varepsilon'_\bullet}\xi_\bullet^{r'_\bullet}}\]
vaut $1$ si $\varepsilon_\bullet=\varepsilon'_\bullet$ et
$r_\bullet=r'_\bullet$ et $0$ sinon.
\end{definition}

Les éléments $\rho(\varepsilon_\bullet,r_\bullet)$ constituent la base de
Milnor de $A^{\star,\star}$ (comme $H^{\star,\star}$-module). Leur
définition a un sens puisque les éléments 
$\tau_\bullet^{\varepsilon_\bullet}\xi_\bullet^{r_\bullet}$ forment une
base de $A_{\star,\star}$ comme $H^{\star,\star}$-module (cf.
corollaire~\ref{corollaire-base-duale}). Que la famille
$\rho(\varepsilon_\bullet,r_\bullet)$ soit libre est évident et qu'elle
soit génératrice résulte de la formule
$F=\sum_{(\varepsilon_\bullet,r_\bullet)}
\acc{F}{\tau_\bullet^{\varepsilon_\bullet}\xi_\bullet^{r_\bullet}}\rho(\varepsilon_\bullet,r_\bullet)$
pour tout $F\in A^{\star,\star}$, une formule qui a un sens et qui est
vraie grâce à des considérations de poids semblables à celles qui
interviennent dans la démonstration du lemme~\ref{lemme-bidualite}.

\begin{definition}
Notons $I\subset A_{\star,\star}$ le sous-module-$H^{\star,\star}$ libre de
base les monômes
$\tau_\bullet^{\varepsilon_\bullet}\xi_\bullet^{r_\bullet}$ tels que
$(r_\bullet)\neq 0$. C'est aussi l'idéal (bilatère) de $A_{\star,\star}$
engendré par les éléments $\xi_k$ pour $k\geq 1$.
\end{definition}

\begin{definition}
On note $B^{\star,\star}\subset A^{\star,\star}$ le
sous-$H^{\star,\star}$-module formé des éléments $F\in A^{\star,\star}$
tels que pour tout $\alpha\in I$, on ait $\acc{F}{\alpha}=0$.
\end{definition}

\begin{proposition}\label{proposition-sous-algebre-b}
$B^{\star,\star}$ est une sous-$H^{\star,\star}$-algèbre de
$A^{\star,\star}$. En tant que $H^{\star,\star}$-module, $B^{\star,\star}$
est libre de base l'ensemble des éléments $\rho(\varepsilon_\bullet,0)$
(notés $Q(\varepsilon_\bullet)$) pour $\varepsilon_\bullet$ parcourant les
suites à valeurs dans $\{0,1\}$ nulles à partir d'un certain rang.
\end{proposition}

Le fait que les $Q(\varepsilon_\bullet)$ forment une base de
$B^{\star,\star}$ comme $H^{\star,\star}$-module est évident à partir de la
définition de $I$. Pour montrer que $B^{\star,\star}$ est un sous-anneau de
$A^{\star,\star}$, il convient en premier lieu d'observer que l'identité
appartient à $B^{\star,\star}$, ce qu'on obtient en faisant la vérification
(facile) de l'identité $Q(0,0,\dots)=\Id$ (utiliser la
proposition~\ref{proposition-matrice-triangulaire} avec $I=(0,\dots)$). Il s'agit en second lieu de
vérifier que si $C$ et $D$ appartiennent à $B^{\star,\star}$, alors $CD$
aussi. Il s'agit alors de montrer que pour tout $\alpha\in I$,
$\acc{CD}{\alpha}=\acccro{C\otimes D}{\Psi_\star \alpha}=0$. Pour
$\beta,\gamma\in A_{\star,\star}$, on a :
\[\acccro{C\otimes D}{\beta\otimes \gamma}=\acc{C\acc{D}{\beta}}{\gamma}=
\acc{C}{\lambda^\star(\acc{D}{\beta}){\gamma}}\;\text{.}\]
Cet accouplement est nul si $\beta$ ou $\gamma$ est dans $I$. C'est évident
si $\beta\in I$ ; si $\gamma\in I$, $\lambda^\star(\acc{D}{\beta})\gamma$
appartient aussi à $I$ puisque $I$ est un idéal et donc $C$ est bien
orthogonal à $\lambda^\star(\acc{D}{\beta})\gamma$.

Pour montrer que $\acccro{C\otimes D}{\Psi_\star\alpha}=0$, il suffit de
montrer que $\Psi_\star \alpha$ appartient au sous-groupe $I\otimes
A_{\star,\star}+A_{\star,\star}\otimes I$. Ce sous-groupe étant un idéal de
$A_{\star,\star}\otimes_{H^{\star,\star},\lambda}A_{\star,\star}$ et
$\Psi_\star$ étant un morphisme d'anneaux, il suffit d'obtenir cette
propriété pour $\alpha=\xi_k$ avec $k\geq 1$, ce que fournit la
proposition~\ref{proposition-calcul-psi-lower-star}. Ceci achève la
démonstration de la proposition et incite à procéder à la
définition suivante :

\begin{definition}
On note $\overline{\Psi}_\star\colon A_{\star,\star}/I\to
(A_{\star,\star}/I)\otimes_{H^{\star,\star},\lambda}(A_{\star,\star}/I)$ le
morphisme induit par $\Psi_\star\colon A_{\star,\star}\to
A_{\star,\star}\otimes_{H^{\star,\star},\lambda}A_{\star,\star}$.
\end{definition}

L'accouplement de la définition~\ref{definition-accouplement-double} induit
un accouplement noté pareillement :
\[\acccro{\cdot}{\cdot}\colon \left(B^{\star,\star}
\otimes_{\rho,H^{\star,\star}} B^{\star,\star}\right)\times
\left((A_{\star,\star}/I)\otimes_{H^{\star,\star},\lambda}(A_{\star,\star}/I)\right)
\to H^{\star,\star}\;\text{.}\]

\begin{definition}
Pour toute partie finie $X\subset \mathbf{N}$, on note
$Q(X)=Q(\varepsilon_\bullet)$ et
$\tau(X)=\tau_\bullet^{\varepsilon_\bullet}$ où
$\varepsilon\colon\mathbf{N}\to\{0,1\}$ est la fonction caractéristique de
$X$. Pour tout $i\in \mathbf{N}$, on note $Q_i=Q(\{i\})\in
A^{2(\ell^i-1)+1,\ell^i-1}$.

Si $\ell=2$, on note aussi $Q(n)=Q(\varepsilon_\bullet)$ et
$\tau(n)=\tau_\bullet^{\varepsilon_\bullet}$ quand $n=\sum_{i=0}^\infty
\varepsilon_i 2^i$. En particulier, $Q(2^i)=Q_i$.
\end{definition}

\begin{proposition}
Pour toute partie finie $X\subset \mathbf{N}$, on a :
\[\overline{\Psi}_\star \tau(X)=\sum_{A\sqcup
B=X}(-1)^{\#\Inv(A,B)}\tau(A)\otimes \tau(B)\;\text{,}\]
où la somme est indexée par l'ensemble des partitions de $X$ en deux
parties $A$ et $B$ et où $\Inv(A,B)=\left\{(a,b)\in A\times B,
a>b\right\}$.
\end{proposition}

La formule est exacte si $X$ est vide ou si $X=\{i\}$ car
$\overline{\Psi}_\star \tau_i=\xi_0\otimes \tau_i+\tau_i\otimes \xi_0$ (cf.
proposition~\ref{proposition-calcul-psi-lower-star}). On peut conclure par
récurrence sur le cardinal de $X$ puisque l'on vérifie facilement que si la
formule est vrai pour une partie finie $X$, elle l'est aussi pour $X\sqcup
\{i\}$ pour tout $i>\max X$.

\begin{corollaire}
Si $A$ et $B$ sont deux parties finies de $\mathbf{N}$, alors 
\[Q(B)Q(A)=\left\{\begin{array}{ll} (-1)^{\#Inv(A,B)}Q(A\sqcup B)
&\text{si $A$ et $B$ sont disjointes}\\
0& \text{sinon}\end{array}\right.\]
En particulier, $Q_i^2=Q_iQ_i=0$ pour tout $i\geq 0$.
\end{corollaire}

Dans $B^{\star,\star}$, on peut écrire \[Q(B)Q(A)=\sum_{\overset{X\subset
\mathbf{N}}{X\text{ finie}}}\lambda_X Q(X)\] avec
$\lambda_X=\acc{Q(B)Q(A)}{\tau(X)}\in
H^{\star,\star}$.
On a :
\[\lambda_X = \acccro{Q(B)\otimes Q(A)}{\overline{\Psi}_\star \tau(X)}=
\sum_{A'\sqcup B'=X}(-1)^{\#\Inv(A',B')}\acccro{Q(B)\otimes
Q(A)}{\tau(A')\otimes \tau(B')}\]
Il vient aussitôt que :
\[\lambda_X=\left\{\begin{array}{ll}(-1)^{\#\Inv(A,B)}&
\text{si }X=A\sqcup B \\ 0 & \text{sinon}\end{array}\right.\]
Ceci permet de conclure.

\begin{corollaire}
La sous-$\Fl$-algèbre de $B^{\star,\star}$ engendrée par les opérations
$Q_i$ pour $i\geq 0$ est une algèbre alternée.
\end{corollaire}

Le lecteur intéressé par les signes prendra garde au fait que si
$X=\{i_1<\dots<i_k\}$ alors $Q(X)=Q_{i_k}\dots Q_{i_1}$ (qui ne coïncide
qu'au signe près avec $Q_{i_1}\dots Q_{i_k}$).

\begin{proposition}
Dans $A^{1,0}$, on a : 
$Q_0=\beta$.
\end{proposition}

Pour des raisons de poids, la décomposition de $\beta$ sur la base de
Milnor est de la forme
$\beta=\acc{\beta}{\xi_0}\Id+\acc{\beta}{\tau_0}Q_0=\acc{\beta}{\tau_0}Q_0$
car $\acc{\beta}{\xi_0}=\beta(1)=0$. La formule pour $\lambda^\star(u)$ du
corollaire~\ref{corollaire-lambda-u-et-v} implique que $\beta
u=\acc{\beta}{\tau_0}v$. Comme on sait que $\beta u=v$, on en déduit que
$\acc{\beta}{\tau_0}=1$ et donc que $\beta=Q_0$.

\begin{proposition}\label{proposition-comultiplication-b-star-star}
La comultiplication $\Psi^\star\colon A^{\star,\star}\to
A^{\star,\star}\otimes_{H^{\star,\star}} A^{\star,\star}$ induit une
comultiplication $B^{\star,\star}\to
B^{\star,\star}\otimes_{H^{\star,\star}} B^{\star,\star}$.
\end{proposition}

Il est évident que $B^{\star,\star}\otimes_{H^{\star,\star}}
B^{\star,\star}$ s'identifie à l'orthogonal de $I\otimes
A_{\star,\star}+A_{\star,\star}\otimes I$ pour l'accouplement de la
définition~\ref{definition-accouplement-double-trivial}. Pour montrer que
$\Psi^\star$ induit une comultiplication sur $B^{\star,\star}$, il suffit
donc de montrer que si $C\in B^{\star,\star}$ (c'est-à-dire que $C$ est
orthogonal à $I$), alors pour tous $\beta$ et $\gamma$ dans
$A_{\star,\star}$, on a $\acc{\Psi^{\star}C}{\beta\otimes \gamma}=0$ si
$\beta$ ou $\gamma$ est dans $I$. Comme $\acc{\Psi^{\star}C}{\beta\otimes
\gamma}=\acc{C}{\beta\gamma}$, il suffit pour cela d'observer que si $\beta$ ou
$\gamma$ est dans $I$, alors $\beta\gamma\in I$ aussi puisque $I$ est un
idéal.

\begin{proposition}
\begin{itemize}
\item Si $\ell\neq 2$, $\Psi^\star Q_i=Q_i\otimes \Id + \Id\otimes Q_i$
pour $i\geq 0$.
\item Si $\ell=2$, pour tout $n\geq 0$,
$\Psi^\star Q(n)=\sum_{p+q=n}\rho^{\sigma(p,q)}Q(p)\otimes Q(q)$ où l'on
note $\sigma(p,q)$ le nombre de retenues qui interviennent lorsque l'on pose
l'addition de $p$ et $q$ en base $2$. En
particulier, pour $n=2^i$ :
\[\Psi^\star Q_i= Q_i\otimes \Id+\Id\otimes Q_i +
\sum_{k=1}^{2^i-1}\rho^{i-v_2(k)} Q(k)\otimes Q(2^i-k)\;\text{,}\]
où $v_2$ désigne la valuation $2$-adique.
\end{itemize}
\end{proposition}

Soit $X$ une partie finie de $\mathbf{N}$. D'après la
proposition~\ref{proposition-comultiplication-b-star-star}, on peut écrire
\[\Psi^\star Q(X)=\sum_{(A,B)}\lambda_{A,B}Q(A)\otimes Q(B)\;\text{,}\]
où $(A,B)$ parcourt les couples de parties finies de $\mathbf{N}$. On a :
\[\lambda_{A,B}=(-1)^{\#A\cdot \#B}\acc{\Psi^\star Q(X)}{\tau(A)\otimes
\tau(B)}=(-1)^{\#A\cdot \#B} \acc{Q(X)}{\tau(A)\tau(B)}\;\text{.}\]

Si $\ell\neq 2$ et $X=\{i\}$, pour qu'un coefficient $\lambda_{A,B}$ soit
non nul, il faut que $\tau(A)\tau(B)$ soit non nul, c'est-à-dire que $A$ et
$B$ soient disjointes auquel cas
$\tau(A)\tau(B)=(-1)^{\#\Inv(A,B)}\tau(A\sqcup B)$ et l'accouplement
$\acc{Q(X)}{\tau(A\sqcup B)}$ ne peut-être non nul que si $A\sqcup
B=\{i\}$. On obtient facilement que $\lambda_{\{i\},\emptyset}=
\lambda_{\emptyset,\{i\}}=1$, ce qui permet de conclure.

Si $\ell=2$, en utilisant le lemme suivant, la même méthode permet
d'obtenir la formule générale pour $\Psi^\star Q(n)$ pour tout $n\geq 0$.

\begin{lemme}
Si $\ell=2$, alors pour tous $p,q\geq 0$, on a l'identité
$\tau(p)\tau(q)=\tau(p+q)\rho^{\sigma(p,q)}$ dans le quotient
$A_{\star,\star}/I$.
\end{lemme}

Si on raisonne modulo $I$, le théorème~\ref{theoreme-presentation-duale}
fournit l'identité $\tau_i^2=\tau_{i+1}\rho$ ce qui est la compatibilité
cherchée dans le cas particulier $p=q=2^i$. Il est par ailleurs
évident que si l'addition de $p$ et de $q$ en base $2$ ne provoque aucune
retenue, on a $\tau(p)\tau(q)=\tau(p+q)$. Le cas général résulte de ces
deux cas particuliers.

\begin{definition}
Les éléments $\rho((0,0,\dots),r_\bullet)$ de la base de Milnor
correspondant aux couples de suites $(\varepsilon_\bullet,r_\bullet)$ tels
que $\varepsilon_i=0$ pour tout $i\geq 0$ sont notés $\mathscr
P^{(r_\bullet)}$.
\end{definition}

\begin{proposition}\label{proposition-q-p}
Pour tout couple $(\varepsilon_\bullet),(r_\bullet)$ de suites d'entiers
nulles à partir d'un certain rang, $\varepsilon_\bullet$ n'étant constituée
que de $0$ et de $1$, on a l'identité suivante dans $A^{\star,\star}$ :
\[\rho(\varepsilon_\bullet,r_\bullet)=Q(\varepsilon_\bullet)\mathscr
P^{(r_\bullet)}\;\text{.}\]
\end{proposition}

Le coefficient $\lambda_{(\varepsilon'_\bullet,r'_\bullet)}$
devant $\rho(\varepsilon',r'_\bullet)$ dans la décomposition de
$Q(\varepsilon_\bullet)\mathscr P^{(r_\bullet)}$ sur la base de Milnor est
\[\lambda_{(\varepsilon'_\bullet,r'_\bullet)}
=\acccro{Q(\varepsilon_\bullet)\otimes
\mathscr
P^{(r_\bullet)}}{\Psi_\star(\tau_\bullet^{\varepsilon'_\bullet}
\xi_\bullet^{r'_\bullet})}\]
Comme $Q(\varepsilon_\bullet)$ est orthogonal à l'idéal $I$ de
$A_{\star,\star}$, il suffit de calculer l'image
$\tilde{\Psi}_\star(\tau_\bullet^{\varepsilon'_\bullet
}\xi_\bullet^{r'_\bullet})$ de
$\Psi_\star(\tau_\bullet^{\varepsilon'_\bullet
}\xi_\bullet^{r'_\bullet})$ de
dans l'anneau quotient
$A_{\star,\star}\otimes_{H^{\star,\star},\lambda}(A_{\star,\star}/I)$
puisque l'accouplement $\acccro{Q(\varepsilon_\bullet)\otimes\mathscr
P^{(r_\bullet)}}{\cdot}$
se factorise par ce quotient (l'argument est le même que
dans la démonstration de la proposition~\ref{proposition-sous-algebre-b}).
Grâce à la proposition~\ref{proposition-calcul-psi-lower-star}, on dispose
des formules suivantes :
\[\tilde{\Psi}_\star \xi_k = \xi_k\otimes \xi_0\qquad
\tilde{\Psi}_\star \tau_k = \xi_0\otimes \tau_k+\tau_k\otimes\xi_0\]
Si on note $I'=\{i\in\mathbf{N},\varepsilon'_i=1\}$, on en déduit :
\[\tilde{\Psi}_\star
(\tau_\bullet^{\varepsilon'_\bullet}\xi_\bullet^{r'_\bullet})=
\sum_{A\sqcup B=I'}
(-1)^{\#\Inv(A,B)}\tau(A)\xi_\bullet^{r'_\bullet}\otimes \tau(B) \]
Ainsi,
\begin{eqnarray*}
\lambda_{(\varepsilon'_\bullet,r'_\bullet)} &=& \sum_{A\sqcup B=I'}
(-1)^{\#\Inv(A,B)}\acccro{Q(\varepsilon_\bullet)\otimes \mathscr
P^{(r_\bullet)}}{\tau(A)\xi_\bullet^{r'_\bullet}\otimes \tau(B)}\\
&=& \sum_{A\sqcup B=I'}
(-1)^{\#\Inv(A,B)}\acc{Q(\varepsilon_\bullet)\acc{\mathscr
P^{(r_\bullet)}}{\tau(A)\xi_\bullet^{r'_\bullet}}}{\tau(B)}
\end{eqnarray*}
On voit immédiatement que seul le terme correspondant à $A=\emptyset$ (et
donc $B=I'$) peut
contribuer à cette somme. Ensuite, il vient que le résultat
$\lambda_{(\varepsilon'_\bullet,r'_\bullet)}$ ne peut être non
nul que si $r'_\bullet=r_\bullet$ et
$\varepsilon'_\bullet=\varepsilon_\bullet$ auquel cas le coefficient vaut
bien $1$.

\begin{proposition}\label{proposition-q-n-commutateur}
Pour tout $n\geq 1$, on a $Q_n=[q_n,\beta]=q_n\beta-\beta q_n$,
où $q_n\in A^{\star,\star}$ est l'élément de la base de Milnor
$q_n=\mathscr P^{(r_\bullet)}$ défini par
$r_i=\delta_{in}$ (autrement dit, $q_n$
est dual du monôme $\xi_n$).
\end{proposition}

Il s'agit de décomposer $q_n\beta$ sur la base de Milnor à laquelle
appartiennent $Q_n$ et $\beta q_n$ (cf. proposition~\ref{proposition-q-p}).
Pour tout $(\varepsilon_\bullet,r_\bullet)$, le coefficient
$\lambda_{(\varepsilon_\bullet,r_\bullet)}$ devant
$\rho(\varepsilon_\bullet,r_\bullet)$ dans la décomposition de $q_n\beta$
est :
\[\lambda_{(\varepsilon_\bullet,r_\bullet)}=
\acccro{q_n\otimes\beta}{\Psi_\star(\tau_\bullet^{\varepsilon_\bullet}
\xi_\bullet^{r_\bullet})}\;\text{.}\]
Notons $J$ l'idéal de $A_{\star,\star}$ engendré par $\tau_k$ pour $k\geq
1$. L'accouplement $\acccro{q_n\otimes \beta}{\cdot}$ se factorise par
l'anneau quotient $(A_{\star,\star}/I+J)\otimes_{H^{\star,\star},\lambda}
A_{\star,\star}$. Notons $\breve{\Psi}_\star\colon A_{\star,\star}\to
(A_{\star,\star}/I+J)\otimes_{H^{\star,\star},\lambda} A_{\star,\star}$
le morphisme induit par $\Psi_\star$. D'après la
proposition~\ref{proposition-calcul-psi-lower-star}, on a :
\[\breve{\Psi}^\star \xi_k=\xi_0\otimes \xi_k\qquad
\breve{\Psi}^\star \tau_k = \xi_0\otimes \tau_k + \tau_0\otimes \xi_k\]
L'anneau quotient $A_{\star,\star}/I+J$ étant un module-$H^{\star,\star}$
libre de base $(\xi_0,\tau_0)$ (avec la relation $\tau_0^2=0$), on obtient
aisément que si
$X=\{i\in\mathbf{N},\varepsilon_i=1\}=\{x_1<x_2<\dots<x_k\}$, alors :
\[\breve{\Psi}^\star(\tau(X))=\xi_0\otimes \tau(X) + \sum_{i=1}^k
(-1)^{i-1} \tau_0\otimes \tau(X-\{x_i\})\xi_{x_i}\;\text{, d'où}\]
\[\breve{\Psi}^\star(\tau_\bullet^{\varepsilon_\bullet}
\xi_\bullet^{r_\bullet})=\xi_0\otimes \tau_\bullet^{\varepsilon_\bullet}
\xi_\bullet^{r_\bullet} + \sum_{i=1}^k
(-1)^{i-1} \tau_0\otimes \tau(X-\{x_i\})\xi_{x_i}
\xi_\bullet^{r_\bullet}\;\text{.}\]
Il en résulte que :
\[\lambda_{(\varepsilon_\bullet,r_\bullet)} = \sum_{i=1}^k (-1)^{i-1}
\acc{q_n}{\tau(X-\{x_i\})\xi_{x_i}\xi_\bullet^{r_\bullet}}\;\text{.}\]
Pour que ce coefficient $\lambda_{(\varepsilon_\bullet,r_\bullet)}$ soit
non nul, il faut donc que $X$ soit un singleton. Dans ce cas, notons $x$
l'unique entier naturel appartenant à $X=\{i,\varepsilon_i=1\}$. On a alors :
\[\lambda_{(\varepsilon_\bullet,r_\bullet)}=\acc{q_n}{\xi_{x}
\xi_\bullet^{r_\bullet}}\;\text{.}\]
L'élément $\xi_{x} \xi_\bullet^{r_\bullet}$ est égal à un monôme en les
$\xi_j$ pour $j\geq 1$ (mais l'unité $\xi_0$ peut se cacher dans cette
expression si $x=0$). Le coefficient
$\lambda_{(\varepsilon_\bullet,r_\bullet)}$ vaut donc $1$ si ce monôme
$\xi_{x} \xi_\bullet^{r_\bullet}$ est $\xi_n$ et $0$ sinon. Le coefficient
$\lambda_{(\varepsilon_\bullet,r_\bullet)}$ vaut $1$ dans les deux cas
suivants :
\begin{itemize}
\item $x=0$ et $\xi_\bullet^{r_\bullet}=\xi_n$, de sorte que
$\tau_\bullet^{\varepsilon_\bullet}\xi_\bullet^{r_\bullet}=\tau_0\xi_n$ et
l'élément correspondant de la base de Milnor est $\beta q_n$ d'après la
proposition~\ref{proposition-q-p}.
\item $x=n$ et $\xi_\bullet^{r_\bullet}=\xi_0$, de sorte que
$\tau_\bullet^{\varepsilon_\bullet}\xi_\bullet^{r_\bullet}=\tau_n$ et
l'élément correspondant de la base de Milnor est $Q_n$.
\end{itemize}
On a ainsi établi l'identité $q_n\beta=Q_n+\beta q_n$.

\begin{proposition}
Pour tout entier $n\geq 0$, l'opération $P^n$ appartient à la base de
Milnor : $P^n=\mathscr P^{(n,0,0,\dots)}$.
\end{proposition}

Il s'agit de montrer que $P^n$ est dual du monôme $\xi_1^n$. D'après la
proposition~\ref{proposition-matrice-triangulaire}, on a
$\acc{P^n}{\tau_\bullet^{\varepsilon_\bullet}\xi_\bullet^{r_\bullet}}=0$ si
$(0,n,0,\dots)<(\varepsilon_0,r_1,\varepsilon_1,r_2,\dots)$ pour l'ordre
lexicographique pour les suites lues de droite à gauche. Les seules suites
$(\varepsilon_0,r_1,\varepsilon_2,\dots)$ qu'il reste à considérer sont
celles de la forme $(0,k,0,\dots)$ avec $k\leq n$ et $(1,k,0,\dots)$ avec
$k<n$. En examinant les bidegrés des éléments de $H^{\star,\star}$ de la
forme $\acc{P^n}{\xi_1^k}$ ou $\acc{P^n}{\tau_0^{\varepsilon_0}\xi_1^k}$,
on peut éliminer un certain nombre de cas. Il ne reste alors plus qu'à
montrer que $\acc{P^n}{\xi_1^n}=1$ et que $\acc{P^n}{\tau_0\xi_1^{n-1}}=0$,
ce qui se démontre aisément par récurrence sur $n$ en utilisant les
formules pour $\Psi^\star P^n$ et les identités suivantes :
\[\acc{P^n}{\xi_1^n}=\acc{\Psi^\star P^n}{\xi_1^{n-1}\otimes \xi_1}\qquad
\acc{P^n}{\tau_0\xi_1^{n-1}}=\acc{\Psi^\star P^n}{\xi_1^{n-1}\otimes
\tau_0}\;\text{.}\]

\subsection{Action sur les classes de Chern et les classes de Thom}

On se propose de décrire l'action des éléments
$\rho(\varepsilon_\bullet,r_\bullet)=Q(\varepsilon_\bullet)\mathscr
P^{(r_\bullet)}$ de la base de Milnor de $A^{\star,\star}$ sur les classes
de Chern $c_i(V)\in H^{2i,i}(X)$ et de Thom $t_V\in \tH^{2r,r}(\Th_X V)$
pour des fibrés vectoriels $V$ de rang $r$ sur $X\in\Sm[k]$. Une évidence à
remarquer dès maintenant est que si $(\varepsilon_\bullet)\neq 0$,
$\rho(\varepsilon_\bullet,r_\bullet)$ s'annule sur toutes ces classes
puisque l'image appartient à un groupe de la forme $H^{2w+e,w}(X)$ ou
$\tH^{2w+e,w}(\Th_X V)$ avec $e=\sum_{i}\varepsilon_i>0$ et qu'un tel groupe
est nul. Le seul cas
intéressant à étudier est celui de l'action des opérations $\mathscr
P^{(r_\bullet)}$ sur ces classes.

\begin{proposition}\label{proposition-lambda-c-1}
Soit $X\in\Sm[k]$. Soit $L$ un fibré en droites sur $X$. Alors :
\[\lambda^\star(c_1(L))=\sum_{k=0}^\infty \xi_k\otimes c_1(L)^{\ell^k}\in
A_{\star,\star}\otimes_{H^{\star,\star}} H^{\star,\star}(X)\;\text{.}\]
\end{proposition}

La formule est déjà connue pour $X=\mathbf{P}^{n-1}$ (pour $n\geq 1$) et
$c_1(\OO(1))$. Le cas général vient de ce que pour $n$ assez grand, il
existe un morphisme $f\colon X\to\mathbf{P}^{n-1}$ dans $\Ho[k]$ tel que
$f^\star(c_1(\OO(1)))=c_1(L)$. Si $L$ est engendré par $n$ sections
globales (par exemple si $X$ est affine et $n$ assez grand), on peut en
effet trouver un morphisme $f\colon X\to \mathbf{P}^{n-1}$ dans $\Sm[k]$
tel que $f^\star(\OO(1))\simeq L$. Le cas général s'en déduit par l'astuce
de Jouanolou.

\begin{proposition}\label{proposition-action-p-r-sur-classes-de-chern}
Soit $(r_\bullet)=(r_1,r_2,\dots)$ une suite d'entiers naturels nulle à
partir d'un certain rang. Soit $i\geq 0$. Soit $d$ un entier naturel. On
note $P\in\Fl[X_1,\dots,X_d]$ le polynôme symétrique
\[P=\sum_{\underset{\#J=i}{J\subset\{1,\dots,d\}}}\sum_{\underset{\prod_{j\in
J}\xi_{k_j}=\xi_\bullet^{r_\bullet}}{k\colon J\to \mathbf{N}}} \prod_{j\in
J} X_j^{\ell^{k_j}}\;\text{.}\]
On note $R\in \Fl[C_1,\dots,C_d]$ l'unique polynôme tel que si on note
$S_j\in \Fl[X_1,\dots,X_d]$ la $j$-ième fonction symétrique élémentaire des
variables $X_1,\dots,X_d$, on ait $P=R(S_1,\dots,S_d)$. Alors, pour tout
fibré vectoriel $V$ sur $X\in\Sm[k]$ de rang $\leq d$, on a :
\[\mathscr P^{(r_\bullet)}(c_i(V))=R(c_1(V),\dots,c_d(V))\in
H^{\star,\star}(X)\;\text{.}\]
\end{proposition}

(Pour des raisons de degrés évidentes, le polynôme $R\in
\Fl[S_1,\dots,S_d]$ se stabilise pour $d$ assez grand, $r_\bullet$ et $i$
étant fixés.)

D'après le principe de scindage, si $V$ est un fibré vectoriel de rang $d$,
il existe un morphisme $p\colon Y\to X$ tel que $p^\star V$ soit isomorphe
à une somme directe de $d$ fibrés en droites et tel que $p^\star\colon
H^{\star,\star}(X)\to H^{\star,\star}(Y)$ soit injectif. On peut donc
supposer que $V=L_1\oplus \dots\oplus L_d$ où pour tout $j\in
\{1,\dots,d\}$, $L_j$ est un fibré en droites sur $X$. On note
$x_j=c_1(L_j)\in H^{\star,\star}(X)$. On a :
\[c_i(V)=S_i(x_1,\dots,x_d)=\sum_{\underset{\#J=i}{J\subset\{1,\dots,d\}}}
\prod_{j\in J}x_j\in H^{\star,\star}(X)\;\text{.}\]
Comme $\lambda^\star(x_j)=\sum_{k=0}^\infty \xi_k\otimes x_j^{\ell^k}$ et
que $\lambda^\star$ est un morphisme d'anneaux, on obtient :
\[\lambda^{\star}(c_i(V))=\sum_{\underset{\#J=i}{J\subset\{1,\dots,d\}}}
\sum_{k\colon J\to\mathbf{N}} \prod_{j\in J}\left(\xi_{k_j}\otimes x_j^
{\ell^{k_j}}\right)\;\text{.}\]
Par conséquent :
\[\mathscr P^{(r_\bullet)}(c_i(V))=\acc{\mathscr
P^{(r_\bullet)}}{\lambda^\star(c_i(V))}=P(x_1,\dots,x_d)=
R(c_1(V),\dots,c_d(V))\;\text{.}\]

\begin{proposition}\label{proposition-action-p-r-sur-classes-de-thom}
Soit $(r_\bullet)=(r_1,r_2,\dots)$ une suite d'entiers naturels nulle à
partir d'un certain rang. Soit $d$ un entier naturel. On
note $P\in\Fl[X_1,\dots,X_d]$ le polynôme symétrique
\[P=\sum_{\underset{\xi_{k_1}\dots\xi_{k_d}=
\xi_\bullet^{r_\bullet}}{(k_1,\dots,k_d)\in\mathbf{N}^d}} \prod_{j=1}^d
X_j^{\ell^{k_j}-1}\;\text{.}\]
On note $R\in \Fl[C_1,\dots,C_d]$ l'unique polynôme tel que si on note
$S_j\in \Fl[X_1,\dots,X_d]$ la $j$-ième fonction symétrique élémentaire des
variables $X_1,\dots,X_d$, on ait $P=R(S_1,\dots,S_d)$. Alors, pour tout
fibré vectoriel $V$ sur $X\in\Sm[k]$ de rang $\leq d$, on a :
\[\mathscr P^{(r_\bullet)}(t_V)=R(c_1(V),\dots,c_d(V))\cdot t_V\in
H^{\star,\star}(\Th_X V)\;\text{.}\]
\end{proposition}

(L'entier $i$ et la suite $(r_\bullet)$ étant fixés, le polynôme $R$ se
stabilise dès que $d\geq \sum_{i} (\ell^i-1)r_i$.)

La démonstration est essentiellement la même que celle de la
proposition~\ref{proposition-action-p-r-sur-classes-de-chern}. Il suffit de
remplacer dans le raisonnement l'expression de $\lambda^{\star}(c_1(L))$
pour $L$ un fibré en droites par la formule du lemme suivant :

\begin{lemme}
Soit $X\in\Sm[k]$. Soit $L$ un fibré en droites sur $X$. Alors :
\[\lambda^\star(t_L)=\sum_{k=0}^\infty \xi_k\otimes c_1(L)^{\ell^k-1}t_L\in
A_{\star,\star}\otimes_{H^{\star,\star}}
\tH^{\star,\star}(\Th_X L)\;\text{.}\]
\end{lemme}

On utilise le plongement canonique de $\tH^{\star,\star}(\Th_X L)$ dans
$H^{\star,\star}(\mathbf{P}(V\oplus \OO_X))$. On a alors
$t_L=c_1(\OO(1))+c_1(L)$. En appliquant la
proposition~\ref{proposition-lambda-c-1} aux fibrés en droites $\OO(1)$ et
$L$, on obtient :
\[\lambda^\star(t_L)=\sum_{k=0}^\infty \xi_k\otimes (c_1(\OO(1))^{\ell^k}
+c_1(L)^{\ell^k})\;\text{.}\]
Il ne reste plus qu'à montrer que 
\[c_1(\OO(1))^{\ell^k}
+c_1(L)^{\ell^k}=c_1(L)^{\ell^k-1}c_1(\OO(1))+
c_1(L)^{\ell^k}=c_1(L)^{\ell^k-1}t_L\;\text{,}\]
ce qui se déduit facilement de l'identité
$c_1(\OO(1))^2+c_1(L)c_1(\OO(1))=0$ dans
$H^{\star,\star}(\mathbf{P}(V\oplus\OO_X))$.

\begin{corollaire}
Pour tout $j\geq 0$, on note $s_j\colon K_0(X)\to \oplus_i H^{2i,i}(X)$ la
transformation naturelle additive pour $X\in \Sm[k]$ telle que pour tout
fibré en droites $L$, $s_j([L])=c_1(L)^j$. Alors, pour tout $n\geq 1$ et
tout fibré vectoriel $V$ sur $X\in\Sm[k]$, on a :
\[q_n(t_V)=s_{\ell^n-1}(V)\cdot t_V\]
où $q_n\in A^{\star,\star}$ est l'élément dual du monôme $\xi_n$ introduit
dans la proposition~\ref{proposition-q-n-commutateur}.
\end{corollaire}
(En effet, avec les notations de la
proposition~\ref{proposition-action-p-r-sur-classes-de-thom}, on a
$P=\sum_{i=1}^d X_i^{\ell^n-1}$.)

\section{Endomorphismes du spectre représentant la cohomologie motivique}
\label{section-endomorphismes}

Soit $\ell$ un nombre premier. Soit $k$ un corps parfait. Notons $\HZl\in
\SH$ le spectre représentant la cohomologie motivique à coefficients
$\mathbf{Z}/\ell\mathbf{Z}$ (cf. \cite[\S{}6.1]{voevodsky-icm}). Il est
constitué d'une suite d'espaces $K(\mathbf{Z}/\ell\mathbf{Z}(n),2n)$ qui
sont exactement les espaces $K_n$ de la
définition~\ref{definition-eilenberg-maclane} (pour
$\Lambda=\mathbf{Z}/\ell\mathbf{Z}$). Les morphismes d'assemblage
$(\mathbf{A}^1/\mathbf{A}^1-\{0\})\wedge K_n \to K_{n+1}$ sont donnés par
la multiplication par les morphismes tautologiques
$\tau_{\mathbf{A}^1}\colon \mathbf{A}^1/\mathbf{A}^1-\{0\}\to K_1$ (cf.
proposition~\ref{proposition-classes-tautologiques}).

Pour ainsi dire par définition des opérations cohomologiques stables et de
la catégorie $\SHnaive$ de \cite[\S{}6]{riou-sh}, le groupe des opérations
cohomologiques stables de bidegré $(p,q)$ sur la cohomologie motivique à
coefficients $\mathbf{Z}/\ell\mathbf{Z}$ (cf.
définition~\ref{definition-operation-stable}) s'identifie au groupe
\[\Hom_{\SHnaive}(\HZl,\HZl\wedge S^{p,q})\]
où $S^{p,q}=(\mathbf{P}^1)^{q}\wedge S^{p-2q}\in \SH$ (ce qui a un sens
quels que soient $p$ et $q$).

\medskip

Le théorème principal de \cite{voevodsky-eilenberg-maclane} est alors le
suivant :

\begin{theoreme1}[Voevodsky \cite{voevodsky-eilenberg-maclane}]
\label{theoreme-voevodsky}
Soit $\ell$ un nombre premier. Soit $k$ un corps de caractéristique zéro.
L'inclusion de $A^{\star,\star}$ dans l'algèbre des opérations
cohomologiques stables (cf.
proposition~\ref{proposition-plongement-a-star-star}) est un isomorphisme.

Autrement dit, pour tout $(p,q)\in\mathbf{Z}^2$, le morphisme canonique
\[A^{p,q}\to \Hom_{\SHnaive}(\HZl,\HZl\wedge S^{p,q})\]
est un isomorphisme.
\end{theoreme1}

Le foncteur évident $\SH\to \SHnaive$ est plein (cf.
\cite[Proposition~6.3]{riou-sh}). En particulier, les applications
\[\Hom_{\SH}(\HZl,\HZl\wedge S^{p,q})\to \Hom_{\SHnaive}(\HZl,\HZl\wedge
S^{p,q})\;\text{.}\] sont surjectives. Un élément du noyau d'une telle
application serait un morphisme \guil{stablement fantôme} $f\colon
\HZl\to\HZl\wedge S^{p,q}$. Ceci signifierait que pour tout $\mathcal X\in
\Ho[k]$, l'application \[H^{\star,\star}(\mathcal X)\to
H^{\star+p,\star+q}(\mathcal X)\] induite par $f$ serait nulle. Pour d'autres
spectres que $\HZl$, il peut exister de tels $f$ stablement fantômes et non
nuls (cf. \cite[Corollary~6.2.3.7]{riou-rr}), c'est-à-dire que $f\neq 0$
mais qu'on ne peut pas détecter la non-nullité de $f$ en l'\guil{évaluant}
sur des spectres de la forme $\mathcal X\wedge S^{i,j}$ avec $\mathcal X\in
\Hopt$.

Nous allons montrer que sous les hypothèses du
théorème~\ref{theoreme-voevodsky}, ce phénomène ne se produit par pour
l'anneau bigradué des endomorphismes de $\HZl$ :

\begin{theoreme1}\label{theoreme-pas-d-operations-fantomes}
Soit $\ell$ un nombre premier. Soit $k$ un corps de caractéristique zéro.
Soit $(p,q)\in\mathbf{Z}^2$.
\begin{itemize}
\item Il n'existe pas de morphisme stablement fantôme non nul
$\HZl\to\HZl\wedge S^{p,q}$ dans $\SH[k]$.
\item On dispose d'un isomorphisme canonique :
\[A^{p,q}\simeq \Hom_{\SH}(\HZl,\HZl\wedge S^{p,q})\;\text{.}\]
\end{itemize}
\end{theoreme1}

Le deuxième résulte énoncé résulte immédiatement du premier et du
théorème~\ref{theoreme-voevodsky}. Nous allons donc démontrer uniquement le
premier énoncé, et ce en utilisant de manière essentielle certaines étapes
de la démonstration du théorème~\ref{theoreme-voevodsky}.

Si le groupe $\Hom_{\SHnaive}(\HZl,\HZl\wedge S^{p,q})$ s'identifie à
\[\lim_n \tH^{2n+p,n+q}(K_n)\;\text{,}\]
on sait d'après \cite[Lemme~6.5]{riou-sh} que le groupe des morphismes
stablement fantômes $\HZl\to\HZl\wedge S^{p,q}$ s'identifie à
\[\R^1 \lim_n \tH^{2n+p-1,n+q}(K_n)\;\text{,}\]
où $\R^1 \lim_n$ est le premier foncteur dérivé à droite du foncteur
limite projective au niveau des systèmes projectifs de groupes abéliens
indexés par $\mathbf{N}$.

Les résultats importants de \cite{voevodsky-eilenberg-maclane} sur les
produits symétriques sont utilisés au début de la démonstration de
\cite[Theorem~3.49]{voevodsky-eilenberg-maclane} (i.e. le
théorème~\ref{theoreme-voevodsky}) pour obtenir l'existence d'objets
$M'_n\in\DMNeg[k,\Fl]$ tels que $\tM(K_n)\simeq M'_n(n)[2n]$.

Les morphismes $\tM(K_n)(1)[2]\to \tM(K_{n+1})$ induits par les morphismes
d'assemblage sur $\HZl$ correspondent donc à des morphismes $M'_n\to
M'_{n+1}$ dans $\DMNeg[k,\Fl]$ (on utilise ici le \guil{théorème de
simplification}).
Pour tout $(p,q)\in\mathbf{Z}^2$, le système projectif $(H^{2n+p,n+q}(K_n))_n$
s'identifie donc au système $(H^{p,q}(M'_n))_n$ où on rappelle que 
$H^{p,q}(M'_n)=\Hom_{\DMmoins[k,\Fl]}(M'_n,\Fl(q)[p])$.

Comme il est établi dans la démonstration de
\cite[Theorem~3.49]{voevodsky-eilenberg-maclane} en utilisant le résultat
de scindage \cite[Corollary~2.71]{voevodsky-eilenberg-maclane}, le triangle
distingué dans $\DMmoins[k,\Fl]$ définissant la colimite homotopique du
système $(M'_n)_n$ est une suite exacte courte scindée :
\[0\to \oplus_n M'_n\vers{i} \oplus_n M'_n \to \hocolim_n M'_n\to 0\;\text{.}\]
Ceci repose sur le fait crucial qu'outre qu'ils existent, les objets
$M'_n$ appartiennent à la sous-catégorie $\overline{SPT}$ de $\DMmoins[k,\Fl]$
formée des sommes diectes d'objets de la forme $\Fl(i)[2i+j]$ avec $i,j\geq
0$.

Pour tout foncteur contravariant $F\colon \DMmoins[k,\Fl]\to \Ab$
transformant sommes directes en produits, on a une suite exacte :
\[0\to \lim_n F(M'_n)\to \prod_n F(M'_n)\vers{F(i)} \prod_n F(M'_n)
\to \R^1\lim_n F(M'_n)\to 0\;\text{.}\]
Le fait que $i$ admette une rétraction implique que $F(i)$ admet une
section, ce qui montre l'annulation de $\R^1\lim_n F(M'_n)$. En appliquant
ceci aux foncteurs $F=H^{p-1,q}$, on obtient l'annulation des groupes 
$\R^1\lim_n H^{p-1,q}(M'_n)$, ce qui achève la démonstration du
théorème~\ref{theoreme-pas-d-operations-fantomes}.


\appendix
\section{Le classifiant $\Bgm G$}

\begin{proposition1}\label{proposition-bgm}
Soit $S$ un schéma noethérien. Soit $G$ un schéma en groupes affine
lisse de type
fini sur $S$. Soit $V$ un fibré vectoriel sur $S$ et $\rho\colon G\to
\GL(V)$ une représentation linéaire fidèle de $G$ (c'est-à-dire que $\rho$
est un monomorphisme). On suppose que pour tout $i\geq 1$, on s'est donné un
ouvert $G$-invariant $U_i$ de l'espace affine
$V^{\oplus i}$ sur lequel $G$ agisse librement et que ceux-ci vérifient que
pour tout $i\geq 1$, $U_i\times\{0\}\subset U_{i+1}$ et $U_i\times
V^{\oplus i}\cup V^{\oplus i}\times U_i\subset U_{2i}$.
On note $U_\infty=\colim_i U_i\in \Sm^\opp\Ens$ et $(G\backslash
U_\infty)_\et$ le faisceau étale quotient de l'action libre de $G$ sur
$U_\infty$.

Si on fait l'hypothèse supplémentaire que pour $i$ assez grand il existe
$u\in U_i(S)$ tel que $\GL(V).u\subset U_i$ alors le morphisme canonique
$(G\backslash U_\infty)_\et\to \Bet G$ est une $\mathbf{A}^1$-équivalence
faible et $U_\infty$ est $\mathbf{A}^1$-contractile.
\end{proposition1}

Il manque un petit argument à \cite{morel-voevodsky} pour justifier ceci,
qui est énoncé sous une forme voisine dans
\cite[p.~133]{morel-voevodsky}. Il s'agit de montrer
que les ouverts $U_i$ des $V^{\oplus i}$ constituent un \guil{gadget
admissible muni d'une bonne action de $G$}. Par le résultat important et
technique \cite[Proposition~2.6, p.~135]{morel-voevodsky},
cela impliquera que
$(G\backslash U_\infty)\to \Bet G$ est une $\mathbf{A}^1$-équivalence
faible.

Il faut montrer que pour tout $G$-torseur étale $T$ sur $X$, le morphisme
évident entre faisceaux étales quotients $(G\backslash (T\times U_i))_\et\to
(G\backslash T)_\et =X$ est un épimorphisme de faisceaux pour la topologie
de Nisnevich. Pour cela, il suffit (et il faut) montrer que localement
Nisnevich sur $X$, il existe un morphisme $G$-équivariant $T\to U_i$.
Le $G$-torseur $T$ sur $X$ se plonge de façon évidente dans le
$\GL(V)$-torseur $T'=((\GL(V)\times_S T)/G)_\et$ (où l'on fait agit $G$ sur
$\GL(V)\times_S T$ par la formule $g.(m,t)=(m\rho(g^{-1}),gt)$). Localement
Nisnevich (et même Zariski) sur $X$, ce torseur $T'$ est trivial. On peut
donc supposer qu'il existe un isomorphisme $T'\simeq \GL(V)\times X$ de
$\GL(V)$-torseurs sur $X$. En composant la première projection puis
le morphisme $\GL(V)\to U_i$ donné par l'action sur l'élément $u$ fourni
par l'hypothèse de la proposition, on obtient le morphisme $G$-équivariant
$T\to U_i$ voulu.

\begin{remarque1}\label{remarque-proposition-bgm-1}
Si des ouverts $U_i\subset V^{\oplus i}$ vérifient les hypothèses de la
proposition~\ref{proposition-bgm} pour $G$, alors ils les
vérifient également pour l'action induite sur un sous-groupe $H\subset G$.
Comme elles sont faciles à satisfaire pour $\GL(V)$ en définissant $U_i$ pour
tout $i\geq 1$ comme étant le plus grand ouvert de $V^{\oplus i}$ sur
lequel $\GL(V)$ agisse librement ($U_i$ paramétrise les épimorphismes
$\mathbf{A}^i\to V$ de fibrés vectoriels), on obtient que la proposition
s'applique pour tout schéma en groupes affine lisse $G$ pour lequel on ait un
monomorphisme $G\to \GL(V)$.
On peut aussi choisir pour $U_i$ le plus grand ouvert
de $V^{\oplus i}$ sur lequel $G$ agisse librement.
Cet ouvert est plus grand que celui induit par
la construction précédente pour $\GL(V)$.
\end{remarque1}

\begin{remarque1}
Dans le cas favorable, les faisceaux étales quotients $(G\backslash
U_i)_\et$ sont représentables par un objet de $\Sm$ qu'il convient de
noter $G\backslash U_i$. C'est évidemment le cas si $G$ est un groupe fini
(puisque ce quotient $G\backslash U_i$ est un ouvert du quotient
$G\backslash V^{\oplus i}$ de SGA~1~V~1).
C'est également le cas si $G=\GL(V)$ et que les ouverts
$U_i$ sont ceux décrits dans la remarque~\ref{remarque-proposition-bgm-1} :
dans ce cas $G\backslash U_i$ est la grassmannienne des
$r$-plans dans $\mathbf{A}^i$ où $r$ est le rang de $V$.

Dans le cas où $S=\Spec k$ et $G$ est un sous-$k$-schéma en groupes lisse
fermé de $\GL_n$, on peut montrer que si on prend pour $U_i$
l'ouvert de $V^{\oplus i}$ où $\GL_n$ agit librement, alors le quotient
$(G\backslash U_i)_\et$ est représentable (ceci se ramène facilement à la
démonstration du fait que le quotient $G\backslash \GL_n$ soit
représentable et que $\GL_n\to G\backslash \GL_n$ soit un $G$-torseur
étale, voir SGA~3~VI${}_{\text{A}}$~3.2). (Pour une présentation
alternative, voir
\cite[\S{}1]{totaro} ou \cite[Lemma~9.2]{colliot-thelene-sansuc}.)
\end{remarque1}

\begin{proposition1}
Soit $k$ un corps parfait.
Soit $i\colon H\subset G$
une inclusion entre groupes finis. Ce morphisme induit un
morphisme $\Bet i\colon \Bet H\to \Bet G$ dans $\Ho[k]$ puis un morphisme
$M(\Bet i)\colon M(\Bet H)\to M(\Bet G)$ dans $\DMNeg[k,\Lambda]$.
On dispose d'un morphisme de transfert $\Tr\colon M(\Bet G)\to M(\Bet H)$
tel que le morphisme composé dans $\DMNeg[k,\Lambda]$
\[\xymatrix{M(\Bet G)\ar[r]^-{\Tr}&
M(\Bet H)\ar[r]^-{M(\Bet i)}& M(\Bet G)}\]
soit la multiplication par l'indice $[G:H]$.
\end{proposition1}

On choisit une représentation fidèle $\rho\colon G\to \GL(V)$ de $G$. Pour
tout $i\geq 1$, on note $U_i$ l'ouvert de $V^{\oplus i}$ où $G$ (et donc
$H$) agit librement. Les faisceaux étales quotients $(G\backslash U_i)_\et$
et $(H\backslash U_i)_\et$ sont représentables par des $k$-schémas lisses
que l'on note $G\backslash U_i$ et $H\backslash U_i$. Les colimites
$G\backslash U_\infty$ et $H\backslash U_\infty$ de ces systèmes (calculées
dans la catégorie des préfaisceaux sur $\Sm[k]$, la notation est ici
quelque peu abusive) s'identifient
respectivement à $\Bet G$ et $\Bet H$ dans $\Ho[k]$ d'après la
proposition~\ref{proposition-bgm}. On dispose du diagramme suivant :
\[\xymatrix{\ar[d]_{H\text{-torseur}}
U_\infty \ar@{=}[r]^{H\text{-équiv.}} &U_\infty \ar[d]^{G\text{-torseur}} \\
H\backslash U_\infty \ar[r] & G\backslash U_\infty}\]
Il permet d'observer, ce qui va de soi, que le morphisme évident
$H\backslash U_\infty \to G\backslash U_\infty$ est une façon de
représenter le morphisme $\Bet i\colon \Bet G\to \Bet H$ \emph{via} les
identifications précédemment faites.

Ainsi, le morphisme $\Tr$ cherché peut être représenté par un morphisme de
faisceaux avec tranferts obtenu en passant à la colimite le morphisme
canonique $\Ltr(G\backslash U_i)\to \Ltr(H\backslash U_i)$ correspondant à la
correspondance finie de $G\backslash U_i$ dans $H\backslash U_i$ associée
au revêtement étale $H\backslash U_i\to G\backslash U_i$. Le morphisme
composé $\Ltr(G\backslash U_\infty)\vers{\Tr} \Ltr(H\backslash U_\infty)\to
\Ltr(G\backslash U_\infty)$ est bien sûr la multiplication par le degré
$[G:H]$ de ce revêtement.

\begin{corollaire1}\label{corollaire-transfert-motif-classifiants}
Soit $k$ un corps parfait. Soit $i\colon H\subset G$ une inclusion entre
groupes finis. On suppose que l'indice $[G:H]$ est inversible dans l'anneau
de coefficients $\Lambda$.
Alors, le morphisme $M(\Bet H)\to M(\Bet G)$ induit par l'inclusion est un
épimorphisme scindé dans $\DMNeg[k,\Lambda]$.

Plus précisément, notons $N$ le normalisateur de $H$ dans $G$. Alors, les
co-invariants de $N$ agissant sur $M(\Bet H)$
sont représentables par un épimorphisme scindé
$M(\Bet H)\to M(\Bet H)_N$ dans $\DMNeg[k,\Lambda]$.
Le morphisme évident $M(\Bet H)\to M(\Bet G)$ induit un morphisme
$M(\Bet H)_N\to M(\Bet G)$ qui est un épimorphisme scindé.
\end{corollaire1}

La seule chose à ajouter pour ce corollaire est que l'action de $N$ par
conjugaison sur $H$ puis sur $M(\Bet H)$ se factorise en une action de
$N/H$ puisque les automorphismes intérieurs de $H$ induisent l'identité sur
$\Bet H\in\Ho[k]$. La prise des co-invariants de $M(\Bet H)$ sous $N$  peut
donc se faire par le procédé standard de moyenne pour l'action du groupe
$N/H$ dont l'ordre est bien inversible dans $\Lambda$.


\bibliography{article}
\bibliographystyle{amsplain}

\end{document}